\documentclass[reqno]{amsart}
\usepackage{amssymb}
\usepackage{hyperref}
\usepackage{graphicx}
\usepackage{subcaption}
\usepackage{color}
\usepackage{textcomp}

\theoremstyle{definition}

\newtheorem{test}{Test}[section]
\theoremstyle{remark}

\begin{document}
\title[Global Convergence and Uniqueness For An Inverse Problem]{Global
Convergence for an $n-$Dimensional Coefficient Inverse Problem for a
Wave-like Equation with a Single Location of the Point Source}
\author{Michael V. Klibanov$^\dagger$}
\author{Jingzhi Li$^\ddagger$}
\author{Tian Niu$^\S $}
\author{Vladimir G. Romanov$^\mathparagraph$}
\thanks{\textbf{Funding:} The work of the second author was partially
supported by the Shenzhen Sci-Tech Fund RCJC20200714114556020, Guangdong
Basic and Applied Research Fund 2023B1515250005. The work of the third
author was partially supported by Guangdong Basic and Applied Research Fund
2023B1515250005. The work of the fourth author was supported by the
Mathematical Center in Akademgorodok under the agreement No. 075-15-2025-348
with the Ministry of Science and Higher Education of the Russian Federation.}
\thanks{$^\dagger$Corresponding author. Department of Mathematics and
Statistics, University of North Carolina at Charlotte, Charlotte, NC 28223,
USA (mklibanv@charlotte.edu).}
\thanks{$^\ddagger$Department of Mathematics, Southern University of Science
and Technology, Shenzhen 518055, China (li.jz@sustech.edu.cn).}
\thanks{$^\S $Department of Mathematics, Southern University of Science and
Technology, Shenzhen 518055, China (niut@sustech.edu.cn).}
\thanks{$^\mathparagraph$Sobolev Institute of Mathematics, Novosibirsk
630090, Russian Federation (romanov@math.nsc.ru).}
\date{}
\dedicatory{}

\begin{abstract}
In this paper, the first globally convergent numerical method is developed
for an $n-$Dimensional coefficient inverse problem for a wave-like equation
with the unknown potential in the most challenging case when the Dirac $%
\delta -$ function is present in one of initial conditions for that equation
and only a single location of the point source is used. Based on a Carleman
estimate, convergence analysis is carried out. This analysis ensures the
global convergence of the proposed numerical method, i.e. its convergence to
the true solution is rigorously guaranteed regardless on the availabilty of
a good first guess about this solution. A certain uniqueness theorem is
proven as well. Extensive computational experiments with noisy data for both
2-Dimensional and 3-Dimensional cases demonstrate a high reconstruction
accuracy.
\end{abstract}

\maketitle

\section{Introduction}

\label{sec:1}

A Coefficient Inverse Problem (CIP) for a Partial Differential Equation
(PDE) is the problem of the determination of an unknown coefficient(s) of
this PDE from boundary data. Those boundary data are also called
\textquotedblleft input data". In applications, the input data are results
of some boundary measurements. Let $m$ be the number of free variables in
the input data for a CIP and $n$ be the number of free variables in the
unknown coefficient of that PDE. That CIP is called \textquotedblleft
formally determined" if 
\begin{equation}
m=n.  \label{1.0}
\end{equation}

This is the first publication, in which a numerical method with a rigorously
guaranteed global convergence property is developed for an $n-$Dimensional
CIP with condition (\ref{1.0}) for a wave-like Partial Differential Equation
(PDE) with an unknown potential. The most challenging case is considered
here: the case when the Dirac $\delta -$function is present in one of
initial conditions of that PDE with only a single location of the point
source. Thus, we work here with a single measurement event. The single
measurement case is the most economical one from the measurement applied
standpoint. On the other hand, it is the most difficult one from the
analytical and numerical standpoints. The global convergence property of the
so-called \textquotedblleft convexification" numerical method of this paper
is ensured by our detailed convergence analysis. The term \textquotedblleft
global convergence" is specified in Definition 1.1. Thus, the
convexification method is not a version of the conventional small
perturbation approach. Globally convergent numerical methods were not
developed for this CIP in the past. Results of extensive numerical testing
for complicated media structures for both 2-Dimensional and 3-Dimensional
cases confirm our theory. It is shown below that the CIP considered here has
an application in acoustics.

\textbf{Definition 1.1.} \emph{We call a numerical method for a Coefficient
Inverse Problem (CIP) globally convergent if there is a theorem, which
claims that this method provides points in a sufficiently small neighborhood
of the true solution of that CIP without an advanced knowledge of any point
in that neighborhood. In other words, convergence of this numerical method
to the true solution of that CIP is rigorously guaranteed without an
availability of a good first guess about that solution.}

In fact, we work below with two CIPs: the main one is for that wave-like
equation and an auxiliary one is for a parabolic equation. To arrange a
focused structure of this paper, we first formulate in section 2 our two
CIPs. In section 3, we describe some previously published relevant results,
provide a brief description of the main new elements of this paper, as well
as outline the key ideas of the above mentioned convexification numerical
method.

Thus, the rest of this paper is structured as follows. In section 2 we
formulate forward and inverse problems for both above mentioned PDEs: the
wave-like equation and the auxiliary parabolic equation. The parabolic
equation is obtained from the wave-like equation via an analog of the
Laplace transform, see \cite[formula (7.130)]{LRS} for that analog. In the
originating CIP for that wave-like equation, we assign the input boundary
data on the infinite time interval $t\in \left( 0,\infty \right) .$ Hence,
the input boundary data for the CIP for that parabolic equation are obtained
from those for the wave-like equation via that analog of the Laplace
transform. After this, we work only with the CIP for the parabolic equation
on a finite time interval. We formulate in section 2 two theorems about some
properties of fundamental solutions of those two PDEs. These two theorems
are proven in Appendices A and B. In section 3, we outline results of some
relevant previous publications. In addition, we briefly describe in section
3 the main new elements of this paper as well as outline the main idea of
the convexification method. In section 4 we describe our transformation
procedure, which transforms the second CIP of section 2 for that parabolic
equation in a boundary value problem for a system of coupled nonlinear
elliptic equations. In section 5, we develop a new version of the
convexification method. In section 6, we carry out a detailed convergence
analysis for this method. In section 7, we prove a certain uniqueness
theorem for CIP2. In section 8 we describe our numerical studies. A summary
of our results is in section 9. All functions considered below are real
valued ones.

\section{Forward and Inverse Problems}

\label{sec:2}

Denote $\mathbf{x=}\left( x_{1},x_{2},x_{3},...,x_{n}\right) \in \mathbb{R}%
^{n}$ points in $\mathbb{R}^{n},$ where $n\geq 1.$ Let $A,B,D>0$ be three
numbers, where $A<B$. We define the domain $\Omega \subset \mathbb{R}^{n}$
as the rectangular prism%
\begin{equation}
\Omega =\left\{ \mathbf{x:}\text{ }x_{1}\in \left( A,B\right) ,\text{ }%
\left\vert x_{j}\right\vert <D,\text{ }j=2,...,n\right\} .  \label{2.1}
\end{equation}
The boundary $\partial \Omega $ of the domain $\Omega $ consists of two
parts,%
\begin{equation}
\left. 
\begin{array}{c}
\partial \Omega =\Gamma _{0}\cup \Gamma _{1}, \\ 
\Gamma _{0}=\left\{ x_{1}=B,\left\vert x_{i}\right\vert <D,i=2,...,n\right\}
, \\ 
\Gamma _{1}=\partial \Omega \diagdown \Gamma _{0}.%
\end{array}%
\right. .  \label{2.2}
\end{equation}%
Let $T>0$ be an arbitrary number. Denote%
\begin{equation}
\left. 
\begin{array}{c}
D_{T}^{n+1}=\mathbb{R}^{n}\times \left( 0,T\right) ,\text{ }D_{\infty
}^{n+1}=\mathbb{R}^{n}\times \left( 0,\infty \right) , \\ 
S_{T}=\partial \Omega \times \left( 0,T\right) ,\text{ }\Gamma _{0,T}=\Gamma
_{0}\times \left( 0,T\right) .\text{ }%
\end{array}%
\right.  \label{2.3}
\end{equation}%
In addition, denote%
\begin{equation}
S_{\infty }=\partial \Omega \times \left( 0,\infty \right) ,\text{ }\Gamma
_{0,\infty }=\Gamma _{0}\times \left( 0,\infty \right) .  \label{2.03}
\end{equation}

\textbf{Remark 2.1.}\emph{\ }In principle, a more general domain $\Omega $
is possible to work with. However, we consider a simple domain $\Omega $\ in
(\ref{2.1}) in order to simplify the presentation, since CIPs of this paper
are quite complicated ones in their own rights. It is also important that
this simple structure of $\Omega $\ helps us to simplify our computational
task of section 8.

Let the function $a\left( \mathbf{x}\right) $ satisfies the following
conditions:%
\begin{equation}
\left. 
\begin{array}{c}
a\in C^{\ell }\left( \mathbb{R}^{n}\right) ,~\Vert a\Vert _{C^{\ell }\left( 
\mathbb{R}^{n}\right) }\leq a_{0}, \\ 
\ell =5\left[ n/2\right] +3,%
\end{array}%
\right.  \label{2.7}
\end{equation}%
\begin{equation}
a\left( \mathbf{x}\right) =0\text{ for }\mathbf{x}\in \mathbb{R}%
^{n}\diagdown \Omega ,  \label{2.8}
\end{equation}%
where $\lfloor n/2\rfloor $ denotes the integer part of $n/2$ and $a_{0}$ is
a positive number.

\textbf{Remark 2.2.} We explore smoothness requirement (\ref{2.7}) in
Theorem 2.1 below in order to use the technique of \cite{Rom1}, which
requires this smoothness. It is well known that extra smoothness
requirements are of a secondary concern in the theory of CIPs, see, e.g. 
\cite{Nov,Rom1}.

\subsection{Statements of forward and inverse problems}

\label{sec:2.1}

\subsubsection{The hyperbolic case}

\label{sec:2.1.1}

Consider the following Cauchy problem%
\begin{equation}
U_{tt}=\Delta U+a\left( \mathbf{x}\right) U,\text{ }\left( \mathbf{x}%
,t\right) \in D_{\infty }^{n+1},  \label{2.9}
\end{equation}%
\begin{equation}
U\left( \mathbf{x},0\right) =0,\text{ }U_{t}\left( \mathbf{x},0\right)
=\delta \left( \mathbf{x}\right) .  \label{2.10}
\end{equation}

Problem (\ref{2.9}), (\ref{2.10}) is the forward problem for hyperbolic
equation (\ref{2.9}). Theorem 2.1 ensures existence and uniqueness of the
solution of this problem. The proof of this theorem can be found in Appendix
A. Denote%
\begin{equation*}
\theta _{0}(t)=\left\{ 
\begin{array}{c}
1,t>0, \\ 
0,t<0,%
\end{array}%
\right.
\end{equation*}%
\begin{equation*}
\theta _{s}(t)=\frac{t^{s}}{s!}\theta _{0}(t),~\theta _{-s}(t)=\frac{d^{s}}{%
dt^{s}}\theta _{0}(t),~s=1,2,\ldots .
\end{equation*}%
Hence, $\theta _{0}(t)$ is the Heaviside function. For each $t>0$ consider
the domain%
\begin{equation}
\Omega _{t}=\left\{ \mathbf{x}\in \mathbb{R}^{n}:\text{ }|\mathbf{x}|\leq
t\right\} .  \label{10}
\end{equation}

\textbf{Theorem 2.1} (existence and uniqueness results for forward problem (%
\ref{2.9}), (\ref{2.10})). \emph{Assume that conditions} \emph{(\ref{2.7})
and (\ref{2.8}) \ hold}. \emph{Then there exists the unique solution} $U(%
\mathbf{x},t)$ \emph{\ of problem} \emph{(\ref{2.9}), (\ref{2.10}),} \emph{%
which has the following properties}

1)\emph{\ If} $n=2m+1$, $m\geq 1$, \emph{then the following representation
is valid with} $S=m+1$ 
\begin{equation}
U(\mathbf{x},t)=\sum\limits_{s=-m}^{S}\alpha _{s}(\mathbf{x})\theta
_{s}(t^{2}-|\mathbf{x}|^{2})+{U}_{S}(\mathbf{x},t),~t>0,  \label{11}
\end{equation}%
\emph{where the} \emph{functions} $\alpha _{s}(\mathbf{x})$ \emph{are:}%
\begin{equation}
\left. 
\begin{array}{c}
\alpha _{-m}(\mathbf{x})=\left( 2\pi ^{m}\right) ^{-1}, \\ 
\alpha _{-m+1}(\mathbf{x})=\left( 8\pi ^{m}\right)
^{-1}\int\limits_{0}^{1}a(y\mathbf{x})\,dy, \\ 
\alpha _{s}(\mathbf{x})=\left( 8\pi ^{m}\right)
^{-1}\int\limits_{0}^{1}y^{s+m-1}(\Delta \alpha _{s-1}(\xi )+a(\xi )\alpha
_{s-1}(\xi ))|_{\xi =y\mathbf{x}}\,dy,\>\text{ } \\ 
s\in \left[ -m+2,S\right] .%
\end{array}%
\right.  \label{12}
\end{equation}%
2) \emph{If} $n=2m,m\geq 1$, \emph{then } 
\begin{equation}
U(\mathbf{x},t)=\sum\limits_{s=-m}^{S}\alpha _{s}(\mathbf{x})\theta
_{s+1/2}(t^{2}-|\mathbf{x}|^{2})+{U}_{S}(\mathbf{x},t),~t>0,  \label{13}
\end{equation}%
\emph{with} $S=m+1$, \emph{where functions} $\theta _{s+1/2}(t)$\emph{\ are
given by the equalities} 
\begin{eqnarray*}
&&\theta _{-1/2}(t)=\frac{1}{\sqrt{t}}\theta _{0}(t),\quad \theta
_{s+1/2}(t)=\frac{2^{s+1}t^{s+1/2}}{(2s+1)!!}\theta _{0}(t),\quad
s=0,1,2,\ldots , \\
&&\theta _{s-1/2}(t)=\frac{d^{s}}{dt^{s}}\theta _{-1/2}(t),~s=1,2,\ldots .
\end{eqnarray*}%
\emph{Here} $(2s+1)!!=1\cdot 3\cdot 5\cdot \ldots \cdot (2s+1)$.

3)\emph{\ The remainder} \emph{term} ${U}_{S}(\mathbf{x},t)=0$ \emph{\ for} $%
t<|\mathbf{x}|$ \emph{in both cases (\ref{11}) and (\ref{13}), and the
function }${U}_{S}(\mathbf{x},t)$ \emph{is twice differentiable with respect
to both }$\mathbf{x}$\emph{\ and }$t$\emph{\ for }$t\geq \left\vert \mathbf{x%
}\right\vert .$ \emph{For any fixed }$t>0$\emph{\ this function is bounded
for }$\mathbf{x}\in \Omega _{t}$\emph{\ together with its derivatives with
respect to }$\mathbf{x}$ \emph{and} $t$ \emph{up to the second order. If} $%
t\rightarrow \infty $\emph{, then the function} ${U}_{S}(\mathbf{x},t)$ 
\emph{together with its derivatives with respect to} $\mathbf{x}$\emph{\ and 
}$t$ \emph{up to the second order, grows not faster than} $e^{Yt}$\emph{\
with a number} $Y>0,$ \emph{which depends only on the number }$a_{0}$ \emph{%
in (\ref{2.7}).}

We now formulate the main CIP of this paper.

\textbf{First Coefficient Inverse Problem (CIP1).} Assume that the
coefficient $a\left( \mathbf{x}\right) $ in (\ref{2.9}) is unknown and
satisfies conditions (\ref{2.7}), (\ref{2.8}). Using conditions (\ref{2.9})
and (\ref{2.10}), find the coefficient $a\left( \mathbf{x}\right) $,
assuming that the following two functions $f_{0}\left( \mathbf{x},t\right) $
and $f_{1}\left( \mathbf{x},t\right) $ are given:%
\begin{equation}
U\left( \mathbf{x},t\right) \mid _{S_{\infty }}=f_{0}\left( \mathbf{x}%
,t\right) ,\text{ }\partial _{x_{1}}U\left( \mathbf{x},t\right) \mid
_{\Gamma _{0,\infty }}=f_{1}\left( \mathbf{x},t\right) ,  \label{2.11}
\end{equation}%
where $S_{\infty }$ and $\Gamma _{0,\infty }$ are defined in (\ref{2.03}).

Thus, the input boundary data (\ref{2.11}) for CIP1 are given on the
infinite time interval $t\in \left( 0,\infty \right) .$

\subsubsection{An application of CIP1 in acoustics}

\label{sec:2.1.2}

Following \cite{Ber}, we now describe an application of CIP1 in acoustics.
Let $n=3.$ Let $\rho \left( \mathbf{x}\right) $ be the acoustical density of
the medium with the speed of propagation of the acoustical waves $c\equiv 1$%
. We assume that 
\begin{equation}
\left. 
\begin{array}{c}
\rho \in C^{10}\left( \mathbb{R}^{3}\right) , \\ 
\rho \left( \mathbf{x}\right) \geq 1\text{ in }\mathbb{R}^{3}, \\ 
\rho \left( \mathbf{x}\right) =1\text{ for }\mathbf{x}\in \mathbb{R}%
^{n}\diagdown \Omega .%
\end{array}%
\right.  \label{2.110}
\end{equation}%
We impose regularity requirement in (\ref{2.110}) since we want the
resulting function $a\left( \mathbf{x}\right) ,$ which is introduced below
in this sub-subsection, to satisfy condition (\ref{2.7}), also, see Remark
2.2. Let $Z\left( \mathbf{x},t\right) $ be the amplitude of the sound wave
at the point $\mathbf{x}$ at the moment of time $t$. The acoustic equation
is:%
\begin{equation}
\frac{1}{\rho }Z_{tt}=\operatorname{div}\left( \frac{1}{\rho }\nabla Z\right) .
\label{2.111}
\end{equation}%
Setting in (\ref{2.111}) 
\begin{equation}
Z=\sqrt{\rho }W,  \label{2.0111}
\end{equation}%
we obtain equation (\ref{2.9}) for the function $W$ with initial conditions
, 
\begin{equation*}
\left. W_{tt}=\Delta W+a\left( \mathbf{x}\right) W,\text{ }\left( \mathbf{x}%
,t\right) \in D_{\infty }^{n+1},\right.
\end{equation*}%
where the coefficient $a\left( \mathbf{x}\right) $ is:%
\begin{equation}
a\left( \mathbf{x}\right) =-\sqrt{\rho \left( \mathbf{x}\right) }\Delta
\left( \frac{1}{\sqrt{\rho \left( \mathbf{x}\right) }}\right) .
\label{2.112}
\end{equation}%
We impose initial conditions (\ref{2.10}) on the function $W$, 
\begin{equation*}
W\left( \mathbf{x},0\right) =0,\text{ }W_{t}\left( \mathbf{x},0\right)
=\delta \left( \mathbf{x}\right) .
\end{equation*}%
Next, using (\ref{2.110}), we obtain that the function $a\left( \mathbf{x}%
\right) $ in (\ref{2.112}) satisfies conditions (\ref{2.7}), (\ref{2.8}). \
Now, it follows from (\ref{2.110}) and (\ref{2.0111}) that 
\begin{equation*}
Z\left( \mathbf{x},t\right) \mid _{S_{\infty }}=W\left( \mathbf{x},t\right)
\mid _{S_{\infty }},\text{ }\partial _{x_{1}}Z\left( \mathbf{x},t\right)
\mid _{\Gamma _{0,\infty }}=\partial _{x_{1}}W\left( \mathbf{x},t\right)
\mid _{\Gamma _{0,\infty }}.
\end{equation*}

Suppose that the function $a\left( \mathbf{x}\right) $ in (\ref{2.112}) is
found. Denote 
\begin{equation}
\mu \left( \mathbf{x}\right) =\frac{1}{\sqrt{\rho \left( \mathbf{x}\right) }}%
.  \label{2.113}
\end{equation}%
Then, using the third line of (\ref{2.110}), (\ref{2.112}) and (\ref{2.113}%
), we obtain 
\begin{equation}
\Delta \mu +a\left( \mathbf{x}\right) \mu =0,  \label{2.114}
\end{equation}%
\begin{equation}
\mu \mid _{\partial \Omega }=1,\text{ }\partial _{\nu }\mu \mid _{\partial
\Omega }=0,  \label{2.115}
\end{equation}%
where $\partial _{\nu }$ is the normal derivative at $\partial \Omega .$
Thus, we have both Dirichlet and Neumann boundary conditions (\ref{2.115})
for elliptic equation (\ref{2.114}). It is well known that there exists at
most one solution $\mu \in C^{2}\left( \overline{\Omega }\right) $ of
problem (\ref{2.114}), (\ref{2.115}), see, e.g. \cite[section 1 of chapter 4]%
{LRS} and \cite[section 2.6.2]{KL}. Furthermore, if this solution exists,
then it can be found numerically via a slight modification of the version of 
\cite{KT} of the so-called \textquotedblleft quasi-reversibility" method.
Convergence of this method for noisy boundary data is proven in \cite{KT}.
We do not discuss that modification here for brevity.

In summary, we have shown above that CIP1 has an application in the problem
of the determination of the acoustical density of the medium from boundary
measurements.

\subsubsection{The parabolic case}

\label{sec:2.1.3}

For any appropriate function $y\left( \mathbf{x,}t\right) $ consider the
following analog of the Laplace transform with respect to $t$ \cite[formula
(7.130)]{LRS}:%
\begin{equation}
\mathcal{L}\left( y\right) \left( \mathbf{x},t\right) =\frac{1}{2\sqrt{\pi
t^{3}}}\int\limits_{0}^{\infty }\exp \left( -\frac{\tau ^{2}}{4t}\right) y(%
\mathbf{x},\tau )\tau d\tau ,\quad t>0.  \label{2.12}
\end{equation}%
Denote%
\begin{equation}
u\left( \mathbf{x},t\right) =\mathcal{L}\left( U\right) \left( \mathbf{x}%
,t\right) ,\text{ }  \label{2.13}
\end{equation}%
\begin{equation}
g_{0}\left( \mathbf{x},t\right) =\mathcal{L}\left( f_{0}\right) \left( 
\mathbf{x},t\right) ,\text{ }g_{1}\left( \mathbf{x},t\right) =\mathcal{L}%
\left( f_{1}\right) \left( \mathbf{x},t\right) .  \label{2.14}
\end{equation}%
Then (\ref{2.9}), (\ref{2.10}), (\ref{2.12}) and (\ref{2.13}) lead to the
following Cauchy problem for any $T>0$ \cite[formulas (7.125), (7.126),
(7.130) and (7.131)]{LRS}:%
\begin{equation}
u_{t}=\Delta u+a\left( \mathbf{x}\right) u,\text{ }\left( \mathbf{x}%
,t\right) \in D_{T}^{n+1},  \label{2.15}
\end{equation}%
\begin{equation}
u\left( \mathbf{x},0\right) =\delta \left( \mathbf{x}\right) .  \label{2.16}
\end{equation}%
In addition, by (\ref{2.11}) and (\ref{2.14})%
\begin{equation}
u\left( \mathbf{x},t\right) \mid _{S_{T}}=g_{0}\left( \mathbf{x},t\right) ,%
\text{ }\partial _{x_{1}}u\left( \mathbf{x},t\right) \mid _{\Gamma
_{0,T}}=g_{1}\left( \mathbf{x},t\right) .  \label{2.17}
\end{equation}%
Therefore, conditions (\ref{2.15})-(\ref{2.17}) lead to the second CIP.

\textbf{Second Coefficient Inverse Problem (CIP2)}. Assume that the function 
$a\left( \mathbf{x}\right) $ is unknown and satisfies conditions (\ref{2.7}%
), (\ref{2.8}). Using conditions (\ref{2.15}), (\ref{2.16}), find this
function$,$ assuming that the functions $g_{0}\left( \mathbf{x},t\right) $
and $g_{1}\left( \mathbf{x},t\right) $ in (\ref{2.17}) are given.

Below we work only with CIP2. We have chosen in (\ref{2.15}) and (\ref{2.17}%
) an arbitrary number $T>0$ instead of $T=\infty $ in CIP1 since we work in
CIP2 only on a finite time interval.

\subsection{Asymptotic behavior of the function $u\left( \mathbf{x},t\right) 
$ at $t\rightarrow 0^{+}$}

\label{sec:2.2}

First, we need to establish an asymptotic behavior of the function $u\left( 
\mathbf{x},t\right) $ at $t\rightarrow 0^{+}.$ This behavior is established
on the basis of Theorem 2.2 and connection (\ref{2.12}), (\ref{2.13})
between solutions of forward problems (\ref{2.9}), (\ref{2.10}) and (\ref%
{2.15}), (\ref{2.16}).

Let $M>0$ and $\rho \in \left( 0,1\right) $ be two arbitrary numbers. Denote%
\begin{equation*}
B_{M}=\{\mathbf{x}\in \mathbb{R}^{n}:\,|\mathbf{x}|\leq M\},~G_{\rho }(M)=\{%
\mathbf{x}\in \mathbb{R}^{n}:\,|x_{i}|\geq \rho ,\text{ }i=1,...,n\}\cap
B_{M}.
\end{equation*}%
Let the function $u_{0}(\mathbf{x},t)$ be the fundamental solution of the
heat equation%
\begin{equation*}
\left. 
\begin{array}{c}
u_{0t}=\Delta u_{0}\text{ in }D_{T}^{n+1}, \\ 
u_{0}(\mathbf{x},0)=\delta \left( \mathbf{x}\right) .%
\end{array}%
\right.
\end{equation*}%
It is well known that 
\begin{equation}
u_{0}(\mathbf{x},t)=\frac{1}{(2\sqrt{\pi t})^{n}}\,\exp \left( -\frac{%
\left\vert \mathbf{x}\right\vert ^{2}}{4t}\right) .  \label{2.18}
\end{equation}

\textbf{Theorem 2.2.} \emph{Assume that conditions (\ref{2.7}) and (\ref{2.8}%
) hold. Then forward problem (\ref{2.15}), (\ref{2.16}) has a unique
solution }%
\begin{equation}
u\in C^{4,2}\left( G_{\rho }\left( M\right) \times \left[ 0,T\right] \right)
\label{2.19}
\end{equation}%
\emph{\ for any }$T>0.$ \emph{Moreover, the following asymptotic formulas
hold}%
\begin{equation}
u(\mathbf{x},t)=u_{0}(\mathbf{x},t)\left[ 1+O(t^{\kappa })\right] ,\text{ }%
t\rightarrow 0^{+},\text{ }\mathbf{x}\in B_{M},  \label{2.20}
\end{equation}%
\begin{equation}
\partial _{t}u(\mathbf{x},t)=\partial _{t}u_{0}(\mathbf{x},t)\left[
1+O(t^{\kappa })\right] ,\text{ }t\rightarrow 0^{+},\text{ }\mathbf{x}\in
B_{M},  \label{2.21}
\end{equation}%
\begin{equation}
u_{x_{i}}(\mathbf{x},t)=\partial _{x_{i}}u_{0}(\mathbf{x},t)\left[
1+O(t^{\kappa })\right] ,\text{ }t\rightarrow 0^{+},\text{ }\mathbf{x}\in
G_{\rho }(M),\   \label{2.22}
\end{equation}%
\begin{equation}
u_{x_{i}x_{j}}(\mathbf{x},t)=\partial _{x_{i}x_{j}}^{2}u_{0}(\mathbf{x},t)%
\left[ 1+O(t^{\kappa })\right] ,\text{ }t\rightarrow 0^{+},\text{ }\mathbf{x}%
\in G_{\rho }(M),  \label{2.23}
\end{equation}%
where $\ i,j=1,\dots ,n.$ \emph{Here }$\kappa =1/2$\emph{\ if }$n=2$\emph{\
and }$\kappa =1$\emph{\ otherwise.}

Formula (\ref{2.20}) was first obtained in \cite{KLF}. However, formulas (%
\ref{2.21})-(\ref{2.23}) were unknown in the past. The proof of Theorem 2.2
can be found in Appendix B. As to the $C^{4,2}-$smoothness of the function $%
u\left( \mathbf{x},t\right) $ in (\ref{2.19}), in fact, this function is
more smooth due to the smoothness requirement in (\ref{2.7}) \cite[chapter 4]%
{Lad}, also, see Remark 2.2. However, we need here only the $C^{4,2}-$%
smoothness of $u$ since this is sufficient for our derivations below.

\section{Some Previous Relevant Publications, The Main New Elements of This
Paper and a Brief Outline of the Convexification Method}

\label{sec:3}

\subsection{Some relevant previous publications}

\label{sec:3.1}

First, we rephrase a statement of section 1: the key challenge of both CIP1
and CIP2 is the presence of the $\delta \left( \mathbf{x}\right) -$function
in initial conditions (\ref{2.10}) and (\ref{2.16}) respectively with the
single location of the point source at $\left\{ 0\right\} $. Numerical
methods were not developed in the past neither for CIP1 nor for CIP2.

Berezansky \cite{Ber} has studied the question of uniqueness for a
3-Dimensional CIP for the Schr\"{o}dinger equation with the unknown
potential. This equation can be obtained from equation (\ref{2.9}) via the
Fourier transform with respect to $t$, see, e.g. the book of Vainberg \cite%
{V} for some sufficient conditions, which guarantee the applicability of
this transform. Condition (\ref{1.0}) does not hold in the uniqueness
theorem of \cite{Ber} since one of assumptions of that theorem is: 
\begin{equation}
m=5>n=3.  \label{3.01}
\end{equation}%
Another assumption of that theorem of \cite{Ber} is that the unknown
potential is a piecewise analytic function. However, it is desirable to
prove uniqueness theorems under a less restrictive assumption that the
unknown coefficient belongs to one of conventional function spaces, such as,
e.g. $C^{r}$.

As to the CIPs with condition (\ref{1.0}) for $n\geq 2$, we refer below only
to those uniqueness and stability results, which are proven under the above
mentioned requirement that an unknown coefficient of a PDEs belongs to one
of conventional function spaces. All such results are obtained by the
method, which was originated in the paper \cite{BukhKlib} with many follow
up publications. The method of \cite{BukhKlib} uses Carleman estimates.
Since the current paper is not a survey, we refer here only to a few
publications \cite{Boul,Isakov,Klib92,Ksurvey,KL}, which use the idea of 
\cite{BukhKlib}. In the case of equation (\ref{2.9}), the method of \cite%
{BukhKlib} requires that one of initial conditions for the function $U\left( 
\mathbf{x},t\right) $ should not vanish in the entire domain of interest,
which is unlike (\ref{2.10}). Still, combining a new idea with the method of 
\cite{BukhKlib}, Rakesh and Salo \cite{RS1} have proven Lipschitz stability
estimate and uniqueness result for a CIP for equation (\ref{2.9}) in the
3-Dimensional case assuming two measurement events generated by two incident
plane waves propagating in two opposite directions. We also refer to, e.g. 
\cite{Ma,RS2} for some extensions of this result.

\subsection{The main new elements of this paper}

\label{sec:3.2}

The key new element of this paper is that a globally convergent numerical
method for CIP2 is developed here for the first time. The same is true for
CIP1 since CIP1 is reduced to CIP2, as established above. An important new
element of our approach is that, using the analog (\ref{2.12}) of the
Laplace transform with respect to $t$, we reduce the original CIP1 for
hyperbolic PDE (\ref{2.9}) to the CIP2 for parabolic PDE (\ref{2.15}) and
work only with CIP2 since then. Computational results of this paper are new.

Another new element of this publication is that condition (\ref{1.0}) holds
in our case, unlike condition (\ref{3.01}) in \cite{Ber}. An additional new
feature here is that, instead of the piecewise analyticity requirement
imposed on the unknown coefficient imposed in the uniqueness theorem of \cite%
{Ber}, we impose in our uniqueness Theorem 7.1 in section 7 a less
restrictive assumption (\ref{2.7}) on the unknown coefficient $a\left( 
\mathbf{x}\right) $.

Let $u\left( \mathbf{x},t\right) $ be the solution of Cauchy problem (\ref%
{2.15}), (\ref{2.16}) considered in Theorem 2.2. We approximate the function 
$u\left( \mathbf{x},\varepsilon \right) $ for a sufficiently small number $%
\varepsilon >0$ by the first term $u_{0}\left( \mathbf{x},\varepsilon
\right) $ of asymptotic expansion (\ref{2.20}) of the function $u\left( 
\mathbf{x},t\right) $ at $t\rightarrow 0^{+}.$ From that step on we work
with the function $u\left( \mathbf{x},t\right) $ only on the interval $t\in
\left( \varepsilon ,T\right) .$ The number $\varepsilon $ is chosen
numerically.

Next, we obtain a nonlinear integral differential equation with Volterra
integrals \ 
\begin{equation}
\displaystyle\int\limits_{\varepsilon }^{t}\left( ...\right) ds,\text{ }t\in
\left( \varepsilon ,T\right)  \label{3.02}
\end{equation}%
in it$.$ This equation is supplied by the Dirichlet boundary condition at $%
\partial \Omega \times \left( \varepsilon ,T\right) $ and the Neumann
boundary condition at $\Gamma _{0}\times \left( \varepsilon ,T\right) ,$ see
(\ref{2.2}) for $\partial \Omega $ and $\Gamma _{0}.$ If the resulting
boundary value problem (BVP) is solved numerically, then the unknown
coefficient can be approximated by an explicit formula.

To solve that BVP, we approximate the $t-$derivative of the solution of that
integral differential equation via finite differences assuming that the grid
step size $h$ satisfies the following inequality%
\begin{equation}
h\geq h_{0}\in \left( 0,1\right) ,  \label{1.1}
\end{equation}%
where the number $h_{0}$ is fixed. In addition, we approximate Volterra
integrals (\ref{3.02}) via their discrete analogs generated by those finite
differences. As a result, we obtain a BVP for a system of coupled nonlinear
elliptic PDEs with the Dirichlet boundary conditions at $\partial \Omega $
and the Neumann boundary conditions at $\Gamma _{0}.$

To solve the latter BVP, we minimize the corresponding weighted least
squares cost functional with the CWF in it, as indicated in subsection 3.3.
Since this functional is new, then our detailed convergence analysis of
section 6 is new as well. That analysis ensures the global strong convexity
of our functional and provides accuracy estimates of its minimizer for both
noisy and noiseless data. The convergence analysis ends up with a theorem,
which ensures the global convergence of the gradient descent method of the
minimization of that functional.

In fact, we work here with an approximate mathematical model for CIP2. This
model consists of two new elements. Its first element is the above mentioned
approximation of the function $u\left( \mathbf{x},\varepsilon \right) $ by
the function $u_{0}\left( \mathbf{x},\varepsilon \right) .$ The second
element is our assumption about finite differences with respect to $t$ with
condition (\ref{1.1}). In section 7, we prove uniqueness Theorem 7.1, which
ensures uniqueness of CIP2 within the framework of this approximate
mathematical model.

\subsection{The concept of the convexification numerical method}

\label{sec:3.3}

We develop a new version of the above mentioned convexification globally
convergent numerical method. The convexification method was first proposed
in \cite{Klib95,Klib97} with the goal to construct globally convergent
numerical procedures for CIPs. The main advantage of this method is its
global convergence property in terms of Definition 1.1. In particular,
convexification does not face the well known phenomenon of multiple local
minima and ravines of conventional least squares mismatch functionals for
CIPs, see, e.g. \cite{Beilina1,Beilina2,Chavent,Gonch1,Gonch2} for those
functionals. We refer to \cite{Scales}\ for a numerical example of multiple
local minima. The originating works \cite{Klib95,Klib97} were purely
theoretical ones. More recently, starting from the publication \cite{Bak},
several versions of the convexification method for many CIPs were developed
in a number of publications, in which the theory is supported by numerical
studies. We refer to, e.g. \cite{Boul,KLZ,KL,Klibgrad,MFG1} as well as to
references cited therein for some samples of those publications.

The convexification method can be considered as a numerical development of
the initially purely theoretical idea of \cite{BukhKlib}. The key element of
the convexification is the presence of the so-called \textquotedblleft
Carleman Weight Function" (CWF) in a certain weighted least squares cost
functional $J$. The CWF is the weight function, which is used in the
Carleman estimate for the corresponding Partial Differential Operator. The
presence of the CWF ensures the strong convexity of $J$ on an \emph{a priory}
chosen convex bounded set $W\left( d\right) \subset H$ of the diameter $d>0,$
where $d$ is a fixed arbitrary number and $H$ is an appropriate Hilbert
space. Since a smallness condition is not imposed on $d$, then $J$ is
globally strongly convex on $W\left( d\right) $. The latter implies the
global convergence, in terms of Definition 1.1, of the gradient descent
method of the minimization of the functional $J$. Furthermore, it is
demonstrated computationally in the numerical Test 8.3 in section 8 that the
presence of the CWF in the functional $J$ affects quite significantly the
quality of numerical results.

\section{The Transformation Procedure}

\label{sec:4}

By (\ref{2.18}) and (\ref{2.20}), we can choose the number $T>0$ so small
that 
\begin{equation}
u\left( \mathbf{x},t\right) \geq \frac{1}{2}u_{0}\left( \mathbf{x},t\right) ,%
\text{ }\left( \mathbf{x},t\right) \in Q_{T}.  \label{3.1}
\end{equation}

\subsection{The first step of the transformation procedure}

\label{sec:4.1}

On this step we obtain an integral differential equation, which does not
contain the target unknown coefficient $a\left( \mathbf{x}\right) .$

Let $\varepsilon \in \left( 0,T\right) $ be a sufficiently small number,
which will be found numerically: we refer to Test 8.2 in section 8. Using
Theorem 2.2 and (\ref{2.18}), we approximate the function $u\left( \mathbf{x}%
,\varepsilon \right) $ by the first term of asymptotic expansion (\ref{2.20}%
) via setting 
\begin{equation}
u\left( \mathbf{x},\varepsilon \right) =u_{0}\left( \mathbf{x},\varepsilon
\right) =\frac{1}{\left( 2\sqrt{\pi \varepsilon }\right) ^{n}}\exp \left( -%
\frac{\left\vert \mathbf{x}\right\vert ^{2}}{4\varepsilon }\right) .
\label{3.2}
\end{equation}%
Formula (\ref{3.2}) is the approximation for the function $u\left( \mathbf{x}%
,\varepsilon \right) ,$ \ which is kept unchanged in all our follow up
derivations. In fact, formula (\ref{3.2}) is the first element of our
approximate mathematical model mentioned in subsection 3.2. Below we
consider the values of the variable $t$ only in the interval $t\in \left(
\varepsilon ,T\right) $.

Following (\ref{2.3}), denote%
\begin{equation}
\left. 
\begin{array}{c}
\text{ }Q_{\varepsilon ,T}=\Omega \times \left( \varepsilon ,T\right) , \\ 
\Gamma _{0,\varepsilon ,T}=\Gamma _{0}\times \left( \varepsilon ,T\right) ,%
\text{ }\partial \Omega _{\varepsilon ,T}=\text{ }\partial \Omega \times
\left( \varepsilon ,T\right) .%
\end{array}%
\right.  \label{3.3}
\end{equation}%
It follows from (\ref{3.1})-(\ref{3.3}) that we can consider the function $%
w\left( \mathbf{x},t\right) ,$%
\begin{equation}
w\left( \mathbf{x},t\right) =\ln \left( u\left( \mathbf{x},t\right) \right) ,%
\text{ }\left( \mathbf{x},t\right) \in Q_{\varepsilon ,T}.  \label{3.4}
\end{equation}%
Substituting (\ref{3.4}) in (\ref{2.15}) and (\ref{2.17}) and using (\ref%
{3.2}) and (\ref{3.3}), we obtain%
\begin{equation}
w_{t}-\Delta w-\left\vert \nabla w\right\vert ^{2}=a\left( \mathbf{x}\right) 
\text{ in }Q_{\varepsilon ,T},  \label{3.5}
\end{equation}%
\begin{equation}
w\left( \mathbf{x},t\right) \mid _{\partial \Omega _{\varepsilon ,T}}=\ln
\left( g_{0}\left( \mathbf{x},t\right) \right) =s_{0}\left( \mathbf{x}%
,t\right) ,\text{ }  \label{3.6}
\end{equation}%
\begin{equation*}
w_{x_{1}}\left( \mathbf{x},t\right) \mid _{\Gamma _{0,\varepsilon ,T}}=\frac{%
g_{1}\left( \mathbf{x},t\right) }{g_{0}\left( \mathbf{x},t\right) }%
=s_{1}\left( \mathbf{x},t\right) .
\end{equation*}%
In addition, using (\ref{3.2}) and (\ref{3.4}), we obtain%
\begin{equation*}
w\left( \mathbf{x},\varepsilon \right) =\ln \left( u_{0}\left( \mathbf{x}%
,t\right) \right) =\ln \left[ \frac{1}{\left( 2\sqrt{\pi \varepsilon }%
\right) ^{n}}\exp \left( -\frac{\left\vert \mathbf{x}\right\vert ^{2}}{%
4\varepsilon }\right) \right] .
\end{equation*}%
Hence, 
\begin{equation}
w\left( \mathbf{x},\varepsilon \right) =-\frac{\left\vert \mathbf{x}%
\right\vert ^{2}}{4\varepsilon }-\ln \left( 2\sqrt{\pi \varepsilon }\right)
^{n}.  \label{3.8}
\end{equation}%
Denote 
\begin{equation}
v\left( \mathbf{x},t\right) =\partial _{t}w\left( \mathbf{x},t\right) .
\label{3.9}
\end{equation}%
By (\ref{3.9})%
\begin{equation}
w\left( \mathbf{x},t\right) =\displaystyle\int\limits_{\varepsilon
}^{t}v\left( \mathbf{x},s\right) ds+w\left( \mathbf{x},\varepsilon \right) ,
\label{3.10}
\end{equation}%
where $w\left( \mathbf{x},\varepsilon \right) $ is given in (\ref{3.8}).
Differentiate both sides of equation (\ref{3.5}) with respect to $t$. Using 
\begin{equation*}
\frac{\partial a\left( \mathbf{x}\right) }{\partial t}\equiv 0
\end{equation*}%
as well as (\ref{3.6})-(\ref{3.10}), we obtain%
\begin{equation}
v_{t}-\Delta v-2\nabla v\left( \displaystyle\int\limits_{\varepsilon
}^{t}\nabla v\left( \mathbf{x},s\right) ds+\nabla w\left( \mathbf{x}%
,\varepsilon \right) \right) =0,\text{ }\left( \mathbf{x},t\right) \in
Q_{\varepsilon ,T},  \label{3.11}
\end{equation}%
\begin{equation}
v\left( \mathbf{x},t\right) \mid _{\partial \Omega _{\varepsilon
,T}}=\partial _{t}s_{0}\left( \mathbf{x},t\right) ,  \label{3.12}
\end{equation}%
\begin{equation}
v_{x_{1}}\left( \mathbf{x},t\right) \mid _{\Gamma _{0,\varepsilon
,T}}=\partial _{t}s_{1}\left( \mathbf{x},t\right) .  \label{3.13}
\end{equation}

We have obtained the desired integral differential equation (\ref{3.11})
with the Dirichlet boundary condition (\ref{3.12}) on the whole lateral
boundary $\partial \Omega _{\varepsilon ,T}$ of the time cylinder $%
Q_{\varepsilon ,T}$ and Neumann boundary condition (\ref{3.13}) on a part $%
\Gamma _{0,\varepsilon ,T}\subset \partial \Omega _{\varepsilon ,T}$.
Suppose that a solution $v\left( \mathbf{x},t\right) \in H^{2,1}\left(
Q_{\varepsilon ,T}\right) $ of problem (\ref{3.11})-(\ref{3.13}) is known.
Then the target coefficient $a\left( \mathbf{x}\right) $ can be computed.
Indeed, integrating (\ref{3.5}) with respect to $t\in \left( \varepsilon
,T\right) ,$ we obtain%
\begin{equation*}
a\left( \mathbf{x}\right) =\frac{1}{T-\varepsilon }\int\limits_{\varepsilon
}^{T}w_{s}\left( \mathbf{x},s\right) ds-\frac{1}{T-\varepsilon }%
\int\limits_{\varepsilon }^{T}\left( \Delta w+\left\vert \nabla w\right\vert
^{2}\right) \left( \mathbf{x},s\right) ds.
\end{equation*}%
Hence, 
\begin{equation}
a\left( \mathbf{x}\right) =\frac{w\left( \mathbf{x},T\right) -w\left( 
\mathbf{x},\varepsilon \right) }{T-\varepsilon }-\frac{1}{T-\varepsilon }%
\int\limits_{\varepsilon }^{T}\left( \Delta w+\left\vert \nabla w\right\vert
^{2}\right) \left( \mathbf{x},s\right) ds,\text{ }\mathbf{x}\in \Omega .
\label{3.15}
\end{equation}

Thus, it follows from (\ref{3.10}) and (\ref{3.15}) that we need to focus
below on the numerical solution of problem (\ref{3.11})-(\ref{3.13}).

\subsection{The second step of the transformation procedure}

\label{sec:4.2}

On this step we replace problem (\ref{3.11})-(\ref{3.13}) with a boundary
value problem for a system of coupled nonlinear elliptic equations. To do
this, we assume that the $t-$derivatives in (\ref{3.11})-(\ref{3.13}) are
written in finite differences with the grid step size $h$ satisfying (\ref%
{1.1}). Respectively, we write Volterra integrals in (\ref{3.11}) in the
discrete form. These two approximations form the second element of our
approximate mathematical model mentioned in subsection 3.2.

The assumption about the finite differences for the $t-$derivative for some
CIPs for parabolic equations was previously used in \cite{TFD}. However,
since the most challenging case of the function $\delta \left( \mathbf{x}%
\right) $ as the initial condition (\ref{2.16}) was not considered in \cite%
{TFD}, then CIPs of \cite{TFD} are significantly simpler than the one we
study here. In addition, unlike the current paper, numerical experiments
were not carried out in \cite{TFD}.

Consider the partition of the interval $\left[ \varepsilon ,T\right] $ in $%
k\geq 3$ subintervals with the grid step size $h$ satisfying (\ref{1.1}), 
\begin{equation}
\varepsilon =t_{0}<t_{1}<...<t_{k-1}<t_{k}=T,\text{ }t_{i}-t_{i-1}=h,\text{ }%
i=1,...,k.  \label{3.16}
\end{equation}%
Denote%
\begin{equation}
Y=\left\{ t_{i}\right\} _{i=0}^{k}.  \label{3.17}
\end{equation}%
Using (\ref{3.17}), define semi-discrete analogs of sets in (\ref{3.3}): 
\begin{equation}
\left. 
\begin{array}{c}
Q_{\varepsilon ,h,T}=\Omega \times Y, \\ 
\Gamma _{0,\varepsilon ,h,T}=\Gamma _{0,\varepsilon ,T}\times Y\text{, }%
\partial \Omega _{\varepsilon ,h,T}=\partial \Omega \times Y.%
\end{array}%
\right.  \label{3.18}
\end{equation}%
Therefore, by (\ref{3.16}) the function $v\left( \mathbf{x},t\right) $
becomes a $\left( k+1\right) -$dimensional vector function $V\left( \mathbf{x%
}\right) $, 
\begin{equation*}
V\left( \mathbf{x}\right) \emph{=}\left( v\left( \mathbf{x},t_{0}\right)
,v\left( \mathbf{x},t_{1}\right) ,...,v\left( \mathbf{x},t_{k}\right)
\right) ^{T}\text{ for }\left( \mathbf{x},t\right) \in Q_{\varepsilon ,h,T}.
\end{equation*}%
We define finite difference derivatives of $v\left( \mathbf{x},t_{j}\right)
, $ $j=0,...,k-1$ with respect to $t$ as:%
\begin{equation}
\partial _{t}^{h}v\left( \mathbf{x},t_{0}\right) =\frac{v\left( \mathbf{x}%
,t_{1}\right) -v\left( \mathbf{x},t_{0}\right) }{h},  \label{3.20}
\end{equation}%
\begin{equation*}
\partial _{t}^{h}v\left( \mathbf{x},t_{i}\right) =\frac{v\left( \mathbf{x}%
,t_{i+1}\right) -v\left( \mathbf{x},t_{i}\right) }{h},\text{ }i=1,...,k-1.
\end{equation*}

We now define the finite difference derivative $\partial _{t}^{h}v\left( 
\mathbf{x},t_{k}\right) $ at the end point of the interval $\left[ 0,T\right]
$. In fact, we want to form such a system of $\left( k+1\right) $ elliptic
PDEs that each equation number $i=0,...,k$ of this system would contain the
term $\Delta v\left( \mathbf{x},t_{i}\right) .$ Furthermore, this term needs
to be involved only once:\ in \textquotedblleft its own" equation number $i$%
. The reason of the latter is that the Carleman estimate of Theorem 5.1
(below) works only for the Laplace operator. This is why we use below a bit
unusual approximation by finite differences of the derivative $\partial
_{t}^{h}v\left( \mathbf{x},t_{k}\right) $. Although some other
approximations of this nature can also be used, they will not change our
theory, as long the above conditions of the involvement of the terms $\Delta
v\left( \mathbf{x},t_{i}\right) $ are met. Using Taylor formula, we obtain
for any function $f\in C^{2}\left[ 0,T\right] $%
\begin{equation*}
f^{\prime }\left( t_{k}\right) =\frac{3f\left( t_{k}\right) -f\left(
t_{k-1}\right) -f\left( t_{k-2}\right) -f\left( t_{k-3}\right) }{6h}+O\left(
h\right) ,\text{ as }h\rightarrow 0^{+}.
\end{equation*}%
Hence, we define $\partial _{t}^{h}v\left( \mathbf{x},t_{k}\right) $ as%
\begin{equation}
\partial _{t}^{h}v\left( \mathbf{x},t_{k}\right) =\frac{3v\left( \mathbf{x}%
,t_{k}\right) -v\left( \mathbf{x},t_{k-1}\right) -v\left( \mathbf{x}%
,t_{k-2}\right) -v\left( \mathbf{x},t_{k-3}\right) }{6h}.  \label{3.22}
\end{equation}

It follows from (\ref{3.16}) that we assume a uniform time grid%
\begin{equation*}
t_{j}=\varepsilon +jh,\text{ }j=0,...,k,
\end{equation*}%
where $t_{0}=\varepsilon $ and $h$ is the above time step. In addition, we
define the discrete analog of the Volterra integral in (\ref{3.11}). Using
the right-rectangle quadratic rule, we approximate the integral as follows: 
\begin{equation}
\int\limits_{\varepsilon }^{t_{i}}\nabla v\left( \mathbf{x},s\right)
ds\approx h\sum\limits_{j=0}^{i}\nabla v\left( \mathbf{x},t_{j}\right) ,%
\text{ }i\geq 1.  \label{3.23}
\end{equation}%
For $i=0$ we have%
\begin{equation}
\int\limits_{\varepsilon }^{t_{0}}\nabla v\left( \mathbf{x},s\right) ds=0.
\label{3.24}
\end{equation}%
It follows from (\ref{3.20})-(\ref{3.24}) that equation (\ref{3.11}) can be
written via finite differences with respect to $t$ as the following system
of $\left( k+1\right) $ coupled elliptic equations:%
\begin{equation}
\left. 
\begin{array}{c}
L_{0}\left( v\left( \mathbf{x},t_{0}\right) ,v\left( \mathbf{x},t_{1}\right)
\right) = \\ 
=\Delta v\left( \mathbf{x},t_{0}\right) +2\nabla v\left( \mathbf{x}%
,t_{0}\right) \nabla w\left( \mathbf{x},\varepsilon \right) + \\ 
+\left( v\left( \mathbf{x},t_{0}\right) -v\left( \mathbf{x},t_{1}\right)
\right) /h=0,\text{ }\mathbf{x}\in \Omega ,%
\end{array}%
\right.  \label{3.25}
\end{equation}%
\begin{equation}
\left. 
\begin{array}{c}
L_{i}\left( v\left( \mathbf{x},t_{0}\right) ,v\left( \mathbf{x},t_{1}\right)
,...,v\left( \mathbf{x},t_{i+1}\right) \right) = \\ 
=\Delta v\left( \mathbf{x},t_{i}\right) +2\nabla v\left( \mathbf{x}%
,t_{i}\right) \nabla w\left( \mathbf{x},\varepsilon \right) +2h\nabla
v\left( \mathbf{x},t_{i}\right) \left( \sum\limits_{j=0}^{i}\nabla v\left( 
\mathbf{x},t_{j}\right) \right) + \\ 
+\left( v\left( \mathbf{x},t_{i}\right) -v\left( \mathbf{x},t_{i+1}\right)
\right) /h=0, \\ 
\text{ }i=1,...,k-1,\text{ }\mathbf{x}\in \Omega .%
\end{array}%
\right.  \label{3.29}
\end{equation}%
Also, we have for $i=k$%
\begin{equation}
\left. 
\begin{array}{c}
L_{k}\left( v\left( \mathbf{x},t_{0}\right) ,v\left( \mathbf{x},t_{1}\right)
,v\left( \mathbf{x},t_{2}\right) ,...,v\left( \mathbf{x},t_{k-1}\right)
,v\left( \mathbf{x},t_{k}\right) \right) = \\ 
=\Delta v\left( \mathbf{x},t_{k}\right) +2\nabla v\left( \mathbf{x}%
,t_{k}\right) \nabla w\left( \mathbf{x},\varepsilon \right) +2h\nabla
v\left( \mathbf{x},t_{k}\right) \left( \sum\limits_{j=0}^{k}\nabla v\left( 
\mathbf{x},t_{j}\right) \right) + \\ 
+\left( v\left( \mathbf{x},t_{k-1}\right) +v\left( \mathbf{x},t_{k-2}\right)
+v\left( \mathbf{x},t_{k-3}\right) -3v\left( \mathbf{x},t_{k}\right) \right)
/\left( 6h\right) =0,\text{ }\mathbf{x}\in \Omega .%
\end{array}%
\right.  \label{3.30}
\end{equation}%
The boundary conditions for this system are derived from (\ref{3.12}), (\ref%
{3.13}), (\ref{3.17}), the second line of (\ref{3.18}) and (\ref{3.20})-(\ref%
{3.22}). These are Dirichlet boundary conditions on the entire boundary $%
\partial \Omega $ and Neumann boundary conditions at the part $\Gamma
_{0}\subset \partial \Omega $. More precisely, 
\begin{equation}
v\left( \mathbf{x},t_{i}\right) =\frac{g_{0}\left( \mathbf{x},t_{i+1}\right)
-g_{0}\left( \mathbf{x},t_{i}\right) }{h},\text{ }\mathbf{x}\in \partial
\Omega ,\text{ }i=0,...,k-1,  \label{3.31}
\end{equation}%
\begin{equation*}
v\left( \mathbf{x},t_{k}\right) =\frac{3g_{0}\left( \mathbf{x},t_{k}\right)
-g_{0}\left( \mathbf{x},t_{k-1}\right) -g_{0}\left( \mathbf{x}%
,t_{k-2}\right) -g_{0}\left( \mathbf{x},t_{k-3}\right) }{6h},\text{ }\mathbf{%
x}\in \partial \Omega ,
\end{equation*}%
\begin{equation*}
v_{x_{1}}\left( \mathbf{x},t_{i}\right) =\frac{g_{1}\left( \mathbf{x}%
,t_{i+1}\right) -g_{1}\left( \mathbf{x},t_{i}\right) }{h},\text{ }\mathbf{x}%
\in \Gamma _{0},\text{ }i=0,...,k-1,
\end{equation*}%
\begin{equation}
v_{x_{1}}\left( \mathbf{x},t_{k}\right) =\frac{3g_{1}\left( \mathbf{x}%
,t_{k}\right) -g_{1}\left( \mathbf{x},t_{k-1}\right) -g_{1}\left( \mathbf{x}%
,t_{k-2}\right) -g_{1}\left( \mathbf{x},t_{k-3}\right) }{6h},\text{ }\mathbf{%
x}\in \Gamma _{0}.  \label{3.34}
\end{equation}

Therefore, we have obtained the following problem:

\textbf{Problem. }Find the $\left( k+1\right) -$dimensional vector function $%
V\left( \mathbf{x}\right) ,$ 
\begin{equation}
V\left( \mathbf{x}\right) =\left( v\left( \mathbf{x},t_{0}\right) ,v\left( 
\mathbf{x},t_{1}\right) ,...,v\left( \mathbf{x},t_{k}\right) \right) ^{T},%
\text{ }\mathbf{x}\in \Omega  \label{3.35}
\end{equation}%
satisfying conditions (\ref{3.25})-(\ref{3.34}).

Suppose that the vector function $V\left( \mathbf{x}\right) $ in (\ref{3.35}%
) is computed. Denote it as $V_{\text{comp}}\left( \mathbf{x}\right) $ 
\begin{equation}
V_{\text{comp}}\left( \mathbf{x}\right) =\left( v_{\text{comp}}\left( 
\mathbf{x},t_{0}\right) ,v_{\text{comp}}\left( \mathbf{x},t_{1}\right)
,...,v_{\text{comp}}\left( \mathbf{x},t_{k}\right) \right) ^{T},\text{ }%
\mathbf{x}\in \Omega .  \label{3.350}
\end{equation}%
Then, using obvious analogs of (\ref{3.10}), (\ref{3.23}) and (\ref{3.24}),
we compute the functions $w_{\text{comp}}\left( \mathbf{x},t_{i}\right) ,$ 
\begin{equation}
w_{\text{comp}}\left( \mathbf{x},t_{i}\right) =h\displaystyle%
\sum\limits_{j=0}^{i}v_{\text{comp}}\left( \mathbf{x},t_{j}\right) +w\left( 
\mathbf{x},\varepsilon \right) ,  \label{3.36}
\end{equation}%
where $w\left( \mathbf{x},\varepsilon \right) $ is given in (\ref{3.8}).
Next, use (\ref{3.36}) and the following analog of (\ref{3.15}): 
\begin{equation*}
a_{\text{comp}}\left( \mathbf{x}\right) =\frac{1}{T-\varepsilon }\left( w_{%
\text{comp}}\left( \mathbf{x},t_{k}\right) -w\left( \mathbf{x},\varepsilon
\right) \right) -
\end{equation*}%
\begin{equation}
-\frac{h}{T-\varepsilon }\sum\limits_{i=0}^{k}\left( \Delta w_{\text{comp}%
}\left( \mathbf{x},t_{i}\right) +\left\vert \nabla w_{\text{comp}}\left( 
\mathbf{\ x},t_{i}\right) \right\vert ^{2}\right) .  \label{3.37}
\end{equation}

\section{Convexification Method for the Numerical Solution of Problem (%
\protect\ref{3.25})-(\protect\ref{3.35})}

\label{sec:5}

Below $C=C\left( \Omega \right) >0$ denotes different positive numbers
depending only on the domain $\Omega $ defined in (\ref{2.1}). We need
functions $v\left( \mathbf{x},t_{i}\right) \in C^{1}\left( \overline{\Omega }%
\right) \cap H^{2}\left( \Omega \right) .$ Since we work in $\mathbb{R}^{n},$
then by Sobolev embedding theorem 
\begin{equation}
\left. 
\begin{array}{c}
H^{m_{n}}\left( \Omega \right) \subset C^{1}\left( \overline{\Omega }\right)
,\text{ }m_{n}=\left[ n/2\right] +2, \\ 
\left\Vert f\right\Vert _{C^{1}\left( \overline{\Omega }\right) }\leq
C\left\Vert f\right\Vert _{H^{m_{n}}\left( \Omega \right) },\text{ }\forall
f\in H^{m_{n}}\left( \Omega \right) .%
\end{array}%
\right.  \label{4.1}
\end{equation}%
Thus, in the most popular cases of $n=2,3$ we have $m_{n}=3.$ Let $H$ be a
Hilbert space with its norm $\left\Vert .\right\Vert _{H}.$ Then we define
the space $H_{k+1}$ as%
\begin{equation*}
H_{k+1}=\left\{ Q=\left( q_{0},...,q_{k}\right) ^{T}:q_{i}\in H,\text{ }%
\left\Vert Q\right\Vert _{H_{k+1}}=\left( \sum\limits_{i=0}^{k}\left\Vert
q_{i}\right\Vert _{H}^{2}\right) ^{1/2}\right\} .
\end{equation*}

Introduce two spaces:%
\begin{equation}
\left. 
\begin{array}{c}
H_{0}^{2}\left( \Omega \right) =\left\{ q\in H^{2}\left( \Omega \right)
:q\mid _{\partial \Omega }=0,\text{ }q_{x_{1}}\mid _{\Gamma _{0}}=0\right\} ,
\\ 
H_{0}^{m_{n}}\left( \Omega \right) =\left\{ q\in H^{m_{n}}\left( \Omega
\right) :\text{ }q\mid _{\partial \Omega }=0,\text{ }q_{x_{1}}\mid _{\Gamma
_{0}}=0\right\} .%
\end{array}%
\right.  \label{4.01}
\end{equation}%
It is convenient to consider below the vector functions $V\left( \mathbf{x}%
\right) $ in the following form%
\begin{equation}
V\left( \mathbf{x}\right) =\left( v_{0}\left( \mathbf{x}\right) ,v_{1}\left( 
\mathbf{x}\right) ,...,v_{k}\left( \mathbf{x}\right) \right) ^{T},
\label{4.100}
\end{equation}%
rather than the one in (\ref{3.35}). Nevertheless, as soon as a vector
function $V\left( \mathbf{x}\right) $ of the form (\ref{4.100}) is computed,
i.e. the vector function

\begin{equation*}
V_{\text{comp}}\left( \mathbf{x}\right) =\left( v_{0,\text{comp}}\left( 
\mathbf{x}\right) ,v_{1,\text{comp}}\left( \mathbf{x}\right) ,...,v_{k,\text{%
comp}}\left( \mathbf{x}\right) \right) ^{T}
\end{equation*}%
is found, we naturally assign below as in (\ref{3.350}): 
\begin{equation*}
\left. 
\begin{array}{c}
V_{\text{comp}}\left( \mathbf{x}\right) =\left( v_{0,\text{comp}}\left( 
\mathbf{x}\right) ,v_{1,\text{comp}}\left( \mathbf{x}\right) ,...,v_{k,\text{
comp}}\left( \mathbf{x}\right) \right) ^{T}= \\ 
=\left( v_{\text{comp}}\left( \mathbf{x},t_{0}\right) ,v_{\text{comp}}\left( 
\mathbf{x},t_{1}\right) ,...,v_{\text{comp}}\left( \mathbf{x},t_{k}\right)
\right) ^{T},\text{ }\mathbf{x}\in \Omega .%
\end{array}%
\right.
\end{equation*}

Let $R>0$ be an arbitrary number. We consider the following set of $\left(
k+1\right) -$dimensional vector functions%
\begin{equation}
B\left( R\right) =\left\{ 
\begin{array}{c}
V\left( \mathbf{x}\right) =\left( v_{0}\left( \mathbf{x}\right) ,v_{1}\left( 
\mathbf{x}\right) ,...,v_{k}\left( \mathbf{x}\right) \right) ^{T}\in
H_{k+1}^{m_{n}}\left( \Omega \right) : \\ 
\left\Vert v_{i}\left( \mathbf{x}\right) \right\Vert _{H^{m_{n}}\left(
\Omega \right) }<R,\text{ }i=0,...,k, \\ 
\text{functions }v_{i}\left( \mathbf{x}\right) \text{ satisfy boundary
conditions } \\ 
\text{(\ref{3.31})-(\ref{3.34}) for corresponding indexes }i%
\end{array}%
\right\} .  \label{4.2}
\end{equation}%
By (\ref{4.1}) all functions $v_{i}\left( \mathbf{x}\right) $ in (\ref{4.2})
have the following properties:%
\begin{equation*}
v_{i}\left( \mathbf{x}\right) \in C^{1}\left( \overline{\Omega }\right) ,%
\text{ }\left\Vert v\left( \mathbf{x}\right) \right\Vert _{C^{1}\left( 
\overline{\Omega }\right) }\leq CR.
\end{equation*}

Let $\lambda \geq 1$ be a parameter. Define the Carleman Weight Function
(CWF) as 
\begin{equation}
\varphi _{\lambda }\left( x_{1}\right) =e^{2\lambda x_{1}^{2}}.  \label{4.4}
\end{equation}

An analog of Theorem 5.1 was proven in \cite{KLZ} and \cite[Theorem 9.4.1]%
{KL} for the case of the parabolic operator $\partial _{t}-\Delta $ and with
the different CWF $\psi _{\lambda }\left( x\right) =e^{2\lambda \left(
x_{1}^{2}-t^{2}\right) }.$ Both the formulation and the proof of Theorem 5.1
follow immediately from the parabolic case if assuming the $t-$independence
of all involved functions.

\textbf{Theorem 5.1} (Carleman estimate). \emph{Let }$\varphi _{\lambda
}\left( x_{1}\right) $\emph{\ be the function defined in (\ref{4.4}) and let 
}$\Omega $\emph{\ be the domain defined in (\ref{2.1}). There exists a
sufficiently large number }$\lambda _{0}=\lambda _{0}\left( \Omega \right)
\geq 1$\emph{\ such that the following Carleman estimate holds:}%
\begin{equation*}
\displaystyle\int\limits_{\Omega }\left( \Delta u\right) ^{2}\varphi
_{\lambda }\left( x_{1}\right) d\mathbf{x}\geq \frac{C}{\lambda }%
\displaystyle\int\limits_{\Omega }\left(
\sum_{i,j=1}^{n}u_{x_{i}x_{j}}^{2}\right) \varphi _{\lambda }\left(
x_{1}\right) d\mathbf{x+}
\end{equation*}%
\begin{equation}
+C\displaystyle\int\limits_{\Omega }\left( \lambda \left\vert \nabla
u\right\vert ^{2}+\lambda ^{3}u^{2}\right) \varphi _{\lambda }\left(
x_{1}\right) d\mathbf{x,}  \label{4.40}
\end{equation}%
\begin{equation*}
\forall u\in H_{0}^{2}\left( \Omega \right) ,\text{ }\forall \lambda \geq
\lambda _{0}.
\end{equation*}

\textbf{Remark 5.1.} A CWF, which is more general than the one in (\ref{4.4}%
), can also be chosen, and this CWF would be involved in a Carleman estimate
for a more general domain $\Omega $ than the one in (\ref{2.1}).
Furthermore, that CWF would provide a Carleman estimate for a general
elliptic operator of the second order with variable coefficients in its
principal part, rather than for the Laplace operator only as in Theorem 5.1.
In this regard we refer to, e.g. \cite[formula (2.30)]{KL} and \cite[\S 1 of
chapter 4]{LRS}. That more general CWF depends on two large parameters.
However, we use CWF (\ref{4.4}) below because the dependence of the CWF of 
\cite{KL,LRS} on two large parameters would make our goal of numerical
studies significantly more complicated than the case of just one large
parameter in (\ref{4.4}). We also refer to Remark 2.1.

The convexification weighted functional for problem (\ref{3.25})-(\ref{3.35}%
) is: 
\begin{equation}
J_{\lambda ,\alpha }:\overline{B\left( R\right) }\rightarrow \mathbb{R},
\label{4.5}
\end{equation}%
\begin{equation}
\left. 
\begin{array}{c}
J_{\lambda ,\alpha }\left( V\right) =e^{-2\lambda
c}\sum\limits_{i=0}^{k-1}\int\limits_{\Omega }\left[ L_{i}\left( v_{0}\left( 
\mathbf{x}\right) ,v_{1}\left( \mathbf{x}\right) ,...,v_{i+1}\left( \mathbf{x%
}\right) \right) \right] ^{2}\varphi _{\lambda }\left( x_{1}\right) d\mathbf{%
x+} \\ 
+e^{-2\lambda c}\int\limits_{\Omega }\left[ L_{k}\left( v_{0}\left( \mathbf{%
\ x}\right) ,v_{1}\left( \mathbf{x}\right) ,...,v_{k}\left( \mathbf{x}%
\right) \right) \right] ^{2}\varphi _{\lambda }\left( x_{1}\right) d\mathbf{%
x+}\alpha \left\Vert V\right\Vert _{H_{k+1}^{m_{n}}\left( \Omega \right)
}^{2}.%
\end{array}%
\right.  \label{4.6}
\end{equation}%
Here $\lambda \geq \lambda _{0},$ where $\lambda _{0}$ is the number of
Theorem 4.1, $\alpha \in \left( 0,1\right) $ is the regularization
parameter, and $c>0$ is a number to be chosen numerically. Indeed,%
\begin{equation}
\max_{\overline{\Omega }}\varphi _{\lambda }\left( x_{1}\right) =e^{2\lambda
B^{2}}>>1.  \label{4.60}
\end{equation}%
Hence, since $\alpha \in \left( 0,1\right) ,$ then we need to balance in (%
\ref{4.6}) integral terms with the regularization term.

\section{Convergence Analysis}

\label{sec:6}

\subsection{The strong convexity of the functional $J_{\protect\lambda , 
\protect\alpha }$ on the set $\overline{B\left( R\right) }$}

\label{sec:6.1} ~\newline

Below $\left( \cdot ,\cdot \right) $\ denotes the scalar product in $%
H_{k+1}^{m_{n}}\left( \Omega \right) $ and $C_{1}=C_{1}\left( R,\alpha
,\varepsilon ,h_{0},T\right) >0$ denotes different numbers depending only on
listed parameters. Recall that the number $h_{0}$ is defined in (\ref{1.1}).
Hence, the number $k$ of grid points in (\ref{3.16}) ultimately depends on $%
\varepsilon ,h_{0}$ and $T$, $k=k\left( \varepsilon ,h_{0},T\right) .$

\textbf{Theorem 6.1} (strong convexity of the functional $J_{\lambda ,\alpha
}$). \emph{Let }$J_{\lambda ,\alpha }\left( V\right) $\emph{\ be the
functional defined in (\ref{4.5}), (\ref{4.6}). Then:}


\begin{enumerate}
\item \emph{At each point }$V\in \overline{B\left( R\right) }$\emph{\ there
exists Fr\'{e}chet derivative }$J_{\lambda ,\alpha }^{\prime }\left(
V\right) \in H_{0,k+1}^{m_{n}}\left( \Omega \right) $\emph{\ of this
functional. Furthermore, this derivative satisfies the Lipschitz continuity
condition} 
\begin{equation}
\left. 
\begin{array}{c}
\left\Vert J_{\lambda ,\alpha }^{\prime }\left( V_{1}\right) -J_{\lambda
,\alpha }^{\prime }\left( V_{2}\right) \right\Vert _{H_{k+1}^{m_{n}}\left(
\Omega \right) }\leq D\left\Vert V_{1}-V_{2}\right\Vert
_{H_{k+1}^{m_{n}}\left( \Omega \right) }^{2}, \\ 
\forall V_{1},V_{2}\in \overline{B\left( R\right) },%
\end{array}%
\right.  \label{5.1}
\end{equation}%
\emph{where the number }$D=D\left( R,\alpha ,\varepsilon ,h_{0},T,\lambda
\right) >0$\emph{\ depends only on listed parameters.}

\item \emph{There exists a sufficiently large number }$\lambda _{1}=\lambda
_{1}\left( R,\alpha ,\varepsilon ,h_{0},T\right) \geq \lambda _{0}\geq 1$%
\emph{\ depending only on listed parameters such that the functional }$%
J_{\lambda ,\alpha }\left( V\right) $\emph{\ is strongly convex on the set }$%
\overline{B\left( R\right) }$\emph{\ in (\ref{4.2}) for all }$\lambda \geq
\lambda _{1}.$\emph{\ \ More precisely, let }%
\begin{equation}
V_{1}\left( \mathbf{x}\right) =\left( v_{0,1}\left( \mathbf{x}\right)
,...,v_{k,1}\left( \mathbf{x}\right) \right) ^{T}\text{ \emph{and} }%
V_{2}\left( \mathbf{x}\right) =\left( v_{0,2}\left( \mathbf{x}\right)
,...,v_{k,2}\left( \mathbf{x}\right) \right) ^{T}  \label{5.100}
\end{equation}
\emph{be two arbitrary points of the set }$\overline{B\left( R\right) }$ 
\emph{. Then there exists such a number }$C_{1}>0$ \emph{that the following
estimate holds:}%
\begin{equation*}
J_{\lambda ,\alpha }\left( V_{2}\right) -J_{\lambda ,\alpha }\left(
V_{1}\right) -\left( J_{\lambda ,\alpha }^{\prime }\left( V_{1}\right)
,V_{2}-V_{1}\right) \geq
\end{equation*}%
\begin{equation}
\geq C_{1}\frac{e^{2\lambda \left( A^{2}-c\right) }}{\lambda }\displaystyle%
\sum\limits_{i,j=1}^{n}\left\Vert V_{2x_{i}x_{j}}-V_{1x_{i}x_{j}}\right\Vert
_{L_{2,k+1}\left( \Omega \right) }^{2}d\mathbf{x+}  \label{5.2}
\end{equation}%
\begin{equation*}
+C_{1}e^{2\lambda \left( A^{2}-c\right) }\left\Vert V_{2}-V_{1}\right\Vert
_{H_{k+1}^{1}\left( \Omega \right) }^{2}+\alpha \left\Vert
V_{2}-V_{1}\right\Vert _{H_{k+1}^{m_{n}}\left( \Omega \right) }^{2},
\end{equation*}%
\begin{equation*}
\forall \lambda \geq \lambda _{1}.
\end{equation*}

\item \emph{For each }$\lambda \geq \lambda _{1}$\emph{\ there exists unique
minimizer }$V_{\min ,\lambda }\in \overline{B\left( R\right) }$\emph{\ of
the functional }$J_{\lambda ,\alpha }\left( V\right) .$\emph{\ Furthermore,
the following inequality holds:} 
\begin{equation}
\left( J_{\lambda ,\alpha }^{\prime }\left( V_{\min ,\lambda }\right)
,V_{\min ,\lambda }-V\right) \leq 0,\text{ }\forall V\in \overline{B\left(
R\right) }.  \label{5.3}
\end{equation}
\end{enumerate}


\textbf{Remark 6.1.} The requirement of \ Theorem 6.1 that the parameter $%
\lambda $\ should be sufficiently large does not affect computations much.
Indeed, the optimal value of $\lambda =3$ is chosen in the numerical section
8.\ In addition, in our past works on computations for the convexification
method, which are cited above, the range of $\lambda $\ is $\lambda \in %
\left[ 1,5\right] .$\ This is similar with many asymptotic theories. Indeed,
any such theory claims that if a parameter $X$\ is sufficiently large, then
a certain formula $Z$\ is accurate. However, in the computational practice,
only specific numerical results can establish which exactly values of $X$\
ensure a good accuracy of $Z$.

\begin{proof}[Proof of Theorem 6.1]
Denote 
\begin{equation*}
P\left( \mathbf{x}\right) =V_{2}\left( \mathbf{x}\right) -V_{1}\left( 
\mathbf{x}\right) .
\end{equation*}%
Then 
\begin{equation}
\left. 
\begin{array}{c}
V_{2}\left( \mathbf{x}\right) =P\left( \mathbf{x}\right) +V_{1}\left( 
\mathbf{x}\right) , \\ 
P\left( \mathbf{x}\right) =\left( p_{0}\left( \mathbf{x}\right) ,p_{1}\left( 
\mathbf{x}\right) ,...,p_{k}\left( \mathbf{x}\right) \right) ^{T}\in
H_{0,k+1}^{m_{n}}\left( \Omega \right) .%
\end{array}%
\right.  \label{5.4}
\end{equation}%
Furthermore, triangle inequality and (\ref{4.2})\emph{\ }imply that%
\begin{equation}
P\left( \mathbf{x}\right) \in \overline{B_{0}\left( 2R\right) }=\left\{ 
\begin{array}{c}
Z\left( \mathbf{x}\right) =\left( z_{0}\left( \mathbf{x}\right) ,z_{1}\left( 
\mathbf{x}\right) ,...,z_{k}\left( \mathbf{x}\right) \right) ^{T}\in
H_{0,k+1}^{m_{n}}\left( \Omega \right) : \\ 
\left\Vert z_{i}\left( \mathbf{x}\right) \right\Vert _{H^{m_{n}}\left(
\Omega \right) }\leq 2R,\text{ }i=0,...,k.%
\end{array}%
\right\} .  \label{5.5}
\end{equation}

Using (\ref{3.29}), (\ref{4.100}), (\ref{5.100}) and (\ref{5.4}), we obtain%
\begin{equation}
\left. 
\begin{array}{c}
L_{i}\left( v_{0,2}\left( \mathbf{x}\right) ,v_{1,2}\left( \mathbf{x}\right)
,...,v_{i+1,2}\left( \mathbf{x}\right) \right) = \\ 
=L_{i}\left( v_{0,1}\left( \mathbf{x}\right) ,v_{1,1}\left( \mathbf{x}%
\right) ,...,v_{i+1,1}\left( \mathbf{x}\right) \right) + \\ 
+\Delta p_{i}\left( \mathbf{x}\right) +2\nabla p_{i}\left( \mathbf{x}\right)
\nabla w\left( \mathbf{x},\varepsilon \right) + \\ 
+2h\nabla p_{i}\left( \mathbf{x}\right) \left( \sum\limits_{j=0}^{i}\nabla
v_{j,1}\left( \mathbf{x}\right) \right) +2h\nabla v_{i,1}\left( \mathbf{x}%
\right) \left( \sum\limits_{j=0}^{i}\nabla p_{j}\left( \mathbf{x}\right)
\right) + \\ 
+2h\nabla p_{i}\left( \mathbf{x}\right) \sum\limits_{j=0}^{i}\nabla
p_{j}\left( \mathbf{x}\right) +\left( p_{i}\left( \mathbf{x}\right)
-p_{i+1}\left( \mathbf{x}\right) \right) /h,\text{ } \\ 
\mathbf{x}\in \Omega ,\text{ }i=0,...,k-1.%
\end{array}%
\right.  \label{5.6}
\end{equation}%
\newline
Hence,%
\begin{equation}
\left. 
\begin{array}{c}
\left[ L_{i}\left( v_{0,2}\left( \mathbf{x}\right) ,v_{1,2}\left( \mathbf{x}%
\right) ,...,v_{i+1,2}\left( \mathbf{x}\right) \right) \right] ^{2}- \\ 
-\left[ L_{i}\left( v_{0,1}\left( \mathbf{x}\right) ,v_{1,1}\left( \mathbf{x}%
\right) ,...,v_{i+1,1}\left( \mathbf{x}\right) \right) \right] ^{2}= \\ 
=M_{i,\text{lin}}\left( p_{0}\left( \mathbf{x}\right) ,p_{1}\left( \mathbf{x}%
\right) ,...,p_{i+1}\left( \mathbf{x}\right) \right) + \\ 
+M_{i,\text{nonlin}}\left( p_{0}\left( \mathbf{x}\right) ,p_{1}\left( 
\mathbf{x}\right) ,...,p_{i+1}\left( \mathbf{x}\right) \right) , \\ 
\text{ }\mathbf{x}\in \Omega ,\text{ }i=0,...,k-1.%
\end{array}%
\right.  \label{5.7}
\end{equation}%
where $M_{i,\text{lin}}$ and $M_{i,\text{nonlin}}$ depend on the vector
function $\left( p_{0}\left( \mathbf{x}\right) ,p_{1}\left( \mathbf{x}%
\right) ,...,p_{i+1}\left( \mathbf{x},t\right) \right) ^{T}$ linearly and
nonlinearly respectively. The formula for the term $M_{i,\text{lin}}$ is:%
\begin{equation}
\left. 
\begin{array}{c}
M_{i,\text{lin}}\left( p_{0}\left( \mathbf{x}\right) ,p_{1}\left( \mathbf{x}%
\right) ,...,p_{i+1}\left( \mathbf{x}\right) \right) = \\ 
=2L_{i}\left( v_{0,1}\left( \mathbf{x}\right) ,v_{1,1}\left( \mathbf{x}%
\right) ,...,v_{i+1,1}\left( \mathbf{x}\right) \right) \times \\ 
\times \left[ 
\begin{array}{c}
\Delta p_{i}\left( \mathbf{x}\right) +2\nabla p_{i}\left( \mathbf{x}\right)
\nabla w\left( \mathbf{x},\varepsilon \right) + \\ 
+2h\nabla p_{i}\left( \mathbf{x}\right) \displaystyle\sum
\limits_{j=0}^{i}\nabla v_{j,1}\left( \mathbf{x}\right) +2h\nabla
v_{i,1}\left( \mathbf{x}\right) \displaystyle\sum \limits_{j=0}^{i}\nabla
p_{j}\left( \mathbf{x}\right) + \\ 
+\left( p_{i}\left( \mathbf{x}\right) -p_{i+1}\left( \mathbf{x}\right)
\right) /h%
\end{array}%
\right] , \\ 
\mathbf{x}\in \Omega ,\text{ }i=0,...,k-1.%
\end{array}%
\right.  \label{5.8}
\end{equation}

The term $M_{i,\text{nonlin}}$ is:%
\begin{equation}
\left. 
\begin{array}{c}
M_{i,\text{nonlin}}= \\ 
=\left[ 
\begin{array}{c}
\Delta p_{i}\left( \mathbf{x}\right) +2\nabla p_{i}\left( \mathbf{x}\right)
\nabla w\left( \mathbf{x},\varepsilon \right) + \\ 
+2h\nabla p_{i}\left( \mathbf{x}\right) \sum\limits_{j=0}^{i}\nabla
v_{j,1}\left( \mathbf{x}\right) +2h\nabla v_{i,1}\left( \mathbf{x}\right)
\sum\limits_{j=0}^{i}\nabla p_{j}\left( \mathbf{x}\right) + \\ 
+\left( p_{i}\left( \mathbf{x}\right) -p_{i+1}\left( \mathbf{x}\right)
\right) /h%
\end{array}%
\right] ^{2}+ \\ 
+2h\nabla p_{i}\left( \mathbf{x}\right) \sum\limits_{j=0}^{i}\nabla
p_{j}\left( \mathbf{x}\right) + \\ 
+4L_{i}\left( v_{0,1}\left( \mathbf{x}\right) ,v_{1,1}\left( \mathbf{x}%
\right) ,...,v_{i+1,1}\left( \mathbf{x}\right) \right) \times \\ 
\times h\nabla p_{i}\left( \mathbf{x}\right) \sum\limits_{j=0}^{i}\nabla
p_{j}\left( \mathbf{x},t\right) ,\text{ } \\ 
\mathbf{x}\in \Omega ,\text{ }i=0,...,k-1.%
\end{array}%
\right.  \label{5.9}
\end{equation}

In the case $i=0,$ we set the sums in formulas (\ref{5.6})-(\ref{5.9}) equal
to zero. Next, using (\ref{3.30}), we similarly obtain the following analog
of (\ref{5.8}) for the case $i=k$%
\begin{equation*}
\left. 
\begin{array}{c}
M_{k,\text{lin}}\left( p_{0}\left( \mathbf{x}\right) ,p_{1}\left( \mathbf{x}%
\right) ,...,p_{k}\left( \mathbf{x}\right) \right) = \\ 
=2L_{k}\left( v_{0,1}\left( \mathbf{x}\right) ,v_{1,1}\left( \mathbf{x}%
\right) ,...,v_{k,1}\left( \mathbf{x}\right) \right) \times \\ 
\times \left[ 
\begin{array}{c}
\Delta p_{k}\left( \mathbf{x}\right) +2\nabla p_{k}\left( \mathbf{x}\right)
\nabla w\left( \mathbf{x},\varepsilon \right) + \\ 
+2h\nabla p_{k}\left( \mathbf{x}\right) \displaystyle\sum\limits_{j=0}^{k}%
\nabla v_{j,1}\left( \mathbf{x}\right) +2h\nabla v_{k,1}\left( \mathbf{x}%
\right) \displaystyle\sum\limits_{j=0}^{k}\nabla p_{j}\left( \mathbf{x}%
\right) + \\ 
+\left( p_{k-1}\left( \mathbf{x}\right) +p_{k-2}\left( \mathbf{x}\right)
+p_{k-3}\left( \mathbf{x}\right) -3p_{k}\left( \mathbf{x}\right) \right)
/\left( 6h\right)%
\end{array}%
\right] , \\ 
\mathbf{x}\in \Omega .%
\end{array}%
\right.
\end{equation*}%
Finally, using (\ref{3.30}) again, we obtain the following analog of (\ref%
{5.9}) for the case $i=k$%
\begin{equation}
\left. 
\begin{array}{c}
M_{k,\text{nonlin}}= \\ 
=\left[ 
\begin{array}{c}
\Delta p_{k}\left( \mathbf{x}\right) +2\nabla p_{k}\left( \mathbf{x}\right)
\nabla w\left( \mathbf{x},\varepsilon \right) + \\ 
+2h\nabla p_{k}\left( \mathbf{x}\right) \sum\limits_{j=0}^{k}\nabla
v_{j,1}\left( \mathbf{x}\right) +2h\nabla v_{k,1}\left( \mathbf{x}\right)
\sum\limits_{j=0}^{k}\nabla p_{j}\left( \mathbf{x}\right) + \\ 
+\left( p_{k-1}\left( \mathbf{x},t_{k-1}\right) +p_{k-2}\left( \mathbf{x}%
,t_{k-2}\right) +p_{k-3}\left( \mathbf{x},t_{k-3}\right) -3p_{k}\left( 
\mathbf{x}\right) \right) /\left( 6h\right)%
\end{array}%
\right] ^{2}+ \\ 
+2h\nabla p_{k}\left( \mathbf{x}\right) \sum\limits_{j=0}^{k}\nabla
p_{j}\left( \mathbf{x}\right) + \\ 
+4L_{k}\left( v_{0,1}\left( \mathbf{x}\right) ,v_{1,1}\left( \mathbf{x}%
\right) ,...,v_{i+1,1}\left( \mathbf{x}\right) \right) \times \\ 
\times h\nabla p_{k}\left( \mathbf{x}\right) \sum\limits_{j=0}^{k}\nabla
p_{j}\left( \mathbf{x},t\right) ,\text{ } \\ 
\mathbf{x}\in \Omega .%
\end{array}%
\right.  \label{5.91}
\end{equation}

It follows from (\ref{4.6}), (\ref{5.4}) and (\ref{5.7})-(\ref{5.91}) that%
\begin{equation*}
J_{\lambda ,\alpha }\left( V_{1}+P\right) -J_{\lambda ,\alpha }\left(
V_{1}\right) =
\end{equation*}%
\begin{equation*}
=e^{-2\lambda c}\sum\limits_{i=0}^{k-1}\int\limits_{\Omega }M_{i,\text{lin}%
}\left( p_{0}\left( \mathbf{x}\right) ,p_{1}\left( \mathbf{x}\right)
,...,p_{i+1}\left( \mathbf{x}\right) \right) \varphi _{\lambda }\left(
x_{1}\right) d\mathbf{x+}
\end{equation*}%
\begin{equation}
+e^{-2\lambda c}\int\limits_{\Omega }M_{k,\text{lin}}\left( p_{0}\left( 
\mathbf{x}\right) ,p_{1}\left( \mathbf{x}\right) ,...,p_{k}\left( \mathbf{x}%
\right) \right) \varphi _{\lambda }\left( x_{1}\right) d\mathbf{x}+2\alpha
\left( P,V_{1}\right) +  \label{5.11}
\end{equation}%
\begin{equation*}
+e^{-2\lambda c}\sum\limits_{i=0}^{k-1}\int\limits_{\Omega }M_{i,\text{
nonlin}}\left( p_{0}\left( \mathbf{x}\right) ,p_{1}\left( \mathbf{x}\right)
,...,p_{i+1}\left( \mathbf{x}\right) \right) \varphi _{\lambda }\left(
x_{1}\right) d\mathbf{x+}
\end{equation*}%
\begin{equation*}
+e^{-2\lambda c}\int\limits_{\Omega }M_{k,\text{nonlin}}\left( p_{0}\left( 
\mathbf{x}\right) ,p_{1}\left( \mathbf{x}\right) ,...,p_{k}\left( \mathbf{x}%
\right) \right) \varphi _{\lambda }\left( x_{1}\right) d\mathbf{x}+\alpha
\left\Vert P\right\Vert _{H_{k+1}^{m_{n}}\left( \Omega \right) }^{2}.
\end{equation*}%
Let $Q\left( \mathbf{x}\right) =\left( q_{0}\left( \mathbf{x}\right)
,q_{1}\left( \mathbf{\ x}\right) ,...,q_{k}\left( \mathbf{x}\right) \right)
^{T}$ be an arbitrary vector function such that $Q\in
H_{0,k+1}^{m_{n}}\left( \Omega \right) $. Consider the expression $K\left(
Q\right) $ in the second and third lines of (\ref{5.11}), in which the
vector function $P\left( \mathbf{x}\right) =\left( p_{0}\left( \mathbf{x}%
\right) ,p_{1}\left( \mathbf{x}\right) ,...,p_{k}\left( \mathbf{x}\right)
\right) ^{T}\in H_{0,k+1}^{m_{n}}\left( \Omega \right) $ is replaced with $Q$%
. Then 
\begin{equation*}
K\left( Q\right) =e^{-2\lambda c}\displaystyle\sum\limits_{i=0}^{k-1}%
\displaystyle\int\limits_{\Omega }M_{i,\text{lin}}\left( q_{0}\left( \mathbf{%
x}\right) ,q_{1}\left( \mathbf{x}\right) ,...,q_{i+1}\left( \mathbf{x}%
\right) \right) \varphi _{\lambda }\left( x_{1}\right) d\mathbf{x+}
\end{equation*}%
\begin{equation*}
+e^{-2\lambda c}\displaystyle\int\limits_{\Omega }M_{k,\text{lin}}\left(
q_{0}\left( \mathbf{x}\right) ,q_{1}\left( \mathbf{x}\right)
,...,q_{k}\left( \mathbf{x}\right) \right) \varphi _{\lambda }\left(
x_{1}\right) d\mathbf{x}+2\alpha \left( Q,V_{1}\right) .
\end{equation*}%
Clearly $K\left( Q\right) :H_{0,k+1}^{m_{n}}\left( \Omega \right)
\rightarrow \mathbb{R}$ is a bounded linear functional. Therefore, by Riesz
theorem there exists a vector function $\widetilde{K}\left( \mathbf{x}%
\right) \in H_{0,k+1}^{m_{n}}\left( \Omega \right) $ such that 
\begin{equation}
\left( \widetilde{K}\left( \mathbf{x}\right) ,Q\right) =K\left( Q\right) ,%
\text{ }\forall Q\in H_{0,k+1}^{m_{n}}\left( \Omega \right) .  \label{5.13}
\end{equation}%
Furthermore, it is clear from (\ref{5.4}) and (\ref{5.9})-(\ref{5.13}) that 
\begin{equation*}
\lim_{\left\Vert P\right\Vert _{H_{k+1}^{m_{n}}\left( \Omega \right)
}\rightarrow 0}\frac{\left\vert J_{\lambda ,\alpha }\left( V_{1}+P\right)
-J_{\lambda ,\alpha }\left( V_{1}\right) -\left( \widetilde{K}\left( \mathbf{%
\ x}\right) ,P\right) \right\vert }{\left\Vert P\right\Vert
_{H_{0,k+1}^{m_{n}}\left( \Omega \right) }}=0.
\end{equation*}%
Hence, $\widetilde{K}\left( \mathbf{x}\right) \in H_{0,k+1}^{m_{n}}\left(
\Omega \right) $ is the Fr\'{e}chet derivative of the functional $J_{\lambda
,\alpha }\left( V\right) $ at the point $V_{1}\in \overline{B\left( R\right) 
},$%
\begin{equation}
\widetilde{K}\left( \mathbf{x}\right) =J_{\lambda ,\alpha }^{\prime }\left(
V_{1}\right) \in H_{0,k+1}^{m_{n}}\left( \Omega \right) .  \label{5.14}
\end{equation}%
The proof of the Lipschitz continuity property (\ref{5.1}) is omitted here
since it is similar with the one in Theorem 3.1 of \cite{Bak} and also in
Theorem 5.3.1 of \cite{KL}.

Thus, (\ref{5.11}) and (\ref{5.14}) imply%
\begin{equation*}
J_{\lambda ,\alpha }\left( V_{1}+P\right) -J_{\lambda ,\alpha }\left(
V_{1}\right) -\left( J_{\lambda ,\alpha }^{\prime }\left( V_{1}\right)
,P\right) =
\end{equation*}%
\begin{equation}
=e^{-2\lambda c}\sum\limits_{i=0}^{k-1}\int\limits_{\Omega }M_{i,\text{
nonlin}}\left( p_{0}\left( \mathbf{x}\right) ,p_{1}\left( \mathbf{x}\right)
,...,p_{i+1}\left( \mathbf{x}\right) \right) \varphi _{\lambda }\left(
x_{1}\right) d\mathbf{x+}  \label{5.15}
\end{equation}%
\begin{equation*}
+e^{-2\lambda c}\int\limits_{\Omega }M_{k,\text{nonlin}}\left( p_{0}\left( 
\mathbf{x}\right) ,p_{1}\left( \mathbf{x}\right) ,...,p_{k}\left( \mathbf{x}%
\right) \right) \varphi _{\lambda }\left( x_{1}\right) d\mathbf{x}+\alpha
\left\Vert P\right\Vert _{H_{k+1}^{m_{n}}\left( \Omega \right) }^{2}.
\end{equation*}

To prove the strong convexity estimate (\ref{5.2}), we need to estimate the
right hand side of (\ref{5.15}) from the below. Denote the right hand side
of (\ref{5.15}) as $RHS$. Using (\ref{5.4}), (\ref{5.5}), (\ref{5.9}), (\ref%
{5.91}) and (\ref{5.15}), we obtain%
\begin{equation}
RHS\geq \frac{e^{-2\lambda c}}{2}\displaystyle\int\limits_{\Omega }\left(
\Delta P\right) ^{2}\varphi _{\lambda }\left( x_{1}\right) d\mathbf{x}%
-C_{1}e^{-2\lambda c}\displaystyle\int\limits_{\Omega }\left[ \left\vert
\nabla P\right\vert ^{2}+P^{2}\right] \varphi _{\lambda }\left( x_{1}\right)
d\mathbf{x+}  \label{5.16}
\end{equation}%
\begin{equation*}
+\alpha \left\Vert P\right\Vert _{H_{k+1}^{m_{n}}}^{2}.
\end{equation*}%
Applying Carleman estimate (\ref{4.40}) to the term with $\left( \Delta
P\right) ^{2}$ in (\ref{5.16}), we obtain 
\begin{equation*}
RHS\geq \frac{C}{\lambda }e^{-2\lambda c}\displaystyle\int\limits_{\Omega
}\left( \displaystyle\sum\limits_{i,j=1}^{n}P_{x_{i}x_{j}}^{2}\right)
\varphi _{\lambda }\left( x_{1}\right) d\mathbf{x}+
\end{equation*}%
\begin{equation}
+Ce^{-2\lambda c}\displaystyle\int\limits_{\Omega }\left( \lambda \left\vert
\nabla P\right\vert ^{2}+\lambda ^{3}P^{2}\right) \varphi _{\lambda }\left(
x_{1}\right) d\mathbf{x}-  \label{5.17}
\end{equation}%
\begin{equation*}
-C_{1}e^{2\lambda \left( A^{2}-c\right) }\displaystyle\int\limits_{\Omega }%
\left[ \left\vert \nabla P\right\vert ^{2}+P^{2}\right] \varphi _{\lambda
}\left( x_{1}\right) d\mathbf{x}+\alpha \left\Vert P\right\Vert
_{H_{k+1}^{m_{n}}}^{2}.
\end{equation*}%
Choose a sufficiently large number $\lambda _{1}=\lambda _{1}\left( R,\alpha
,\varepsilon ,h_{0},T\right) \geq \lambda _{0}$ such that $C$ $\lambda
_{1}\geq 2C_{1}.$ Hence, using (\ref{5.15}) and (\ref{5.17}) and keeping in
mind that by (\ref{2.1}) and (\ref{4.4}) $\varphi _{\lambda }\left(
x_{1}\right) \geq e^{2\lambda A^{2}}$ in $\Omega ,$ we obtain 
\begin{equation*}
J_{\lambda ,\alpha }\left( V_{1}+P\right) -J_{\lambda ,\alpha }\left(
V_{1}\right) -\left( J_{\lambda ,\alpha }^{\prime }\left( V_{1}\right)
,P\right) \geq
\end{equation*}%
\begin{equation*}
\geq \frac{C}{\lambda }e^{2\lambda \left( A^{2}-c\right)
}\int\limits_{\Omega }\left(
\sum\limits_{i,j=1}^{n}P_{x_{i}x_{j}}^{2}\right) d\mathbf{x}%
+C_{1}e^{2\lambda \left( A^{2}-c\right) }\int\limits_{\Omega }\left(
\left\vert \nabla P\right\vert ^{2}+P^{2}\right) d\mathbf{x}+
\end{equation*}%
\begin{equation*}
+\alpha \left\Vert P\right\Vert _{H_{k+1}^{m_{n}}\left( \Omega \right) }^{2},%
\text{ }\forall \lambda \geq \lambda _{1},
\end{equation*}%
which proves (\ref{5.2}).

The existence and uniqueness of the minimizer $V_{\min ,\lambda }\in 
\overline{B\left( R\right) }$ of the functional $J_{\lambda ,\alpha }\left(
V\right) $ on the set $\overline{B\left( R\right) }$ as well as inequality (%
\ref{5.3}) follow immediately from a combination of either Lemma 2.1 with
Theorem 2.1 of \cite{Bak} or Lemma 5.2.1 with Theorem 5.2.1 of \cite{KL}. $%
\square $
\end{proof}

\subsection{Accuracy estimates in the case of noisy data}

\label{sec:6.2}

It is always assumed in the theory of Ill-Posed problems that there exists a
true solution of a CIP at hands for the case of the \textquotedblleft
ideal", i.e. noiseless data \cite{T}.

Thus, let the function $a^{\ast }\left( \mathbf{x}\right) $ be the true
solution of the finite difference version of our CIP. In other words, we
assume that $a^{\ast }\left( \mathbf{x}\right) $ generates the following
analog of the vector function $V\left( \mathbf{x}\right) $ in (\ref{3.35})%
\begin{equation}
\left. V^{\ast }\left( \mathbf{x}\right) =\left( v^{\ast }\left( \mathbf{x}%
,t_{0}\right) ,v^{\ast }\left( \mathbf{x},t_{1}\right) ,...,v^{\ast }\left( 
\mathbf{x},t_{k}\right) \right) ^{T},\text{ }\mathbf{x}\in \Omega .\right. 
\label{5.22}
\end{equation}%
More precisely, we assume that functions $v^{\ast }\left( \mathbf{x}%
,t_{i}\right) $ satisfy equations (\ref{3.25})-(\ref{3.30}) i.e. we assume
that%
\begin{equation}
\left. 
\begin{array}{c}
L_{0}\left( v^{\ast }\left( \mathbf{x},t_{0}\right) ,v^{\ast }\left( \mathbf{%
\ x},t_{1}\right) \right) = \\ 
=\Delta v^{\ast }\left( \mathbf{x},t_{0}\right) +2\nabla v^{\ast }\left( 
\mathbf{x},t_{0}\right) \nabla w\left( \mathbf{x},\varepsilon \right) + \\ 
+\left( v^{\ast }\left( \mathbf{x},t_{0}\right) -v^{\ast }\left( \mathbf{x}%
,t_{1}\right) \right) /h,\text{ }\mathbf{x}\in \Omega ,%
\end{array}%
\right.   \label{5.23}
\end{equation}%
\begin{equation*}
\left. 
\begin{array}{c}
L_{i}\left( v^{\ast }\left( \mathbf{x},t_{0}\right) ,v^{\ast }\left( \mathbf{%
\ x},t_{1}\right) ,...,v^{\ast }\left( \mathbf{x},t_{i+1}\right) \right) =
\\ 
=\Delta v^{\ast }\left( \mathbf{x},t_{i}\right) +2\nabla v^{\ast }\left( 
\mathbf{x},t_{i}\right) \nabla w\left( \mathbf{x},\varepsilon \right)
+2h\nabla v^{\ast }\left( \mathbf{x},t_{i}\right)
\sum\limits_{j=0}^{i}\nabla v^{\ast }\left( \mathbf{x},t_{j}\right) + \\ 
+\left( v^{\ast }\left( \mathbf{x},t_{i}\right) -v^{\ast }\left( \mathbf{x}%
,t_{i+1}\right) \right) /h=0,\text{ } \\ 
\mathbf{x}\in \Omega ,\text{ }i=1,...,k-1,%
\end{array}%
\right. 
\end{equation*}%
\begin{equation}
\left. 
\begin{array}{c}
L_{k}\left( v^{\ast }\left( \mathbf{x},t_{0}\right) ,v^{\ast }\left( \mathbf{%
\ x},t_{1}\right) ,v^{\ast }\left( \mathbf{x},t_{2}\right) ,...,v^{\ast
}\left( \mathbf{x},t_{k-1}\right) ,v^{\ast }\left( \mathbf{x},t_{k}\right)
\right) = \\ 
=\Delta v^{\ast }\left( \mathbf{x},t_{k}\right) +2\nabla v^{\ast }\left( 
\mathbf{x},t_{k}\right) \nabla w\left( \mathbf{x},\varepsilon \right) + \\ 
+2h\nabla v^{\ast }\left( \mathbf{x},t_{k}\right) \left(
\sum\limits_{j=0}^{k}\nabla v^{\ast }\left( \mathbf{x},t_{j}\right) \right) +
\\ 
+\left( v^{\ast }\left( \mathbf{x},t_{k-1}\right) +v^{\ast }\left( \mathbf{x}%
,t_{k-2}\right) +v^{\ast }\left( \mathbf{x},t_{k-3}\right) -3v^{\ast }\left( 
\mathbf{x},t_{k}\right) \right) /\left( 6h\right) =0, \\ 
\mathbf{x}\in \Omega .%
\end{array}%
\right.   \label{5.26}
\end{equation}%
In addition, we assume that functions $v^{\ast }\left( \mathbf{x}%
,t_{i}\right) $ satisfy boundary conditions (\ref{3.31})-(\ref{3.34}) with
noiseless data $g_{0}^{\ast },g_{1}^{\ast },$ 
\begin{equation}
v^{\ast }\left( \mathbf{x},t_{i}\right) =\frac{g_{0}^{\ast }\left( \mathbf{x}%
,t_{i+1}\right) -g_{0}^{\ast }\left( \mathbf{x},t_{i}\right) }{h},\text{ }%
\mathbf{x}\in \partial \Omega ,\text{ }i=0,...,k-1,  \label{5.27}
\end{equation}%
\begin{equation*}
v^{\ast }\left( \mathbf{x},t_{k}\right) =\frac{3g_{0}^{\ast }\left( \mathbf{x%
},t_{k}\right) -g_{0}^{\ast }\left( \mathbf{x},t_{k-1}\right) -g_{0}^{\ast
}\left( \mathbf{x},t_{k-2}\right) -g_{0}^{\ast }\left( \mathbf{x}%
,t_{k-3}\right) }{6h},\text{ }\mathbf{x}\in \partial \Omega ,
\end{equation*}%
\begin{equation*}
v_{x_{1}}^{\ast }\left( \mathbf{x},t_{i}\right) =\frac{g_{1}^{\ast }\left( 
\mathbf{x},t_{i+1}\right) -g_{1}^{\ast }\left( \mathbf{x},t_{i}\right) }{h},%
\text{ }\mathbf{x}\in \Gamma _{0},\text{ }i=0,...,k-1,
\end{equation*}%
\begin{equation}
v_{x_{1}}^{\ast }\left( \mathbf{x},t_{k}\right) =\frac{3g_{1}^{\ast }\left( 
\mathbf{x},t_{k}\right) -g_{1}^{\ast }\left( \mathbf{x},t_{k-1}\right)
-g_{1}^{\ast }\left( \mathbf{x},t_{k-2}\right) -g_{1}^{\ast }\left( \mathbf{x%
},t_{k-3}\right) }{6h},\text{ }\mathbf{x}\in \Gamma _{0}.  \label{5.30}
\end{equation}%
We define the set $B^{\ast }\left( R\right) $ analogously with the set $%
B^{\ast }\left( R\right) $ in (\ref{4.2}), 
\begin{equation}
B^{\ast }\left( R\right) =\left\{ 
\begin{array}{c}
W\left( \mathbf{x}\right) =\left( w_{0}\left( \mathbf{x}\right) ,w_{1}\left( 
\mathbf{x}\right) ,...,w_{k}\left( \mathbf{x}\right) \right) ^{T}\in
H_{k+1}^{m_{n}}\left( \Omega \right) : \\ 
\left\Vert w_{k}\left( \mathbf{x}\right) \right\Vert _{H^{m_{n}}\left(
\Omega \right) }<R,\text{ }i=0,...,k, \\ 
\text{functions }w_{i}\left( \mathbf{x}\right) \text{ satisfy boundary
conditions } \\ 
\text{(\ref{5.27})-(\ref{5.30}) for corresponding indexes }i%
\end{array}%
\right\} .  \label{5.31}
\end{equation}%
Hence, using (\ref{5.22}), we assume that 
\begin{equation}
V^{\ast }\left( \mathbf{x}\right) \in B^{\ast }\left( R\right) .
\label{5.32}
\end{equation}%
Similarly with (\ref{3.36}) let 
\begin{equation}
w^{\ast }\left( \mathbf{x},t_{i}\right) =h\displaystyle\sum%
\limits_{j=0}^{i}v^{\ast }\left( \mathbf{x},t_{j}\right) +w\left( \mathbf{x}%
,\varepsilon \right) ,  \label{5.33}
\end{equation}%
where $w\left( \mathbf{x},\varepsilon \right) $ is given in (\ref{3.8}).
Thus, similarly with (\ref{3.350})-(\ref{3.37}) we assume that 
\begin{equation*}
a^{\ast }\left( \mathbf{x}\right) =\frac{1}{T-\varepsilon }\left( w^{\ast
}\left( \mathbf{x},t_{k}\right) -w\left( \mathbf{x},\varepsilon \right)
\right) -
\end{equation*}%
\begin{equation}
-\frac{h}{T-\varepsilon }\sum\limits_{i=0}^{k}\left( \Delta w^{\ast }\left( 
\mathbf{x},t_{i}\right) +\left( \nabla w^{\ast }\left( \mathbf{x}%
,t_{i}\right) \right) ^{2}\right) ,  \label{5.34}
\end{equation}

The minimizer $V_{\min ,\lambda }\left( \mathbf{x}\right) ,$ which was found
in Theorem 6.1, is called \textquotedblleft regularized solution" in the
theory of Ill-Posed Problems \cite{T}. An estimate of the distance between
the regularized and true solutions of a CIP naturally represents an
important task. We provide such an estimate in this subsection for an analog
of $V_{\min ,\lambda }\left( \mathbf{x}\right) $. Our estimate involves two
small parameters: the level of the noise $\sigma \in \left( 0,1\right) $ in
the input data and the regularization parameter $\alpha \in \left(
0,1\right) $ in (\ref{4.6}). In parallel we estimate the accuracy of the
reconstruction of our target coefficient $a^{\ast }\left( \mathbf{x}\right)
. $

We assume the existence of such a pair of vector functions $G\left( \mathbf{x%
}\right) $ and $G^{\ast }\left( \mathbf{x}\right) $ that%
\begin{equation}
\left. 
\begin{array}{c}
G\in B\left( R\right) ,\text{ }G^{\ast }\in B^{\ast }\left( R\right) , \\ 
\left\Vert G-G^{\ast }\right\Vert _{H_{k+1}^{m_{n}}}<\sigma ,%
\end{array}%
\right.  \label{5.38}
\end{equation}%
where a sufficiently small number $\sigma \in \left( 0,1\right) $
characterizes the level of the noise in the boundary data (\ref{3.31})-(\ref%
{3.34}) and sets $B\left( R\right) $ and $B^{\ast }\left( R\right) $ were
defined in (\ref{4.2}) and (\ref{5.31}) respectively. Since $G^{\ast }\in
B^{\ast }\left( R\right) ,$ then (\ref{5.31}) and (\ref{5.32}) imply that 
\begin{equation}
V^{\ast }-G^{\ast }\in H_{0,k+1}^{m_{n}}\left( \Omega \right) .
\label{5.380}
\end{equation}%
For each vector function $V\in \overline{B\left( R\right) }$ consider the
vector function 
\begin{equation}
Q\left( V\right) =V-G\in \overline{B_{0}\left( 2R\right) },  \label{5.39}
\end{equation}%
where the set $\overline{B_{0}\left( 2R\right) }$ is defined in (\ref{5.5}).
By (\ref{5.38})-(\ref{5.39}) 
\begin{equation}
Q\left( V^{\ast }\right) =V^{\ast }-G^{\ast }\in \overline{B_{0}\left(
2R\right) }.  \label{5.40}
\end{equation}

Hence, for each $V\in \overline{B\left( R\right) }$ the vector function $%
Q\left( V\right) $ satisfies exactly the same boundary conditions (\ref{3.31}%
)-(\ref{3.34}) but with zeros in their right hand sides. Hence, we can apply
now an analog of Theorem 5.1. More precisely, consider the following
functional 
\begin{equation}
\left. 
\begin{array}{c}
I_{\lambda ,\alpha }:\overline{B_{0}\left( 2R\right) }\rightarrow \mathbb{R}%
\text{,} \\ 
I_{\lambda ,\alpha }\left( W\right) =J_{\lambda ,\alpha }\left( W+G\right) .%
\end{array}%
\right.  \label{5.41}
\end{equation}%
An obvious analog of Theorem 6.1 is valid for $I_{\lambda ,\alpha }$. We are
not formulating this analog here for brevity. We only note that the number $%
\lambda _{1}=\lambda _{1}\left( R,\alpha ,\varepsilon ,h_{0},T\right) $ of
Theorem 5.1 should obviously be replaced with the number $\lambda _{2},$%
\begin{equation}
\lambda _{2}=\lambda _{1}\left( 2R,\alpha ,\varepsilon ,h_{0},T\right) .
\label{5.42}
\end{equation}

\textbf{Theorem 6.2} (accuracy estimates for noisy data). \emph{Assume that
conditions (\ref{5.22})-(\ref{5.38}) are satisfied. Let }$I_{\lambda ,\alpha
}$\emph{\ be the functional defined in (\ref{5.41}) and }$\lambda _{2}$\emph{%
\ be the number in (\ref{5.42}). For any }$\lambda \geq \lambda _{2}$\emph{\
let }$W_{\min ,\lambda }\in \overline{B_{0}\left( 2R\right) }$\emph{\ be the
unique minimizer of the functional }$I_{\lambda ,\alpha }$\emph{\ on the set 
}$\overline{B_{0}\left( 2R\right) },$\emph{\ the existence and uniqueness of
which is guaranteed by the above mentioned analog of Theorem 6.1. Let }$a_{%
\text{comp,min,}\lambda }\left( \mathbf{x}\right) $\emph{\ be the function
in the left hand side of equality (\ref{3.37}), the right hand side of which
is formed as in (\ref{3.350}), (\ref{3.36}), where the vector function }$V_{%
\text{comp}}\left( \mathbf{x}\right) $ \emph{is replaced with }$W_{\min
,\lambda }+G$\emph{\ with corresponding replacements of all components of} $%
V_{\text{comp}}\left( \mathbf{x}\right) $. \emph{Let }$a^{\ast }\left( 
\mathbf{x}\right) $\emph{\ be the function constructed in (\ref{5.33}), (\ref%
{5.34}). Then the following accuracy estimates hold for all }$\lambda \geq
\lambda _{2}$%
\begin{equation}
\left\Vert W_{\min ,\lambda }-Q\left( V^{\ast }\right) \right\Vert
_{H_{k+1}^{2}\left( \Omega \right) }\leq C_{1}\sqrt{\lambda }\sqrt{\sigma }%
\cdot e^{\lambda \left( B^{2}-A^{2}\right) }+C_{1}\sqrt{\lambda \alpha }%
\cdot e^{-\left( A^{2}-c\right) },  \label{5.43}
\end{equation}%
\begin{equation}
\left\Vert a_{\text{comp,min,}\lambda }-a^{\ast }\right\Vert _{L_{2}\left(
\Omega \right) }\leq C_{1}\sqrt{\lambda }\sqrt{\sigma }\cdot e^{\lambda
\left( B^{2}-A^{2}\right) }+C_{1}\sqrt{\lambda \alpha }e^{-\left(
A^{2}-c\right) }.  \label{5.44}
\end{equation}

\begin{proof}
It follows from (\ref{5.380})-(\ref{5.40}) that $Q\left( V^{\ast }\right)
,W_{\min ,\lambda }\in \overline{B_{0}\left( 2R\right) }.$ Hence, applying
the above mentioned analog of Theorem 6.1 to the functional $I_{\lambda
,\alpha }$ in (\ref{5.41}) and ignoring the non-negative term $\alpha
\left\Vert Q\left( V^{\ast }\right) -W_{\min ,\lambda }\right\Vert
_{H_{k+1}^{m_{n}}\left( \Omega \right) }^{2},$ we obtain%
\begin{equation*}
\lambda I_{\lambda ,\alpha }\left( Q\left( V^{\ast }\right) \right) -\lambda
I_{\lambda ,\alpha }\left( W_{\min ,\lambda }\right) -\lambda \left[
I_{\lambda ,\alpha }^{\prime }\left( W_{\min ,\lambda }\right) ,Q\left(
V^{\ast }\right) -W_{\min ,\lambda }\right] \geq
\end{equation*}%
\begin{equation}
\geq C_{1}e^{2\lambda \left( A^{2}-c\right) }\left\Vert Q\left( V^{\ast
}\right) -W_{\min ,\lambda }\right\Vert _{H_{k+1}^{2}\left( \Omega \right)
}^{2},\text{ }\forall \lambda \geq \lambda _{2}.  \label{5.45}
\end{equation}%
Using (\ref{5.3}), we obtain%
\begin{equation*}
-\lambda \left[ I_{\lambda ,\alpha }^{\prime }\left( W_{\min ,\lambda
}\right) ,Q\left( V^{\ast }\right) -W_{\min ,\lambda }\right] \leq 0.
\end{equation*}%
Also, obviously $-\lambda I_{\lambda ,\alpha }\left( W_{\min ,\lambda
}\right) \leq 0.$ Hence, (\ref{5.45}) implies 
\begin{equation}
\lambda I_{\lambda ,\alpha }\left( Q\left( V^{\ast }\right) \right) \geq
C_{1}e^{2\lambda A^{2}}\left\Vert Q\left( V^{\ast }\right) -W_{\min ,\lambda
}\right\Vert _{H_{k+1}^{2}\left( \Omega \right) }^{2}.  \label{5.46}
\end{equation}

We now estimate the left hand side of (\ref{5.46}). By (\ref{5.40}) and (\ref%
{5.41})%
\begin{equation*}
I_{\lambda ,\alpha }\left( Q\left( V^{\ast }\right) \right) =J_{\lambda
,\alpha }\left( Q\left( V^{\ast }\right) +G\right) =J_{\lambda ,\alpha
}\left( V^{\ast }+\left( G-G^{\ast }\right) \right) .
\end{equation*}%
Hence, applying (\ref{4.6}), (\ref{4.60}), (\ref{5.38}) and Cauchy-Schwarz
inequality and also using the fact that \textquotedblleft $=0"$ is present
in each 
of equalities (\ref{5.23})-(\ref{5.26}), we obtain%
\begin{equation}
\lambda I_{\lambda ,\alpha }\left( Q\left( V^{\ast }\right) \right) =\lambda
J_{\lambda ,\alpha }\left( V^{\ast }+\left( G-G^{\ast }\right) \right) \leq
C_{1}\lambda e^{2\lambda \left( B^{2}-c\right) }\sigma ^{2}+C_{1}\lambda
\alpha .  \label{5.47}
\end{equation}%
The first target estimate (\ref{5.43}) of this theorem follows immediately
from (\ref{5.46}) and (\ref{5.47}). Finally, using (\ref{5.43}), the above
described constructions of functions $a_{\text{comp}}\left( \mathbf{x}%
\right) $ and $a^{\ast }\left( \mathbf{x}\right) $ respectively in (\ref%
{3.36})-(\ref{3.37}) and (\ref{5.33})-(\ref{5.34}) as well as (\ref{5.38})-(%
\ref{5.40}), we obtain the second target estimate (\ref{5.44}). $\square $
\end{proof}

\subsection{Accuracy estimates in the case of noiseless data}

\label{sec:6.3}

Recall that $\sigma $ is the level of the noise in the data. In the
noiseless case the vector function $V^{\ast }\left( \mathbf{x}\right) $
satisfies exactly the same boundary conditions as the ones for the minimizer 
$V_{\min ,\lambda }\left( \mathbf{x}\right) \in \overline{B\left( R\right) }%
, $ which was found in Theorem 6.1. In other words, we assume now that%
\begin{equation}
\left. 
\begin{array}{c}
\sigma =0, \\ 
g_{0}^{\ast }\left( \mathbf{x},t_{i}\right) =g_{0}\left( \mathbf{x}%
,t_{i}\right) ,\text{ }g_{1}^{\ast }\left( \mathbf{x},t_{i}\right)
=g_{1}\left( \mathbf{x},t_{i}\right) ,\text{ }i=0,...,k,%
\end{array}%
\right.  \label{5.48}
\end{equation}%
where functions $g_{0}\left( \mathbf{x},t_{i}\right) ,g_{1}\left( \mathbf{x}%
,t_{i}\right) $ and $g_{0}^{\ast }\left( \mathbf{x},t_{i}\right)
,g_{1}^{\ast }\left( \mathbf{x},t_{i}\right) $ are involved in boundary
conditions (\ref{3.31})-(\ref{3.34}) and (\ref{5.27})-(\ref{5.30})
respectively. Hence, there is no need to construct the functional $%
I_{\lambda ,\alpha }$ in (\ref{5.41}). Rather, we can replace assumption (%
\ref{5.32}) with 
\begin{equation}
V^{\ast }\left( \mathbf{x}\right) \in B\left( R\right) .  \label{5.49}
\end{equation}%
We omit the proof of Theorem 6.3 since it is completely similar with the
proof of Theorem 6.2 in the case $\sigma =0.$

\textbf{Theorem 6.3.} \emph{Assume that conditions (\ref{5.22})-(\ref{5.30}%
), (\ref{5.48}) and (\ref{5.49}) hold. Let }$\lambda _{1}=\lambda _{1}\left(
R,\alpha ,\varepsilon ,h_{0},T\right) \geq \lambda _{0}\geq 1$ \emph{\ be
the number of that theorem. For any }$\lambda \geq \lambda _{1}$\emph{\ let }%
$V_{\min ,\lambda }\left( \mathbf{x}\right) \in \overline{B\left( R\right) }$%
\emph{\ be the unique minimizer of the functional }$J_{\lambda ,\alpha }$%
\emph{\ on the set }$\overline{B\left( R\right) },$\emph{\ which was found
in Theorem 5.1. Let }$a_{\text{comp,min,}\lambda }\left( \mathbf{x}\right) $%
\emph{\ be the function in the left hand side of equality (\ref{3.37}), the
right hand side of which is formed as in (\ref{3.350}), (\ref{3.36}), where
the vector function }$V_{\text{comp}}\left( \mathbf{x}\right) $ \emph{is
replaced with }$V_{\min ,\lambda }\left( \mathbf{x}\right) $\emph{\ \ with
corresponding replacements of all components of} $V_{\text{comp}}\left( 
\mathbf{x}\right) $. \emph{Let }$a^{\ast }\left( \mathbf{x}\right) $\emph{\
be the function constructed in (\ref{5.33}), (\ref{5.34}). Then the
following accuracy estimates hold for all }$\lambda \geq \lambda _{1}$%
\begin{equation}
\left\Vert V_{\min ,\lambda }-V^{\ast }\right\Vert _{H_{k+1}^{2}\left(
\Omega \right) }\leq C_{1}\sqrt{\lambda \alpha }\cdot e^{-\lambda \left(
A^{2}-c\right) },  \label{5.51}
\end{equation}%
\begin{equation}
\left\Vert a_{\text{comp,min,}\lambda }-a^{\ast }\right\Vert _{L_{2}\left(
\Omega \right) }\leq C_{1}\sqrt{\lambda \alpha }\cdot e^{-\lambda \left(
A^{2}-c\right) }.  \label{5.52}
\end{equation}

\subsection{Global convergence of the gradient descent method}

\label{sec:6.4}

To simplify the presentation, we consider in this section only the case of
the noiseless data with the noise level $\sigma =0$ as in the first line of (%
\ref{5.48}). The case of noisy data can be handled along the same lines,
see, e.g. \cite[Theorem 4.5]{MFG1}. Assume that in (\ref{4.6})%
\begin{equation*}
c\in \left( 0,A^{2}\right) .
\end{equation*}%
Let $\lambda \geq \lambda _{1}$ is so large and the regularization parameter 
$\alpha $ is so small that 
\begin{equation}
\frac{R}{3}>C_{1}\sqrt{\lambda \alpha }\cdot e^{-\lambda \left(
A^{2}-c\right) },  \label{5.53}
\end{equation}%
Assume that 
\begin{equation}
V^{\ast }\in B\left( \frac{R}{3}-C_{1}\sqrt{\lambda \alpha }\cdot
e^{-\lambda \left( A^{2}-c\right) }\right) .  \label{5.54}
\end{equation}%
Consider an arbitrary vector function $V_{0}$ such that 
\begin{equation}
V_{0}\in B\left( \frac{R}{3}-C_{1}\sqrt{\lambda \alpha }\cdot e^{-\lambda
\left( A^{2}-c\right) }\right) .  \label{5.55}
\end{equation}%
Let $\gamma >0$ be a number and let 
\begin{equation}
J_{\lambda ,\alpha }^{\prime }\left( V\right) \in H_{0,k+1}^{m_{n}}\left(
\Omega \right) ,\text{ }\forall V\in \overline{B\left( R\right) }
\label{5.56}
\end{equation}%
be the Fr\'{e}chet derivative of the functional $J_{\lambda ,\alpha }\left(
V\right) ,$ which was found in Theorem 5.1. The sequence of the gradient
descent method is 
\begin{equation}
V_{n}=V_{n-1}-\gamma J_{\lambda ,\alpha }^{\prime }\left( V_{n-1}\right) ,%
\text{ }n=1,2,...  \label{5.57}
\end{equation}%
Note that since by (\ref{5.56}) $J_{\lambda ,\alpha }^{\prime }\left(
V_{n-1}\right) \in H_{0,k+1}^{m_{n}}\left( \Omega \right) $ and since (\ref%
{5.55}) holds, then all vector functions $V_{n}$ have the same boundary
conditions (\ref{3.31})-(\ref{3.34}), which are the same as ones in (\ref%
{5.27})-(\ref{5.30}).

\textbf{Theorem 6.4} (global convergence of the gradient descent method (\ref%
{5.57})). \emph{Let }$\lambda \geq \lambda _{1},$\emph{\ where the number }$%
\lambda _{1}$\emph{\ was chosen in Theorem 5.1. Assume that conditions of
Theorem 6.3 hold, so as conditions (\ref{5.53})-(\ref{5.55}). Then }%
\begin{equation}
V_{\min ,\lambda }\in B\left( \frac{R}{3}\right) .  \label{5.58}
\end{equation}%
\emph{Furthermore, there exists a sufficiently small number }$\gamma \in
\left( 0,1\right) $\emph{\ such that the sequence }%
\begin{equation}
\left\{ V_{n}\right\} _{n=1}^{\infty }\subset B\left( R\right) ,
\label{5.580}
\end{equation}%
\emph{and this sequence converges to }$V_{\min ,\lambda }.$\emph{\ More
precisely,} \emph{there exists a number }$\theta =\theta \left( \gamma
\right) \in \left( 0,1\right) $\emph{\ such that }%
\begin{equation}
\left\Vert V_{\min ,\lambda }-V_{n}\right\Vert _{H_{0,k+1}^{m_{n}}\left(
\Omega \right) }\leq \theta ^{n}\left\Vert V_{\min ,\lambda
}-V_{0}\right\Vert _{H_{0,k+1}^{m_{n}}\left( \Omega \right) },\text{ }%
\forall n\geq 1,  \label{5.59}
\end{equation}%
\begin{equation}
\left. 
\begin{array}{c}
\left\Vert a_{n}-a^{\ast }\right\Vert _{L_{2}\left( \Omega \right) }\leq
C_{1}\sqrt{\lambda \alpha }\cdot e^{-\lambda \left( A^{2}-c\right) }+ \\ 
+\theta ^{n}\left\Vert V_{\min ,\lambda }-V_{0}\right\Vert
_{H_{0,k+1}^{m_{n}}\left( \Omega \right) },\text{ }\forall n\geq 1,%
\end{array}%
\right.  \label{5.60}
\end{equation}%
\emph{where functions }$a_{n}$\emph{\ are constructed as in (\ref{3.350})-( %
\ref{3.37}), where the vector function }$V_{\text{comp}}\left( \mathbf{x}%
\right) $\emph{\ as well as its components are replaced with }%
\begin{equation*}
V_{n}\left( \mathbf{x}\right) =\left( v_{n}\left( \mathbf{x},t_{0}\right)
,v_{n}\left( \mathbf{x},t_{1}\right) ,...,v_{n}\left( \mathbf{x}%
,t_{k}\right) \right) ^{T},\text{ }\mathbf{x}\in \Omega
\end{equation*}%
\emph{with corresponding replacements of all components of} $V_{\text{comp}%
}\left( \mathbf{x}\right) $. \emph{\ }

\begin{proof}
Formula (\ref{5.58}) follows immediately from (\ref{5.51}), (\ref{5.54}) and
triangle inequality.\emph{\ }Relations (\ref{5.580}) and (\ref{5.59}) follow
from \cite[Theorem 6]{Klibgrad}. Finally, (\ref{5.60}) easily follows from (%
\ref{5.52}), (\ref{5.59}) and triangle inequality.
\end{proof}

\textbf{Remark 6.2.} Since a smallness condition is not imposed on the
number $R$, then Theorem 6.4 ensures the global convergence of the gradient
descent method (\ref{5.57}) in terms of Definition 1.1.

\section{Uniqueness Theorem for CIP2}

\label{sec:7}

As stated in subsection 3.2, the approximation of the function $u\left( 
\mathbf{x},\varepsilon \right) $ by formula (\ref{3.2}) combined with the
above assumption of finite differences with respect to $t$ form our
approximate mathematical model. In this section we prove uniqueness theorem
for CIP2 within the framework of that model.

\textbf{Theorem 7.1} (uniqueness). \emph{Assume that conditions (\ref{5.31}%
)-(\ref{5.33}) hold. Then there exists at most one function }$a^{\ast
}\left( \mathbf{x}\right) $\emph{\ satisfying (\ref{5.34}).}

\begin{proof}
The function $a^{\ast }\left( \mathbf{x}\right) $\emph{\ }in (\ref{5.34}) is
generated by the vector function $V^{\ast }\left( \mathbf{x}\right) $ in (%
\ref{5.22}). Therefore, it is sufficient to prove that there exists at most
one vector function $V^{\ast }\left( \mathbf{x}\right) $ satisfying
conditions (\ref{5.22})-(\ref{5.32}). Assume that there exist two such
vector functions: $V_{1}^{\ast }\left( \mathbf{x}\right) $ and $V_{2}^{\ast
}\left( \mathbf{x}\right) $. Denote%
\begin{equation*}
D\left( \mathbf{x}\right) =V_{1}^{\ast }\left( \mathbf{x}\right)
-V_{2}^{\ast }\left( \mathbf{x}\right) =\left( d_{0}\left( \mathbf{x}\right)
,...,d_{k}\left( \mathbf{x}\right) \right) ^{T}.
\end{equation*}%
Then (\ref{4.01}) and (\ref{5.27})-(\ref{5.32}) imply that%
\begin{equation}
D\in H_{0,k+1}^{m_{n}}\left( \Omega \right) .  \label{5.35}
\end{equation}%
Subtracting equations (\ref{5.23})-(\ref{5.26}) for components of the vector
function $V_{2}^{\ast }\left( \mathbf{x}\right) $ from the same equations
for components of the vector function $V_{1}^{\ast }\left( \mathbf{x}\right) 
$ and using (\ref{4.01}) and (\ref{5.35}), we obtain 
\begin{equation}
\Delta D+F_{1}\left( \mathbf{x}\right) \nabla D+F_{2}\left( \mathbf{x}%
\right) D=0,\text{ }\mathbf{x}\in \Omega ,  \label{5.36}
\end{equation}%
\begin{equation}
D\mid _{\partial \Omega }=0,\text{ }D_{x_{1}}\mid _{\partial \Omega }=0,
\label{5.37}
\end{equation}%
where $F_{1}\left( \mathbf{x}\right) $ is a $\left( k+1\right) \times \left(
n\left( k+1\right) \right) $ matrix, $F_{2}\left( \mathbf{x}\right) $ is a $%
\left( k+1\right) \times \left( k+1\right) $ matrix, and all components of
both matrices belong to $C\left( \overline{\Omega }\right) .$ Square both
sides of (\ref{5.36}). Then multiply the result by the CWF $\varphi
_{\lambda }\left( x_{1}\right) $ in (\ref{4.4}) and use Cauchy-Schwarz
inequality. We obtain%
\begin{equation}
C_{1}\int\limits_{\Omega }\left( \left\vert \nabla D\right\vert
^{2}+D^{2}\right) \varphi _{\lambda }\left( x_{1}\right) d\mathbf{x\geq }%
\int\limits_{\Omega }\left( \Delta D\right) ^{2}\varphi _{\lambda }\left(
x_{1}\right) d\mathbf{x.}  \label{5.21}
\end{equation}%
Applying Carleman estimate (\ref{4.40}) to the right hand side of (\ref{5.21}%
) and using (\ref{5.37}), we obtain%
\begin{equation*}
C_{1}\int\limits_{\Omega }\left( \left\vert \nabla D\right\vert
^{2}+D^{2}\right) \varphi _{\lambda }\left( x_{1}\right) d\mathbf{x\geq }
\end{equation*}%
\begin{equation*}
\geq \frac{1}{\lambda }\int\limits_{\Omega }\left(
\sum\limits_{i,j=1}^{n}D_{x_{i}x_{j}}^{2}\right) \varphi _{\lambda }\left(
x_{1}\right) d\mathbf{x}+\int\limits_{\Omega }\left( \lambda \left\vert
\nabla D\right\vert ^{2}+\lambda ^{3}D^{2}\right) \varphi _{\lambda }\left(
x_{1}\right) d\mathbf{x,}\text{ }\forall \lambda \geq \lambda _{0}\geq 1.
\end{equation*}

Choosing $\widetilde{\lambda }\geq \lambda _{0}$ so large that $\widetilde{%
\lambda }\geq 2C_{1},$ we obtain%
\begin{equation*}
\frac{1}{\lambda }\int\limits_{\Omega }\left(
\sum\limits_{i,j=1}^{n}D_{x_{i}x_{j}}^{2}\right) \varphi _{\lambda }\left(
x_{1}\right) d\mathbf{x}+\int\limits_{\Omega }\left( \lambda \left\vert
\nabla D\right\vert ^{2}+\lambda ^{3}D^{2}\right) \varphi _{\lambda }\left(
x_{1}\right) d\mathbf{x}\leq 0,\text{ }\forall \lambda \geq \widetilde{%
\lambda }.
\end{equation*}%
Thus, $D\left( \mathbf{x}\right) \equiv 0$ in $\Omega .$ $\square $
\end{proof}

\section{Numerical Studies}

\label{sec:8}

To demonstrate a high robustness of our numerical method, we test it
numerically for rather complicated letter-like shapes of inclusions. Indeed,
letters are non-convex and have voids. We are not concerned here with some
blur in our images. Indeed, the main thing for us is to accurately image
both shapes of targets and values of the coefficient $a\left( \mathbf{x}%
\right) $ inside of them. Blur can be reduced on a later stage via applying
one of standard \textquotedblleft image cleaning" procedures. However, the
latter is outside of the scope of this paper. We also note that even though
our theory requires a sufficient smoothness of the coefficient $a\left( 
\mathbf{x}\right) $, see (\ref{2.7}), our numerical results are sort of
\textquotedblleft less pessimistic". In other words, our numerical studies
demonstrate a quite well performance for piecewise continuous functions $%
a\left( \mathbf{x}\right) .$ Such observations often take place in numerical
studies.

First, we need to select proper values of four parameters: the value of $%
\varepsilon $ in (\ref{3.2}), the step size $h$ of the finite differences in
the $t-$direction in (\ref{3.16}), the parameter $\lambda \geq 1$ in the CWF
(\ref{4.4}) and the regularization parameter $\alpha \in \left( 0,1\right) $
in (\ref{4.6}). In principle, some insignificant, although space consuming
modifications of the above theory allow us to make these choices
analytically. Then, however, the values of those parameters would be
significantly under/over estimated. Thus, based on our rich experience of
the above cited publications about the convexification method, we do these
choices by trial and error, since exactly this procedure has worked well in
those publications. We demonstrate our trial and error attempts below.

\begin{figure}[tbph]
\centering
\subcaptionbox{$a=2$
\label{fig:a2A}}{\includegraphics[width=0.23\textwidth]{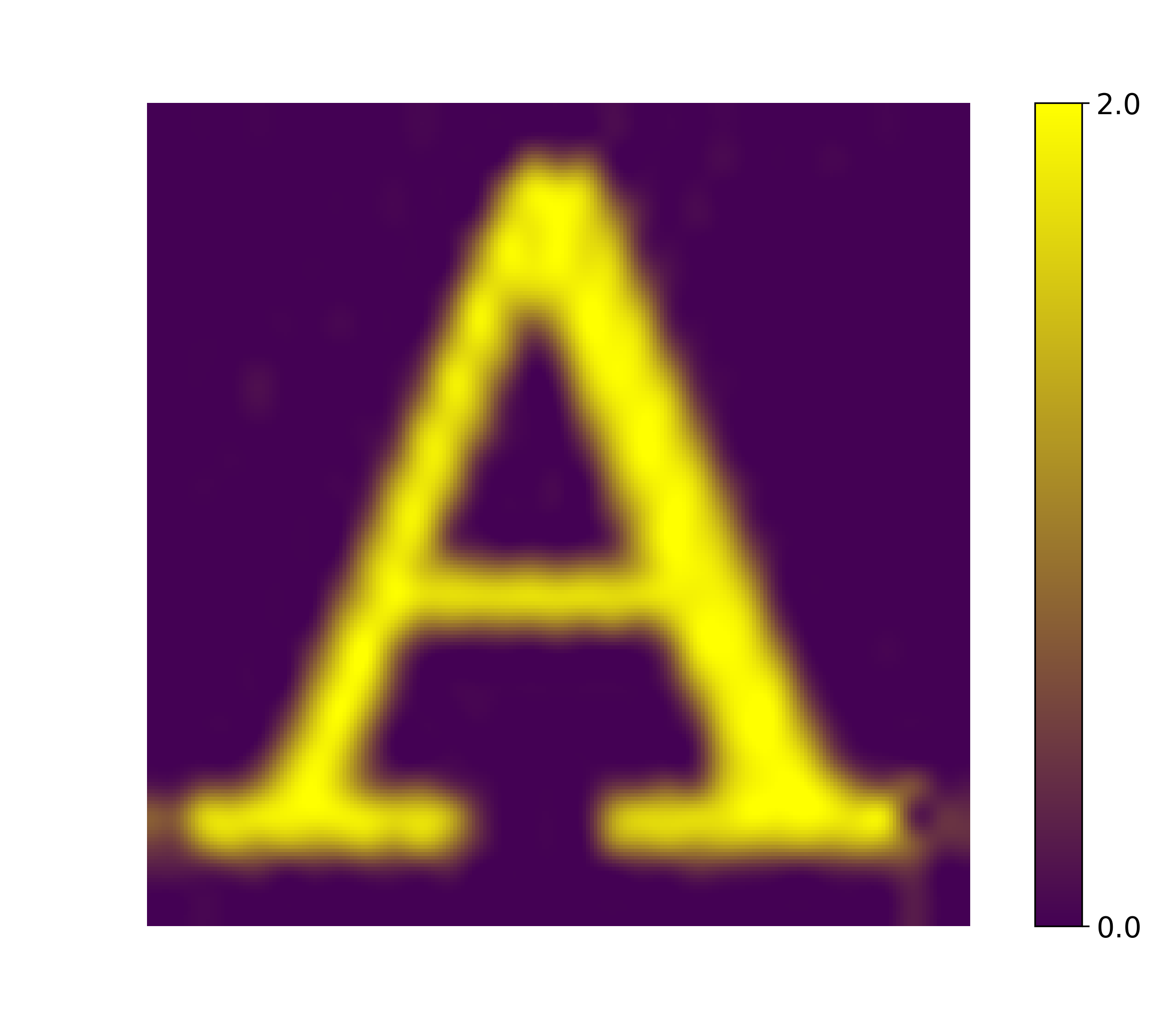}} 
\subcaptionbox{$a=3$
\label{fig:a3A}}{\includegraphics[width=0.23\textwidth]{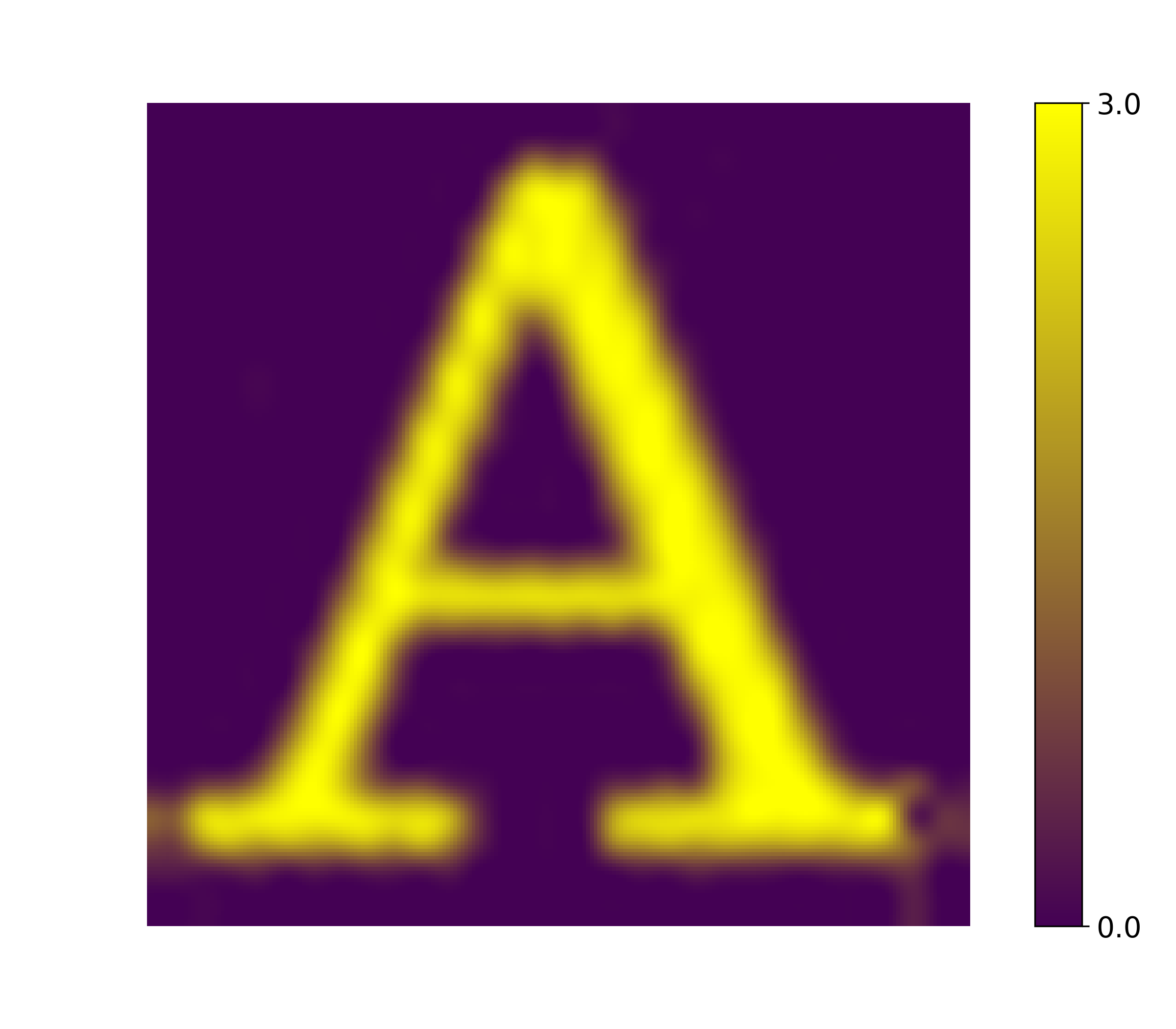}} 
\subcaptionbox{$a=5$
\label{fig:a5A}}{\includegraphics[width=0.23\textwidth]{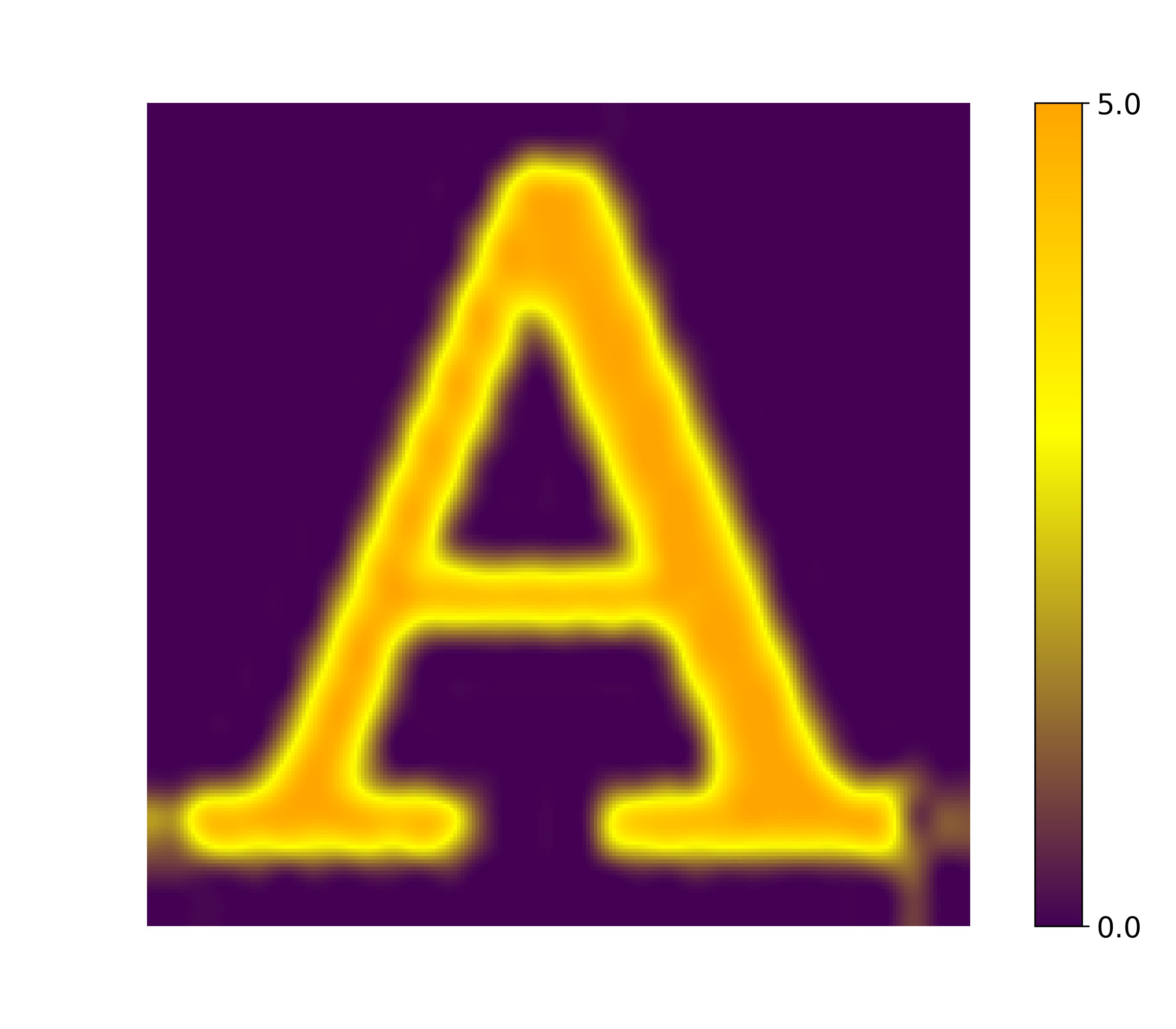}} %
\subcaptionbox{$a=10$\label{fig:a10A}}{\includegraphics[width=0.23%
\textwidth]{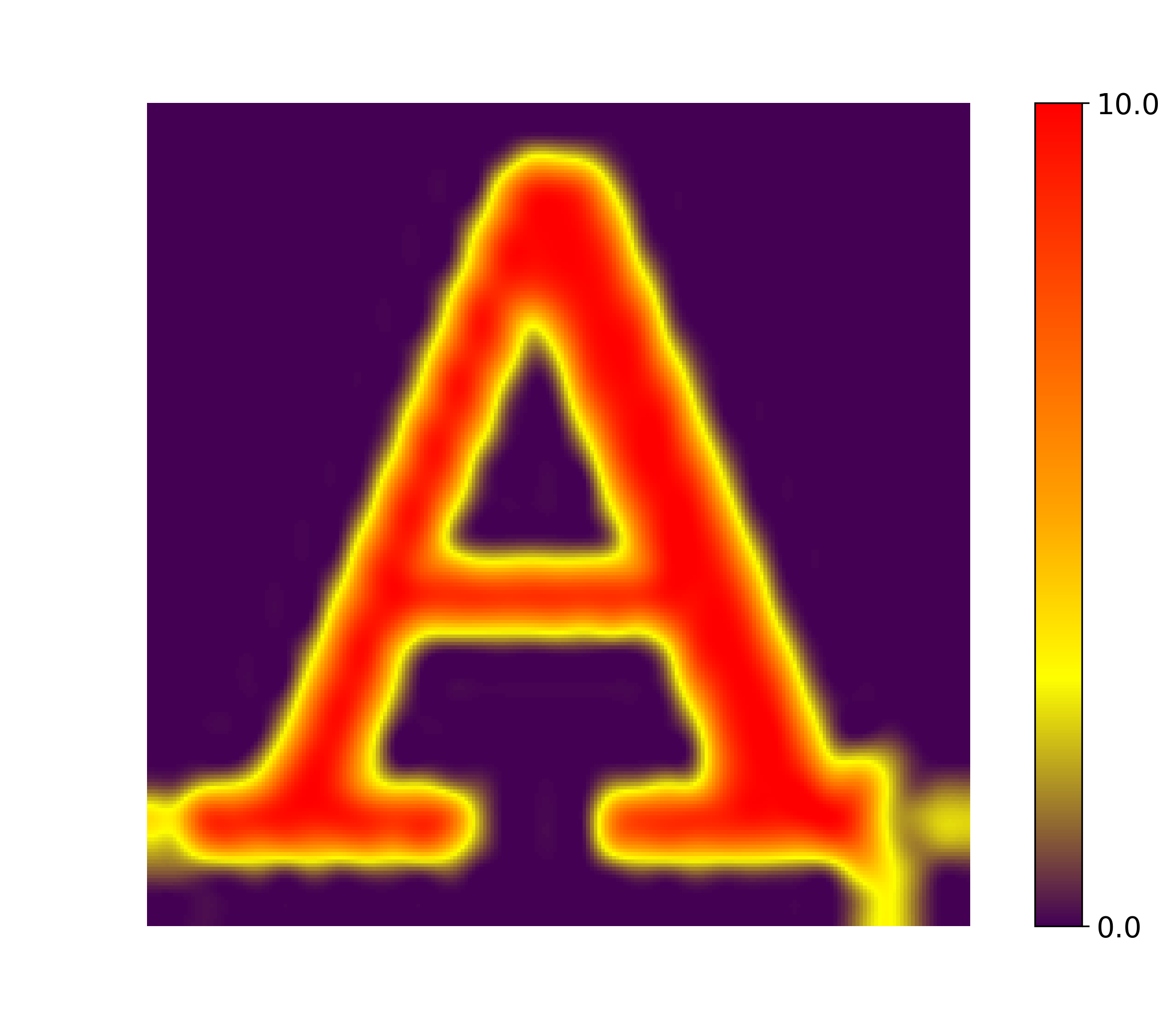}}
\par
\subcaptionbox{$a=2$
\label{fig:a2SZ}}{\includegraphics[width=0.23\textwidth]{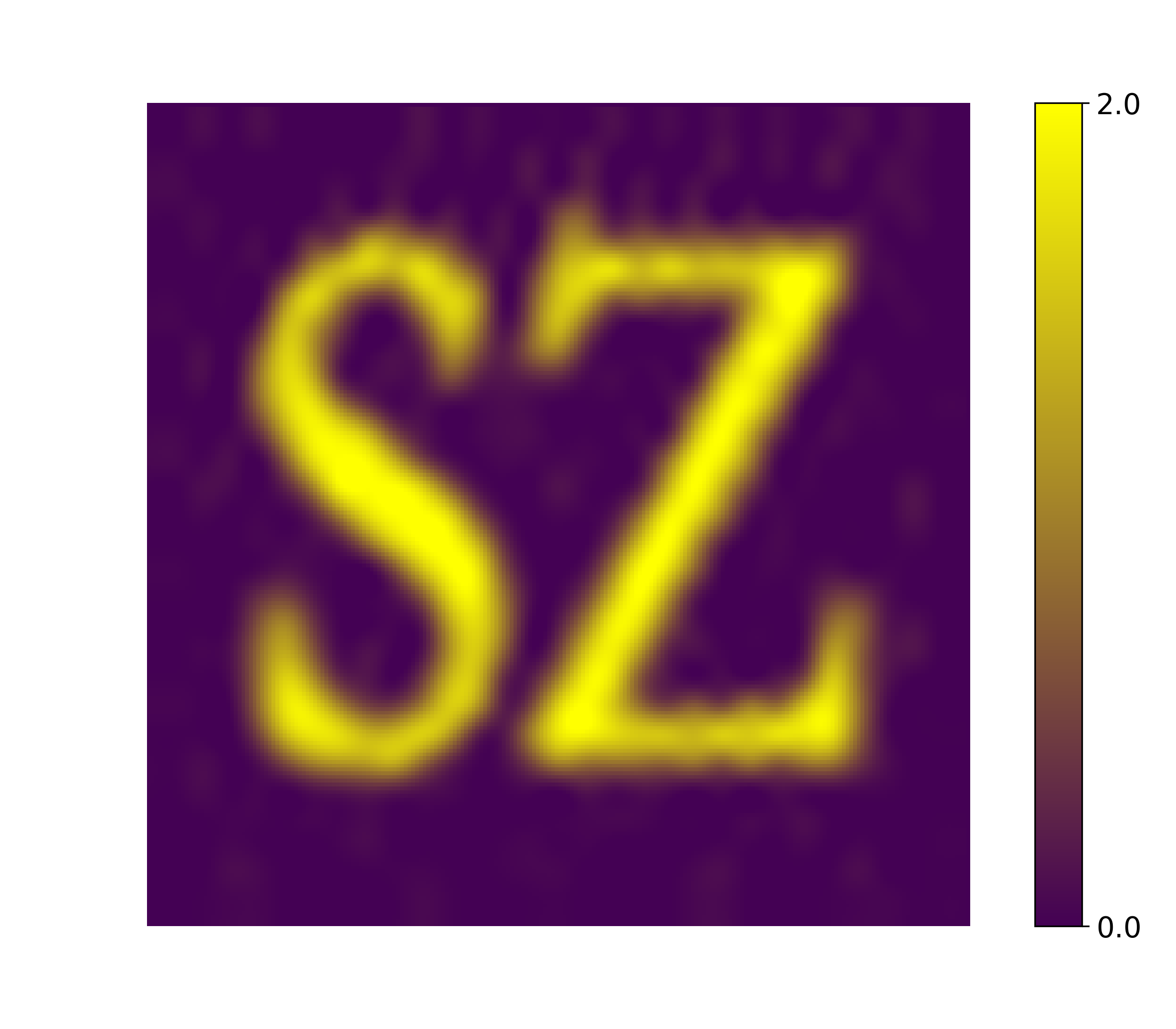}}
\subcaptionbox{$a=3$
\label{fig:a3SZ}}{\includegraphics[width=0.23\textwidth]{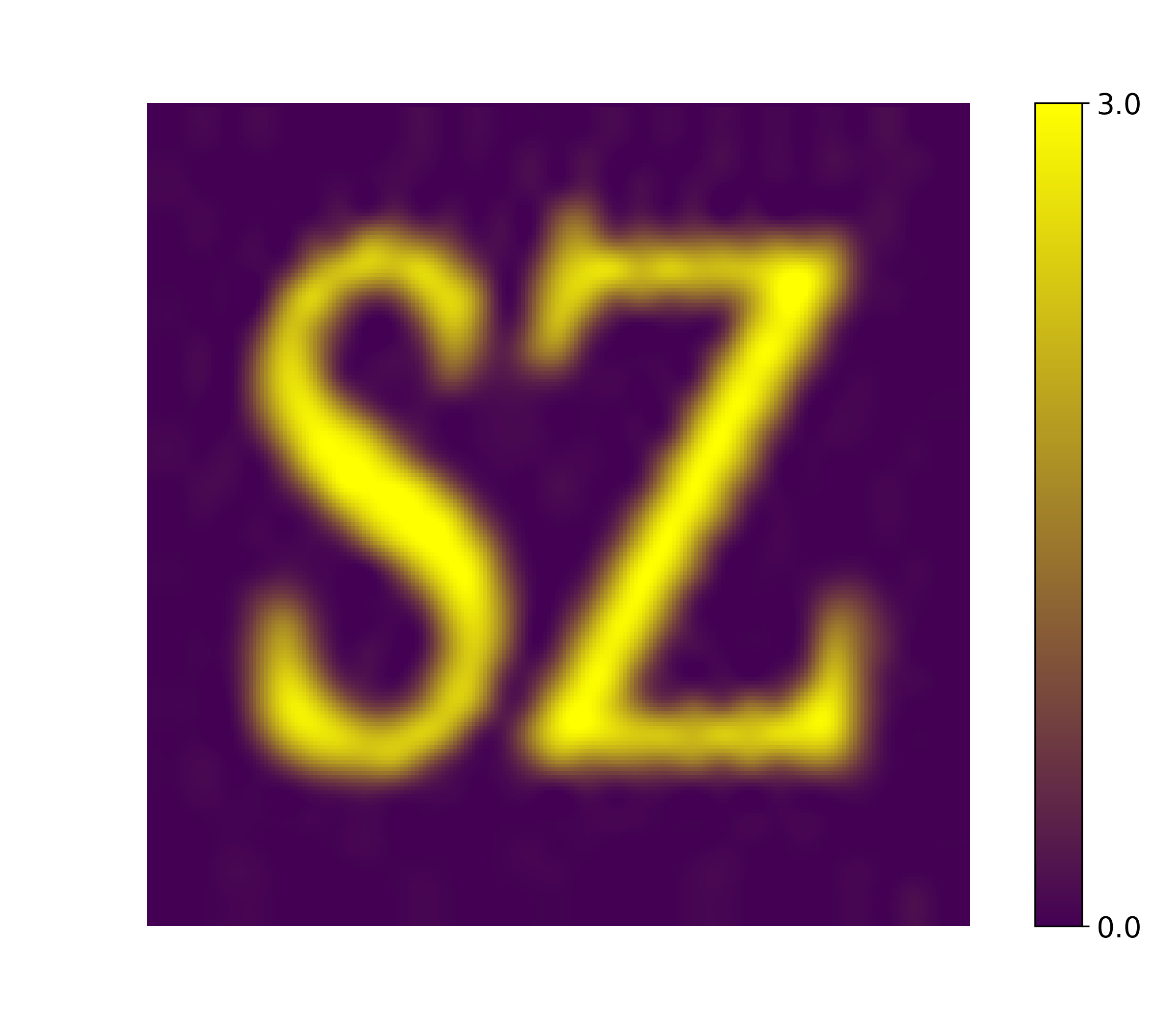}}
\subcaptionbox{$a=5$
\label{fig:a5SZ}}{\includegraphics[width=0.23\textwidth]{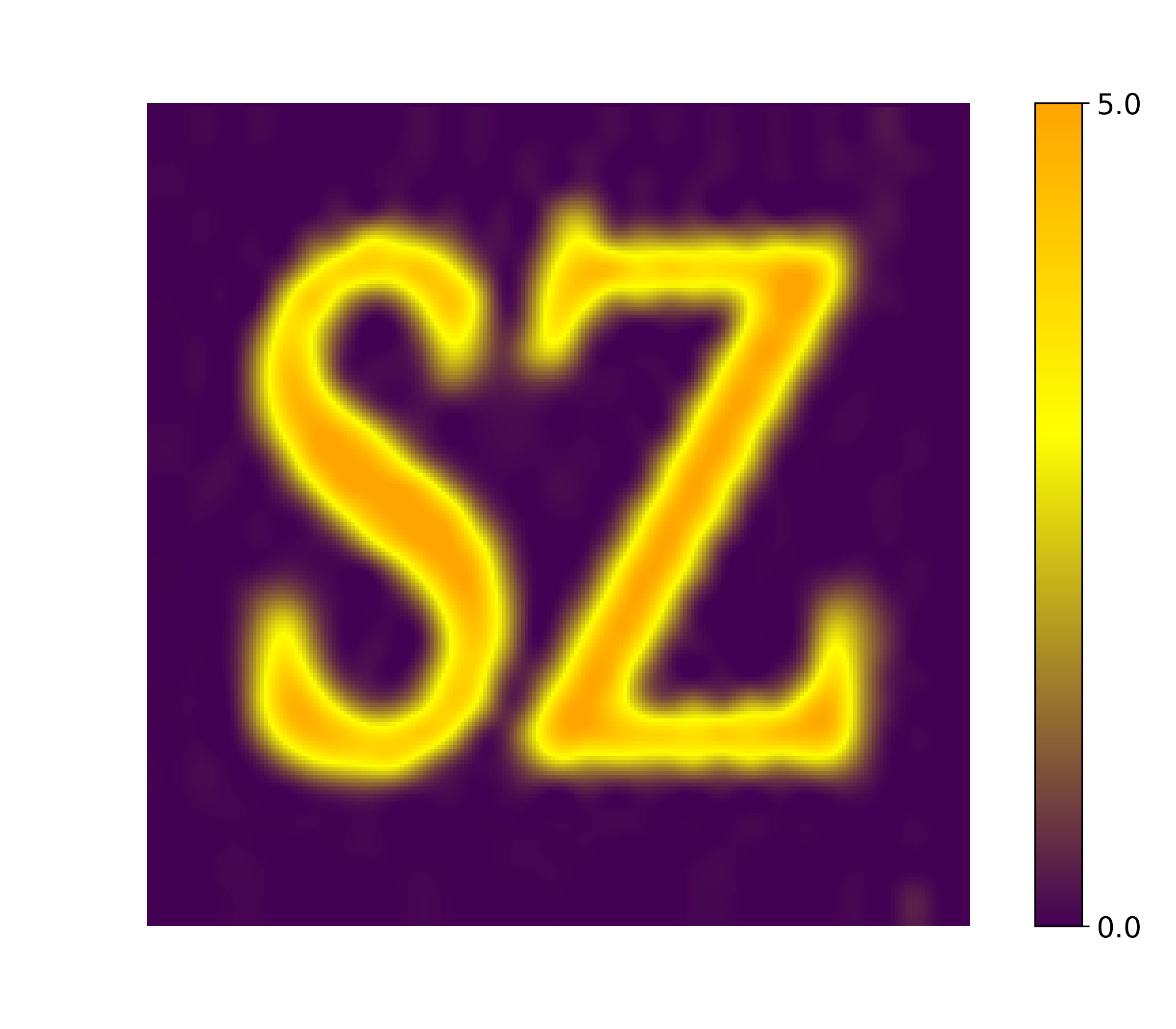}}
\subcaptionbox{$a=10$\label{fig:a10SZ}}{\includegraphics[width=0.23%
\textwidth]{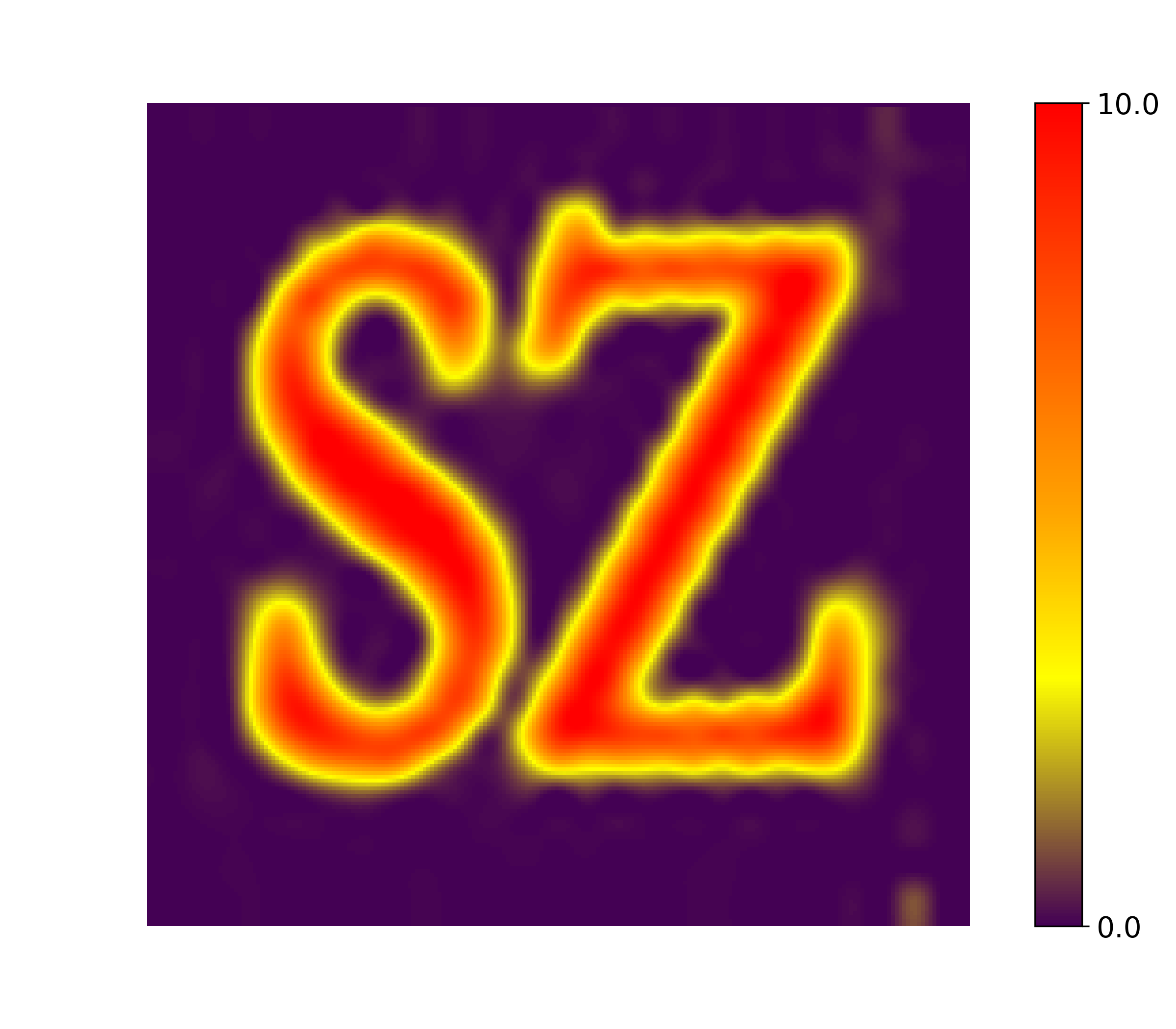}}
\caption{Results of the reconstructions illustrating the sensitivity of our
method to the value of the number $a$ in (\protect\ref{9.3}). Both shapes of
inclusions and the maximal values of the unknown coefficient $a(\mathbf{x})$
are accurately reconstructed.}
\label{fig:coeff_values}
\end{figure}

We conduct numerical experiments both in 2-d and 3-d cases. We specify the
domain $\Omega \subset \mathbb{R}^{n},n=2,3$ in (\ref{2.1}) as%
\begin{equation}
\left. 
\begin{array}{c}
\Omega =\left\{ \mathbf{x=}\left( x_{1},x_{2}\right)
:1<x_{1},x_{2}<2\right\} \text{ in }\mathbb{R}^{2}, \\ 
\Omega =\left\{ \mathbf{x=}\left( x_{1},x_{2}\right)
:1<x_{1},x_{2},x_{3}<2\right\} \text{ in }\mathbb{R}^{3}.%
\end{array}%
\right.  \label{9.1}
\end{equation}%
Then by (\ref{2.2}) and (\ref{9.1}) we have in the 3-d case:%
\begin{equation}
\left. 
\begin{array}{c}
\partial \Omega =\Gamma _{0}\cup \Gamma _{1}, \\ 
\Gamma _{0}=\left\{ x_{1}=2,x_{2}\in \left( 1,2\right) ,x_{3}\in \left(
1,2\right) \right\} , \\ 
\Gamma _{1}=\partial \Omega \diagdown \Gamma _{0}%
\end{array}%
\right. .  \label{9.2}
\end{equation}%
In the 2-D case $x_{3}$ is not present in (\ref{9.2}). In our data
simulations we have chosen the function $a\left( \mathbf{x}\right) $ as:%
\begin{equation}
a\left( \mathbf{x}\right) =\left\{ 
\begin{array}{c}
a=const.\geq 2\text{ inside of an inclusion,} \\ 
0\text{ otherwise.}%
\end{array}%
\right.  \label{9.3}
\end{equation}%
In our studies we have taken in (\ref{9.3})%
\begin{equation}
a=2,3,5,10.  \label{9.4}
\end{equation}%
The value $a=10$ is considered to be large. We got accurate images for all
these four values of $a$, see Figures \ref{fig:coeff_values}.

\textbf{Remark 8.1.} To better demonstrate a high degree of robustness of
our numerical method, we select letters-like shapes of inclusions we image.
Indeed, letters are non-convex and have voids in them.

We have worked numerically only with CIP2. First, we need to generate the
boundary data $g_{0}\left( \mathbf{x},t\right) $ and $g_{1}\left( \mathbf{x}%
,t\right) $ in (\ref{2.17}), which are our computationally simulated data.
To do this, we need to solve numerically Cauchy problem (\ref{2.15}), (\ref%
{2.16}). In our numerical studies, we approximate the $\delta -$function in (%
\ref{2.16}) by the $C^{\infty }\left( \mathbb{R}^{n}\right) -$function 
\begin{equation}
\delta _{\xi }\left( \mathbf{x}\right) =C_{\xi }\left\{ 
\begin{array}{c}
\exp \left( \frac{\left\vert \mathbf{x}\right\vert ^{2}}{\left\vert \mathbf{x%
}\right\vert ^{2}-\xi ^{2}}\right) ,\text{ }\left\vert \mathbf{x}\right\vert
<\xi , \\ 
0,\text{ }\left\vert \mathbf{x}\right\vert \geq \xi ,%
\end{array}%
\right.  \label{9.5}
\end{equation}%
where the parameter $\xi =0.05$, and the constant $C_{\xi }>0$ is chosen
such that 
\begin{equation*}
\displaystyle\int\limits_{\left\vert \mathbf{x}\right\vert <\xi }\delta
_{\xi }\left( \mathbf{x}\right) d\mathbf{x}=1.
\end{equation*}%
It is well known that \cite[formula (14.2) in chapter 4]{Lad} 
\begin{equation}
\lim_{\left\vert \mathbf{x}\right\vert \rightarrow \infty }u\left( \mathbf{x}%
,t\right) =0,  \label{9.6}
\end{equation}%
where $u\left( \mathbf{x},t\right) $ is the solution of problem (\ref{2.15}%
), (\ref{2.16}). Since we cannot perform computations in the infinite domain 
$\mathbb{R}^{n},$ then, using (\ref{9.5}) and (\ref{9.6}), we numerically
approximate the solution of problem (\ref{2.15}), (\ref{2.16}) by solving
its analog in a finite domain. We now specify this statement. Let $\Psi
\subset \mathbb{R}^{n}$ be a ball of the radius $r>0$ with its center at $%
(1.5,1.5,1.5)$ in the 3-d case and such that $\left\{ 0\right\} \in \Psi $.
In the 2-Dimensional case $\Psi $ is the disk with its center at $(1.5,1.5)$%
, see (\ref{9.1}). We solve equation (\ref{2.15}) with initial condition (%
\ref{9.5}) in $\Psi \times \left( 0,T\right) $ with the zero Dirichlet
boundary condition at $\partial \Psi \times \left( 0,T\right) ,$%
\begin{equation}
S_{t}=\Delta S+a\left( \mathbf{x}\right) S\text{, }\left( \mathbf{x}%
,t\right) \in \Psi \times \left( 0,T\right) ,  \label{9.7}
\end{equation}%
\begin{equation*}
S\left( \mathbf{x},0\right) =\delta _{\xi }\left( \mathbf{x}\right) ,
\end{equation*}%
\begin{equation}
S\mid _{\partial \Psi \times \left( 0,T\right) }=0.  \label{9.9}
\end{equation}%
We use the Finite Element Method for these computations with the spatial
mesh size $\widetilde{h}=0.01667$. The time interval $[0,T]$ is discretized
with 800 steps, where we set the final time 
\begin{equation}
T=4.  \label{9.90}
\end{equation}

The next question is: How to choose an appropriate radius $r_{appr}$\ of $\
\Psi ?$ Since we know the explicit form (\ref{2.18}) of the solution of
problem (\ref{2.15}), (\ref{2.16}) for $a\left( \mathbf{x}\right) \equiv 0$,
then we choose such a value of $r_{appr},$ for which our numerical solution
of problem (\ref{9.7})-(\ref{9.9}) with $a\left( \mathbf{x}\right) \equiv 0$
approximates well the right hand side of (\ref{2.18}) for $\left( \mathbf{x}%
,t\right) \in \Psi \times \left( 0,T\right) .$ These considerations led us
to find the optimal value $r_{appr}=6$ in both 3-Dimensional and
2-Dimensional cases. Next, we assign $u\left( \mathbf{x},t\right) :=S\left( 
\mathbf{x},t\right) $ for $\left( \mathbf{x},t\right) \in \Psi \times \left(
0,T\right) .$

The value of the parameter $\varepsilon $ in (\ref{3.2}) is chosen
numerically based on the analysis of Test \ref{ex:epsilon_sensitivity}.
Random noise is added to the Dirichlet and Neumann boundary data $g_{0}$ and 
$g_{1}$ (\ref{2.17}) to simulate realistic measurement conditions.
Specifically, the noisy data $g_{0}^{\sigma }$ and $g_{1}^{\sigma }$ are
generated by 
\begin{equation}
g_{0}^{\sigma }(\mathbf{x},t)=g_{0}(\mathbf{x},t)(1+\sigma \zeta ^{(0)}),
\label{9.10}
\end{equation}%
\begin{equation}
g_{1}^{\sigma }(\mathbf{x},t)=g_{1}(\mathbf{x},t)(1+\sigma \zeta ^{(1)}),
\label{9.11}
\end{equation}%
where $\sigma \in \left( 0,1\right) $ is the noise level and $\zeta
^{(0)},\zeta ^{(1)}$ are uniformly distributed random variables in $[-1,1]$.
Hence, e.g. $\sigma =0.01$ corresponds to the 1\% noise level.

For the inverse problem, we discretize the spatial domain $\Omega =(1,2)^{n}$
using a uniform grid with 20 points in each coordinate direction:
specifically, $20\times 20$ for the 2-Dimensional case and $20\times
20\times 20$ for the 3-Dimensional case. The discretization of the time
interval $[\varepsilon ,T]$ to obtain the above system of elliptic PDEs
requires a careful consideration, as discussed in Test \ref%
{ex:TNode_sensitivity}. We write the differential operators in the
functional $J_{\lambda ,\alpha }\left( V\right) $ in (\ref{4.6}) in the form
of finite differences with the above grid points inside of the domain $%
\Omega $. Next, we minimize the resulting discretized functional $J_{\lambda
,\alpha }^{disc}\left( V\right) $ with respect to the values of the vector
function $V\left( \mathbf{x}\right) $ at those grid points.

The starting point $V^{(0)}(\mathbf{x})$ must satisfy the Dirichlet boundary
conditions generated by $g_{0}$ and capture the short-time asymptotic
behavior of the function $v\left( \mathbf{x},t\right) $ at $t=\varepsilon $.
To handle the latter, we temporary assume that $v\left( \mathbf{x}%
,\varepsilon \right) \approx v_{bg}\left( \mathbf{x},\varepsilon \right)
=\partial _{t}\left( \ln u_{0}\left( \mathbf{x},t\right) \right) \mid
_{t=\varepsilon },$ see (\ref{2.18}), (\ref{2.20}), (\ref{3.4}) and (\ref%
{3.9}). However, we do not use this assumption neither in the above theory
nor in the iterations, which follow our first step. Consider the 2-d case
first and temporary denote $\mathbf{x}=\left( x,y\right) .$ Accordingly, our
starting vector $V^{(0)}(\mathbf{x})=(v_{0}^{(0)}(\mathbf{x}),v_{1}^{(0)}(%
\mathbf{x}),\dots ,v_{k}^{(0)}(\mathbf{x}))^{T}$, where each component for
the time step $t_{j}$ is constructed as: 
\begin{equation*}
\begin{aligned} v_j^{(0)}(x,y) &= \gamma_j v_{\mathrm{bg}}(x,y,\epsilon) +
\\ &\quad +
\left(\frac{x_{\max}-x}{x_{\max}-x_{\min}}\tilde{g}_{0,j}(x_{\min},y) +
\frac{x-x_{\min}}{x_{\max}-x_{\min}}\tilde{g}_{0,j}(x_{\max},y)\right) \\
&\quad +
\left(\frac{y_{\max}-y}{y_{\max}-y_{\min}}\tilde{g}_{0,j}(x,y_{\min}) +
\frac{y-y_{\min}}{y_{\max}-y_{\min}}\tilde{g}_{0,j}(x,y_{\max})\right) - \\
&\quad - \tilde{v}_{\mathrm{corners},j}(x,y), \end{aligned}
\end{equation*}%
where 
\begin{equation*}
\gamma _{j}=\frac{T-t_{j}}{T-\varepsilon },\text{ }v_{\mathrm{bg}%
}(x,y,\varepsilon )=\frac{x^{2}+y^{2}}{4\varepsilon ^{2}}-\frac{1}{%
\varepsilon },\text{ }\tilde{g}_{0,j}(x,y)=g_{0}(x,y,t_{j})-\gamma _{j}v_{%
\mathrm{bg}}(x,y,\varepsilon ),
\end{equation*}
and $\tilde{v}_{\mathrm{corners},j}(x,y)$ is the standard bilinear
interpolation of $\tilde{g}_{0,j}$ at the four corners of the spatial
domain. For our numerical studies, $x_{\min }=y_{\min }=1$, $x_{\max
}=y_{\max }=2$. The initial guess for the 3-d case is constructed in a
similar manner.

We enforce the Neumann boundary condition via the second-order forward
difference. Let the boundary $\Gamma _{0}$ be indexed by $I$ and let $\tilde{%
h}$ be the spatial step. The values at the first interior layer $I-1$ are
given by 
\begin{equation*}
v_{I-1,j,\ell }=\frac{3}{4}g_{0}(\mathbf{x}_{I,j},t_{\ell })+\frac{1}{4}%
v_{I-2,j,\ell }-\frac{\tilde{h}}{2}\,g_{1}(\mathbf{x}_{I,j},t_{\ell }),
\end{equation*}%
where $v_{I-2,j,\ell }$ are free variables. This relation is applied after
each iteration, ensuring the Neumann condition holds exactly. By embedding
all boundary constraints in this manner, the constrained minimization
problem is transformed into an unconstrained one. We then apply the L-BFGS
quasi-Newton algorithm, implemented in Python via PyTorch, to solve this
unconstrained minimization problem for the functional $J_{\lambda ,\alpha
}\left( V\right) $ in (\ref{4.6}). This approach is highly advantageous for
large-scale optimization owing to its low memory requirements.

\begin{figure}[htbp]
\centering
\includegraphics[width=0.8\textwidth]{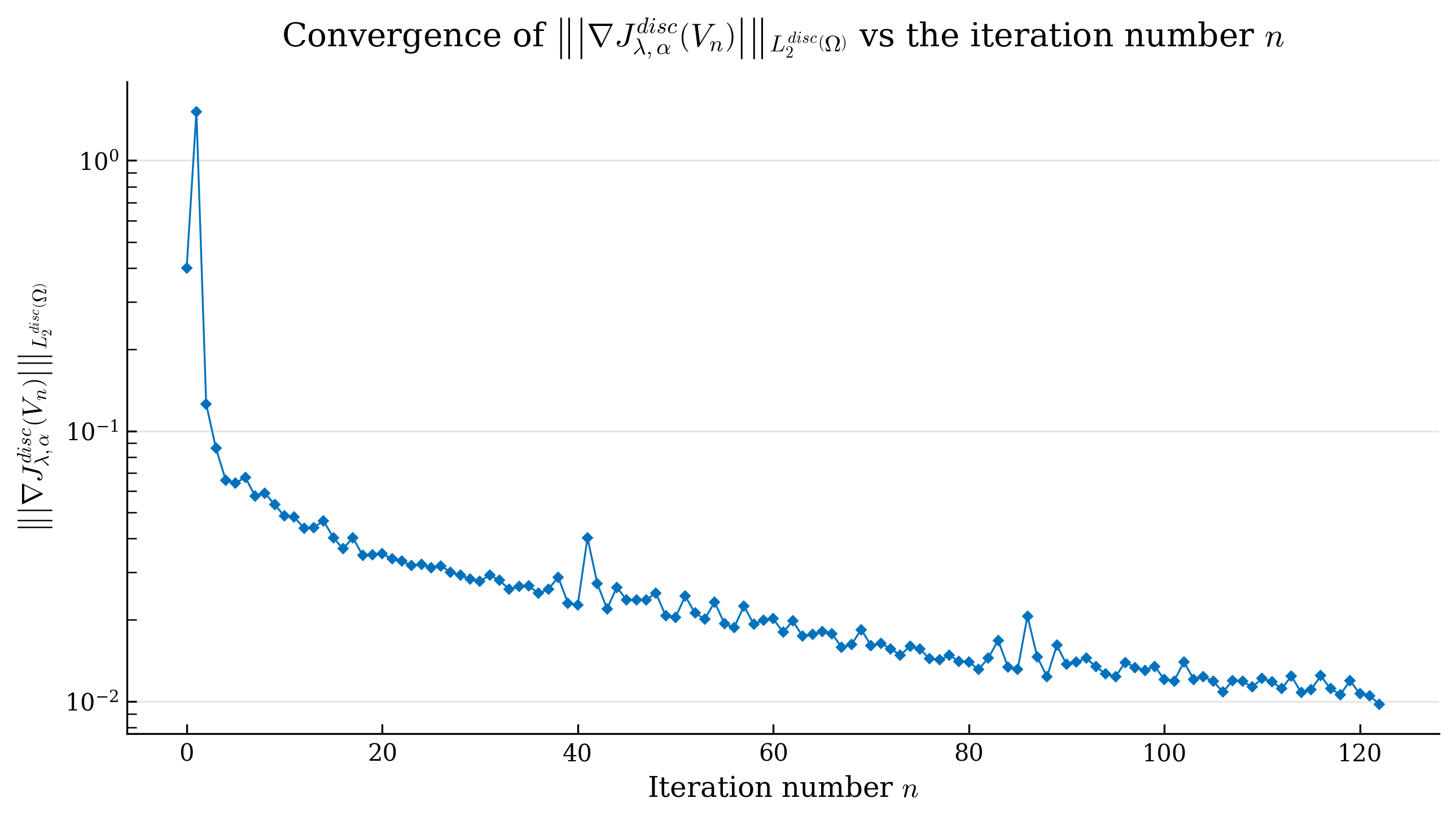}
\caption{A typical convergence behavior of $\left\Vert \left\vert \protect%
\nabla J_{\protect\lambda ,\protect\alpha }^{disc}\left( V_{m}\right)
\right\vert \right\Vert _{L_{2}^{disc}\left( \Omega \right) }$ with respect
to the iteration number n of iterations of the L-BFGS algorithm. Note that
the value of the norm $\left\Vert \left\vert \protect\nabla J_{\protect%
\lambda ,\protect\alpha }^{disc}\left( V_{m}\right) \right\vert \right\Vert
_{L_{2}^{disc}\left( \Omega \right) }$ decreases by the factor of 100 due to
the global convergence. This figure explains our stopping criterion (\protect
\ref{9.100}).}
\label{fig:gradient_convergence}
\end{figure}

The stopping criterion of our minimization procedure was 
\begin{equation}
\left\Vert \left\vert \nabla J_{\lambda ,\alpha }^{disc}\left( V_{m}\right)
\right\vert \right\Vert _{L_{2}^{disc}\left( \Omega \right) }\leq 0.01,
\label{9.100}
\end{equation}%
where $m$ is the stopping iteration number, $\left\vert \nabla J_{\lambda
,\alpha }^{disc}\left( V_{m}\right) \right\vert $ is the magnitude of the
gradient of the discretized functional at $V_{m}$, and $\left\Vert \cdot
\right\Vert _{L_{2}^{disc}\left( \Omega \right) }$ is the discrete analog of
the $L_{2}\left( \Omega \right) -$norm. Figure \ref{fig:gradient_convergence}
displays a typical convergence behavior for all the tests below and explains
our stopping criterion (\ref{9.100}). Note that the value of the norm $%
\left\Vert \left\vert \nabla J_{\lambda ,\alpha }^{disc}\left( V_{m}\right)
\right\vert \right\Vert _{L_{2}^{disc}\left( \Omega \right) }$ decreases by
the factor of 100 due to the global convergence.

We begin by investigating the choice of the grid step size $h$ in (\ref{3.16}%
) with respect to time, which is a critical factor for the quality of the
reconstruction. Indeed, $h$ determines the number of resulting elliptic PDEs
(\ref{3.25})-(\ref{3.30}). We note that we have selected an optimal value of
the parameter $c$ in (\ref{4.6}) the same way as we select optimal values of
other parameters below. Hence, we do not describe this selection for
brevity. More precisely that value is $c=5$ for all tests below.

\begin{test}
\label{ex:TNode_sensitivity} We investigate the impact of the temporal
discretization by recovering an inclusion in the shape of the letter `$B$',
for which the true coefficient is $a(\mathbf{x})=2$, see (\ref{9.3}) and (%
\ref{9.4}). In this test, $\lambda =3$ in CWF (\ref{4.4}), and the Tikhonov
regularization parameter is set to $\alpha =3\times 10^{-5}$. We set the
parameter $\varepsilon =0.01$ in (\ref{3.2}), $c=5$ in (\ref{4.6}) and the
final time $T=4,$ as in (\ref{9.90}). To simulate a realistic scenario, we
add $1\%$ noise to the measurement data, i.e. $\sigma =0.01$ in (\ref{9.10}%
), (\ref{9.11}). We perform reconstructions for the CIP2 using three
different numbers of time steps: $N_{t}\in \{10,20,40\}$. The results are
displayed on Figure~\ref{fig:timestep_test}.

\begin{figure}[tbph]
\centering
\begin{subfigure}[b]{0.23\textwidth}
			\centering
			\includegraphics[width=\textwidth]{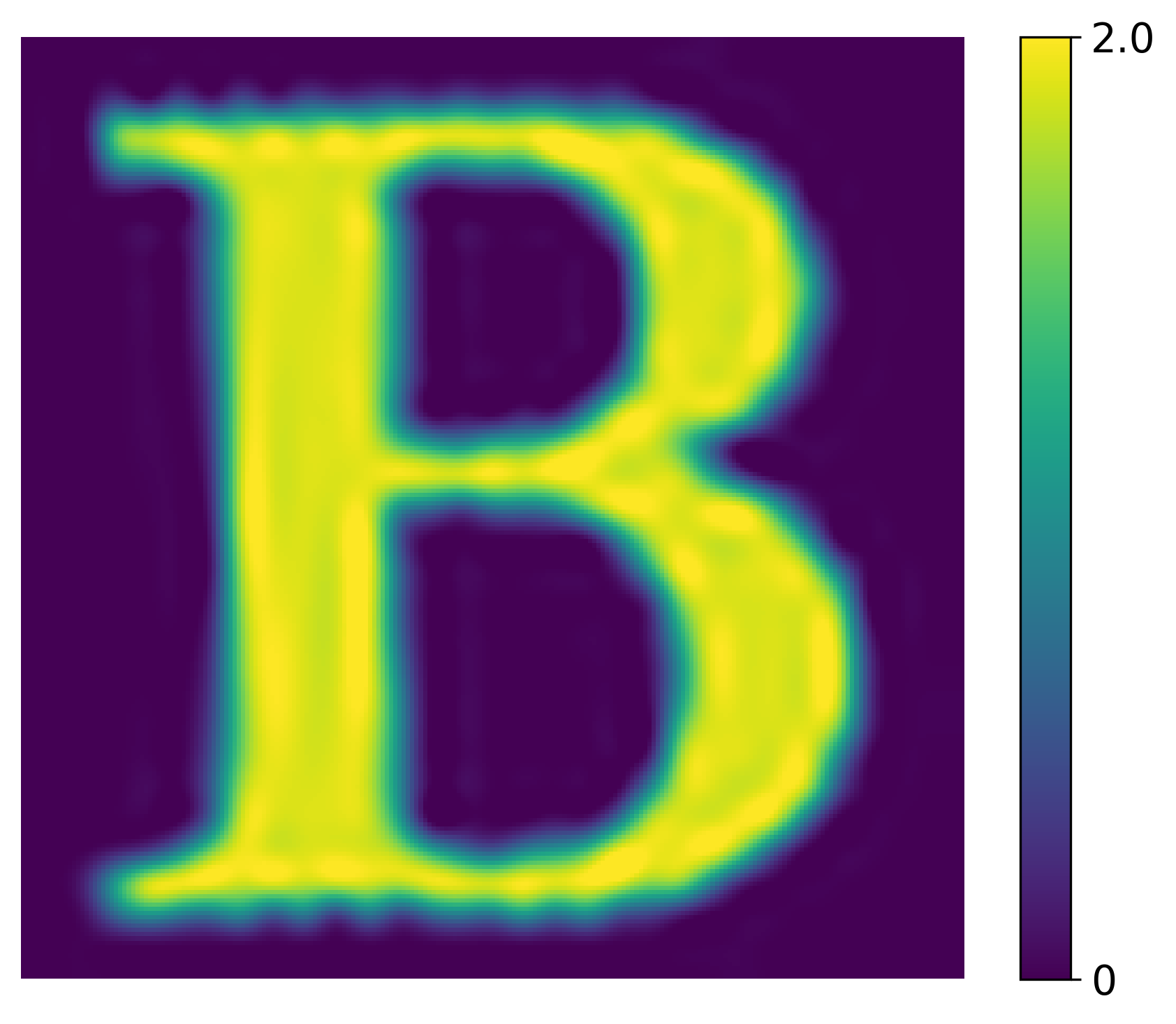}
			\caption{True image}
		\end{subfigure}
\hfill 
\begin{subfigure}[b]{0.23\textwidth}
			\centering
			\includegraphics[width=\textwidth]{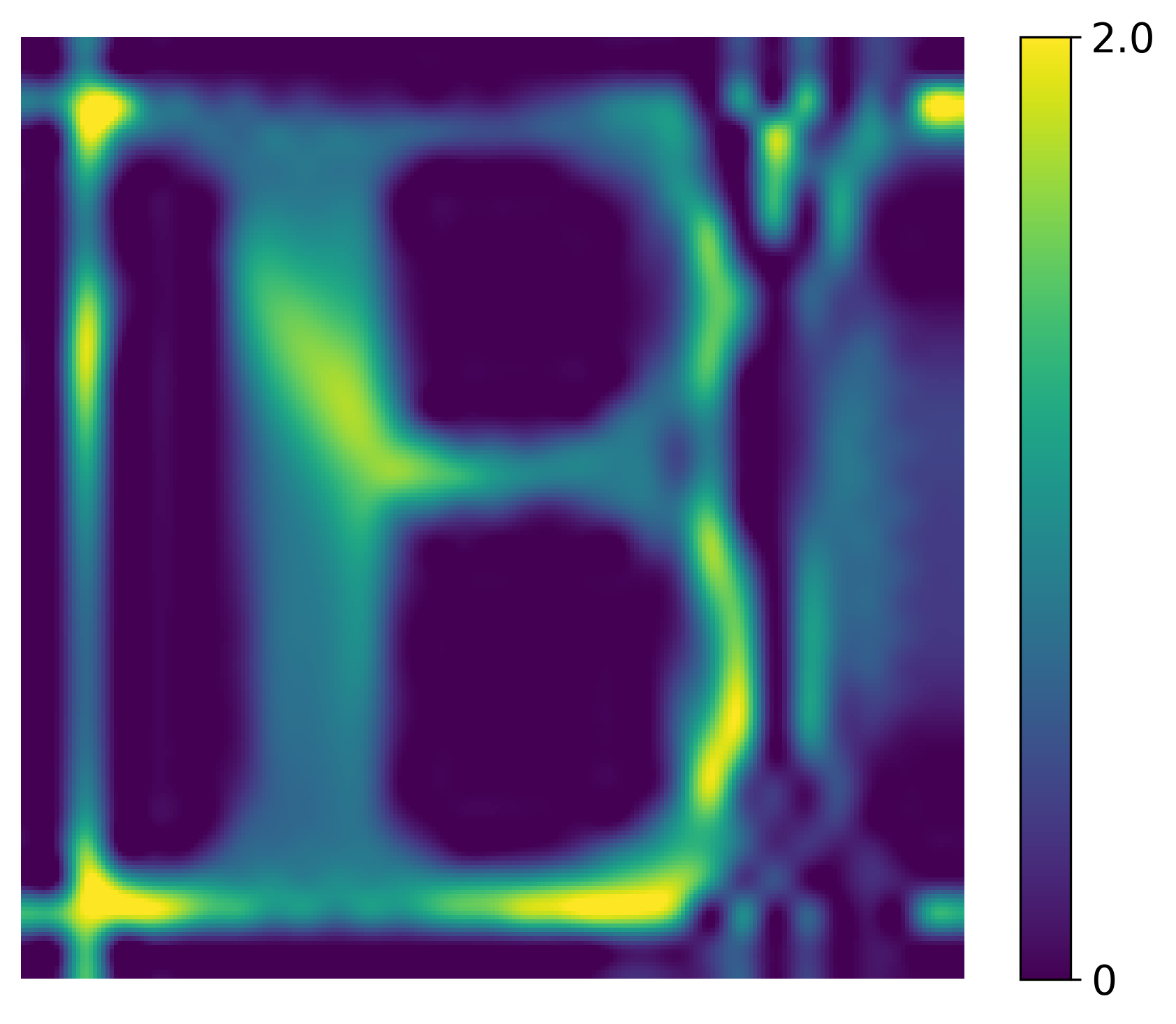}
			\caption{$N_t = 10$}
		\end{subfigure}
\hfill 
\begin{subfigure}[b]{0.23\textwidth}
			\centering
			\includegraphics[width=\textwidth]{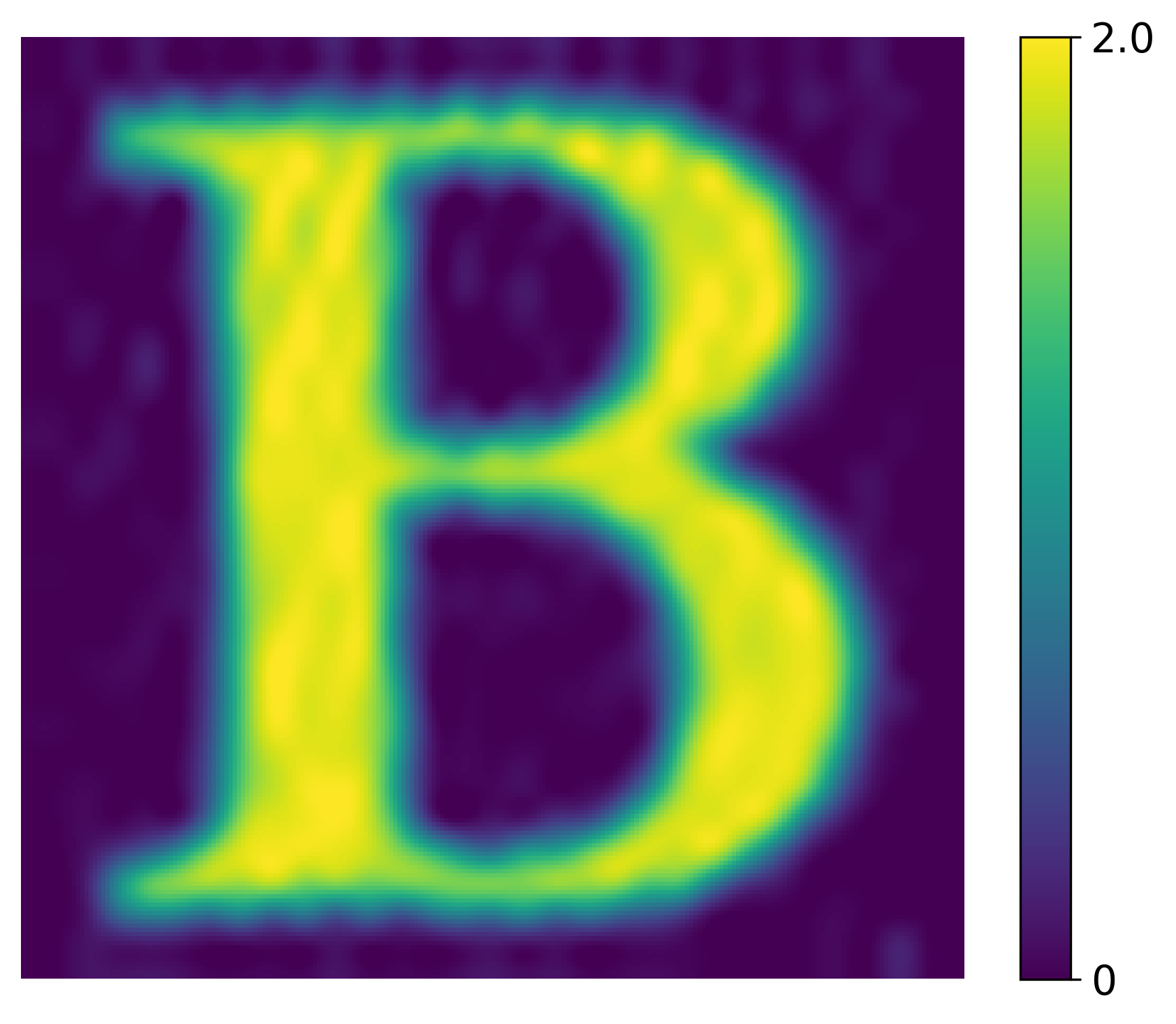}
			\caption{$N_t = 20$}
		\end{subfigure}
\hfill 
\begin{subfigure}[b]{0.23\textwidth}
			\centering
			\includegraphics[width=\textwidth]{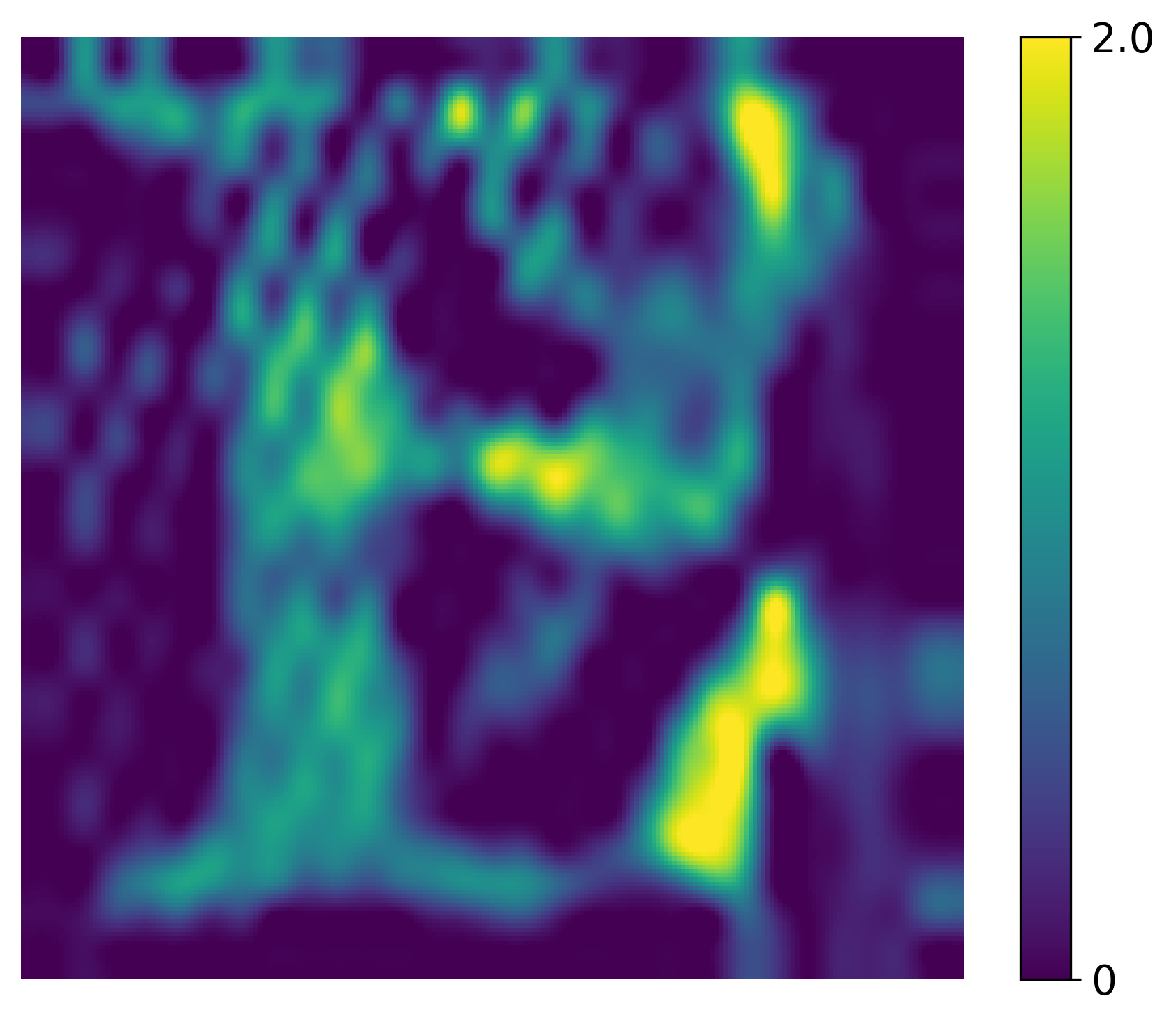}
			\caption{$N_t = 40$}
		\end{subfigure}
\caption{Reconstruction results for different numbers of time steps $%
N_{t}=T/h=4/h,$ where $T=4$ as in (\protect\ref{9.90}), and $h$ is as in (%
\protect\ref{3.16}). Clearly, $N_{t}=20$ is the optimal number.}
\label{fig:timestep_test}
\end{figure}
\end{test}

For $N_{t}=40$, the result is heavily contaminated by the numerical noise,
rendering the image heavily corrupted by artifacts. For $N_{t}=10$, the
reconstruction suffers from a significant blur. The clearest and the most
accurate image is obtained with $N_{t}=20$. This finding highlights a
counter-intuitive phenomenon of CIP2. Indeed, the idea that by decreasing
the grid step size $h=4/N_{t}$ in (\ref{3.16}) one would automatically
improve the image quality does not work here. In fact, we demonstrate that
the opposite is true for this ill-posed problem. Apparently, the step size $%
h $ in (\ref{3.16}) acts as a sort of a regularization parameter here by
filtering out high-frequency instabilities. Consequently, we fix $N_{t}=20$
for all subsequent numerical tests, which means $h=4/20=0.2.$

In Test 8.2 we numerically find the value of the optimal parameter $%
\varepsilon $ in (\ref{3.2}).

\begin{test}
\label{ex:epsilon_sensitivity}

We now find an optimal value of the parameter $\varepsilon $ in formula (\ref%
{3.2}). We test the case when the inclusion has the shape of the letter $%
\Omega $, for which the true coefficient $a(\mathbf{x})=2$ inside of this
inclusion and $a(\mathbf{x})=0$ outside of it, see (\ref{9.3}) and (\ref{9.4}%
). In this test, $\lambda =3$ in CWF (\ref{4.4}), the Tikhonov
regularization parameter is $\alpha =3\times 10^{-5}$ and $N_{t}=20$. We
again add 1\% noise to the measurement data, i.e. $\sigma =0.01$ in (\ref%
{9.10}), (\ref{9.11}). The parameter $\varepsilon $ is varied over the set $%
\{0.001,0.01,0.03,0.05,0.1\}$, while the final time is fixed at $T=4$ as in (%
\ref{9.90}) and $c=5$.

\begin{figure}[th]
\centering
\begin{subfigure}[b]{0.32\textwidth}
			\centering
			\includegraphics[width=\textwidth]{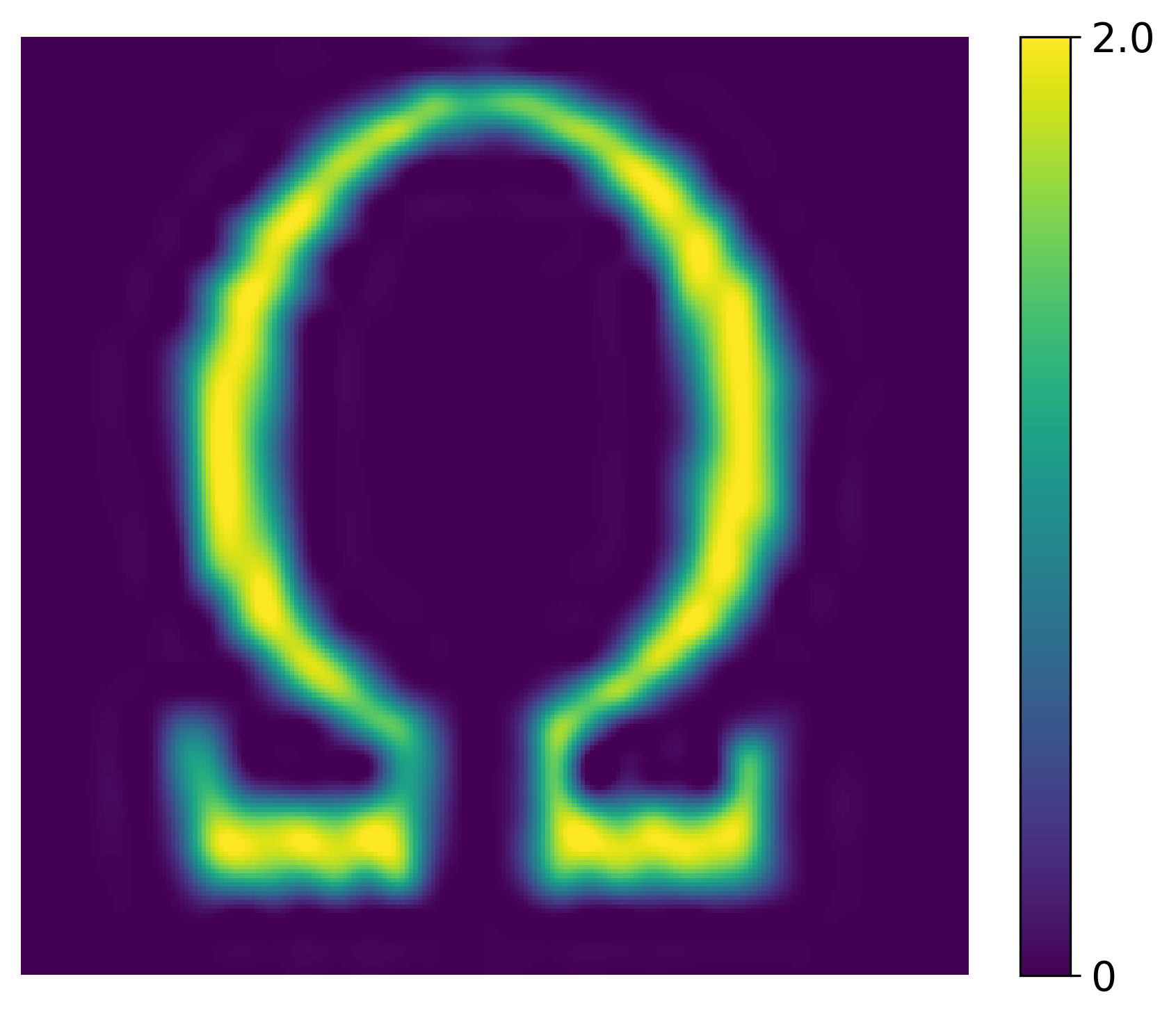}
			\caption{True image}
			\label{original letter}
		\end{subfigure}
\hfill 
\begin{subfigure}[b]{0.32\textwidth}
			\centering
			\includegraphics[width=\textwidth]{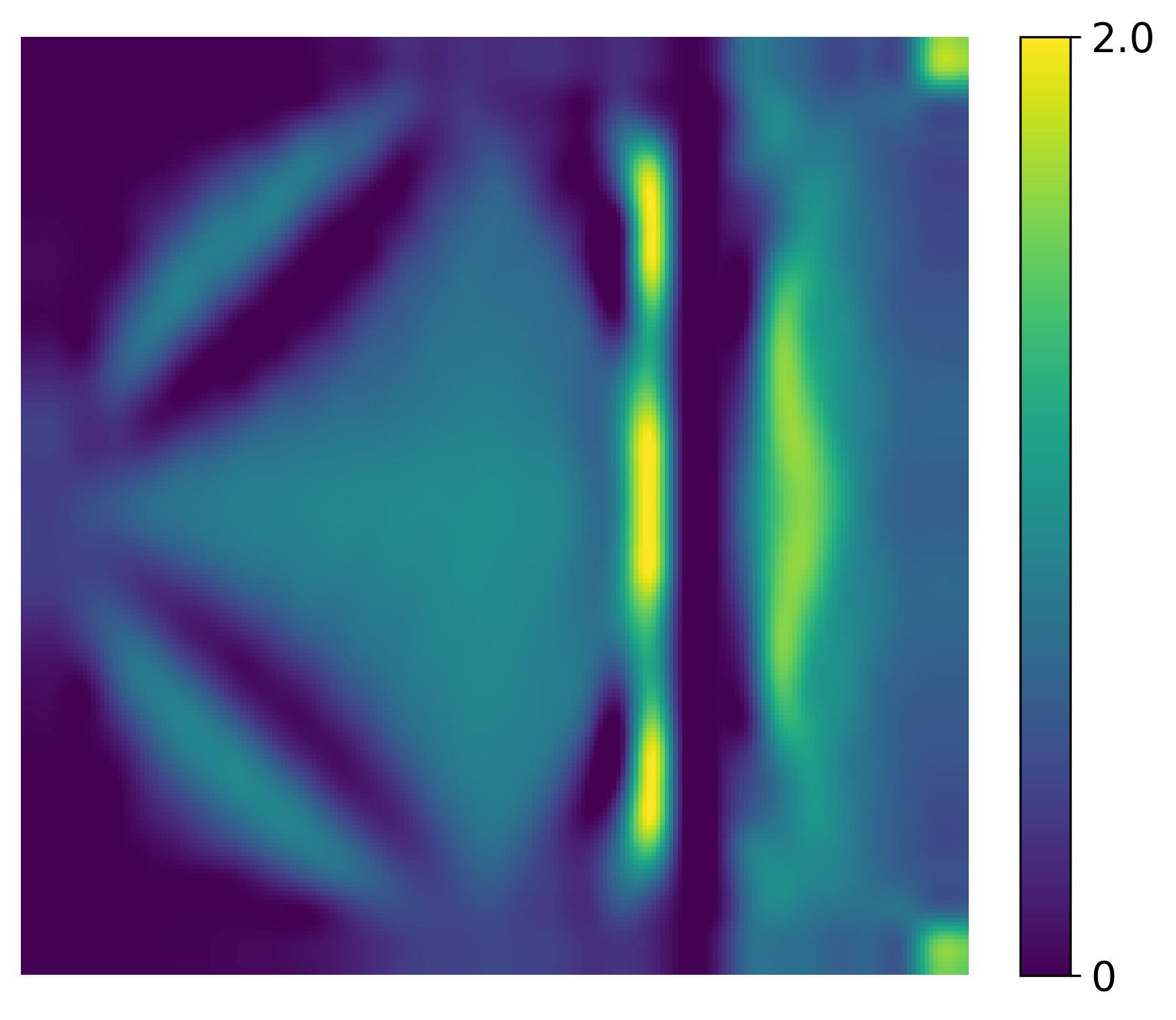}
			\caption{$\varepsilon=0.001$}
			\label{fig:eps_0.01}
		\end{subfigure}
\hfill 
\begin{subfigure}[b]{0.32\textwidth}
			\centering
			\includegraphics[width=\textwidth]{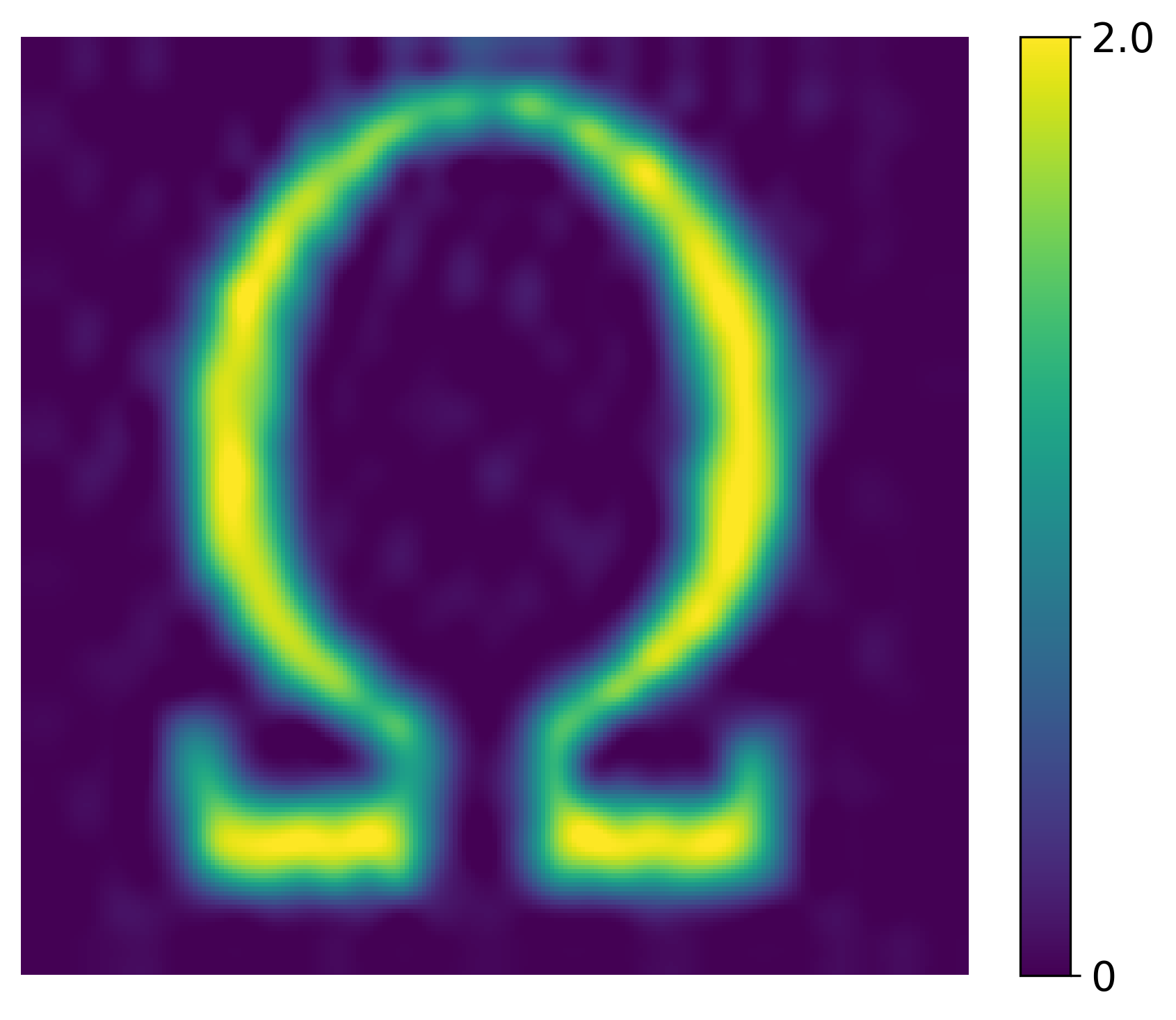}
			\caption{$\varepsilon=0.01$}
			\label{fig:eps_0.03}
		\end{subfigure}
\vspace{1em} 
\begin{subfigure}[b]{0.32\textwidth}
			\centering
			\includegraphics[width=\textwidth]{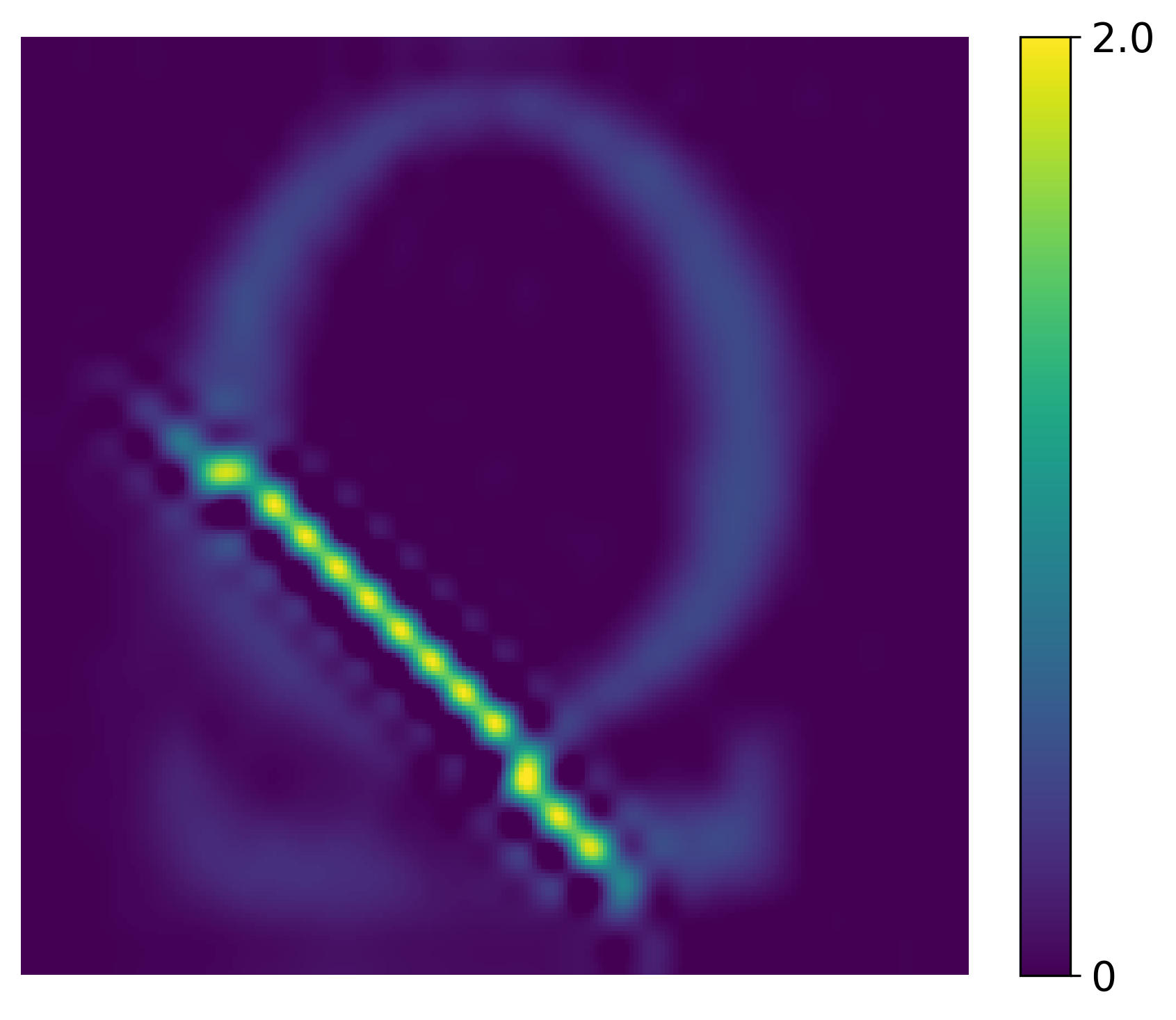}
			\caption{$\varepsilon=0.03$}
			\label{fig:eps_0.05}
		\end{subfigure}
\hfill 
\begin{subfigure}[b]{0.32\textwidth}
			\centering
			\includegraphics[width=\textwidth]{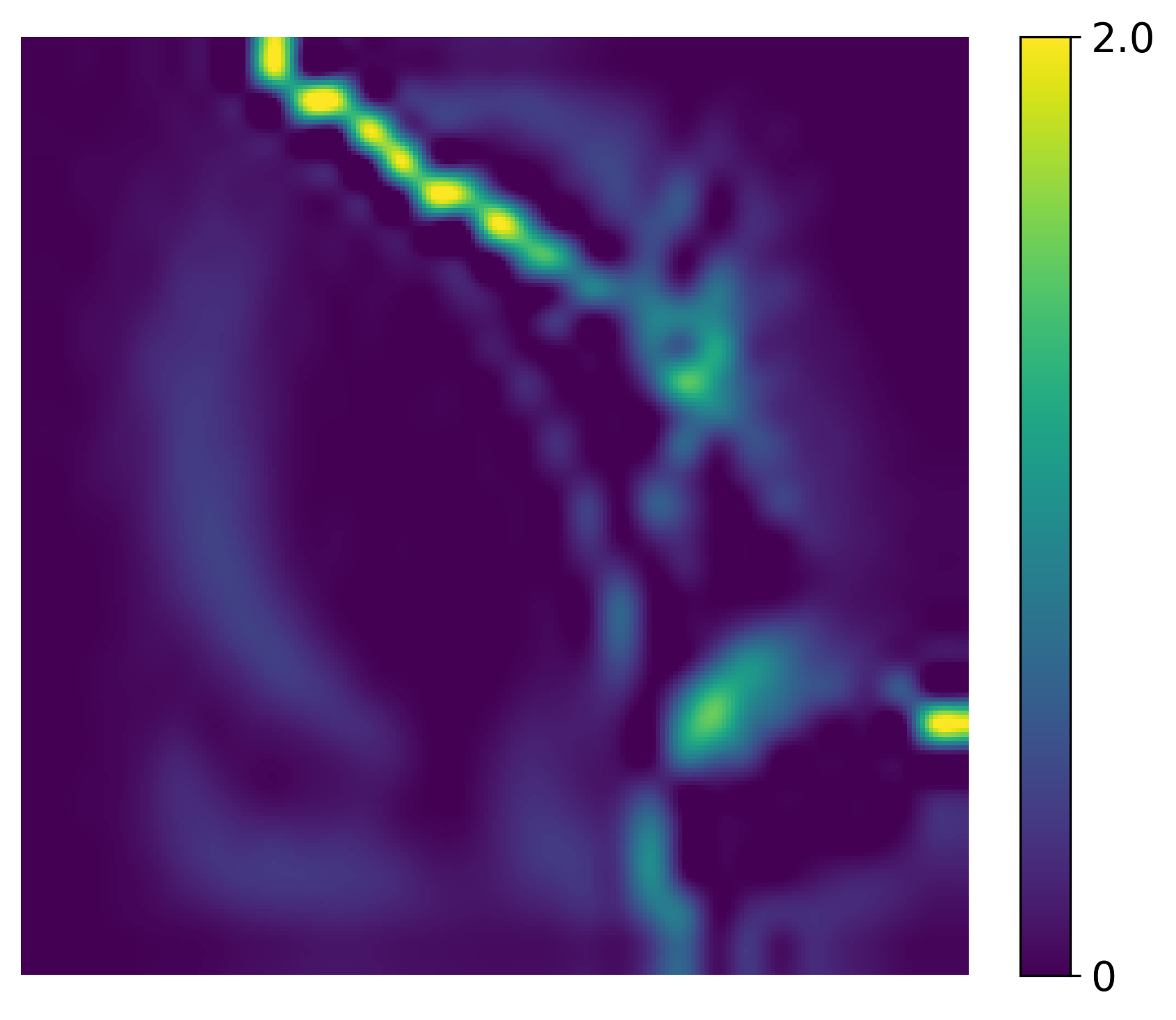}
			\caption{$\varepsilon=0.05$}
			\label{fig:eps_0.15}
		\end{subfigure}
\hfill 
\begin{subfigure}[b]{0.32\textwidth}
			\centering
			\includegraphics[width=\textwidth]{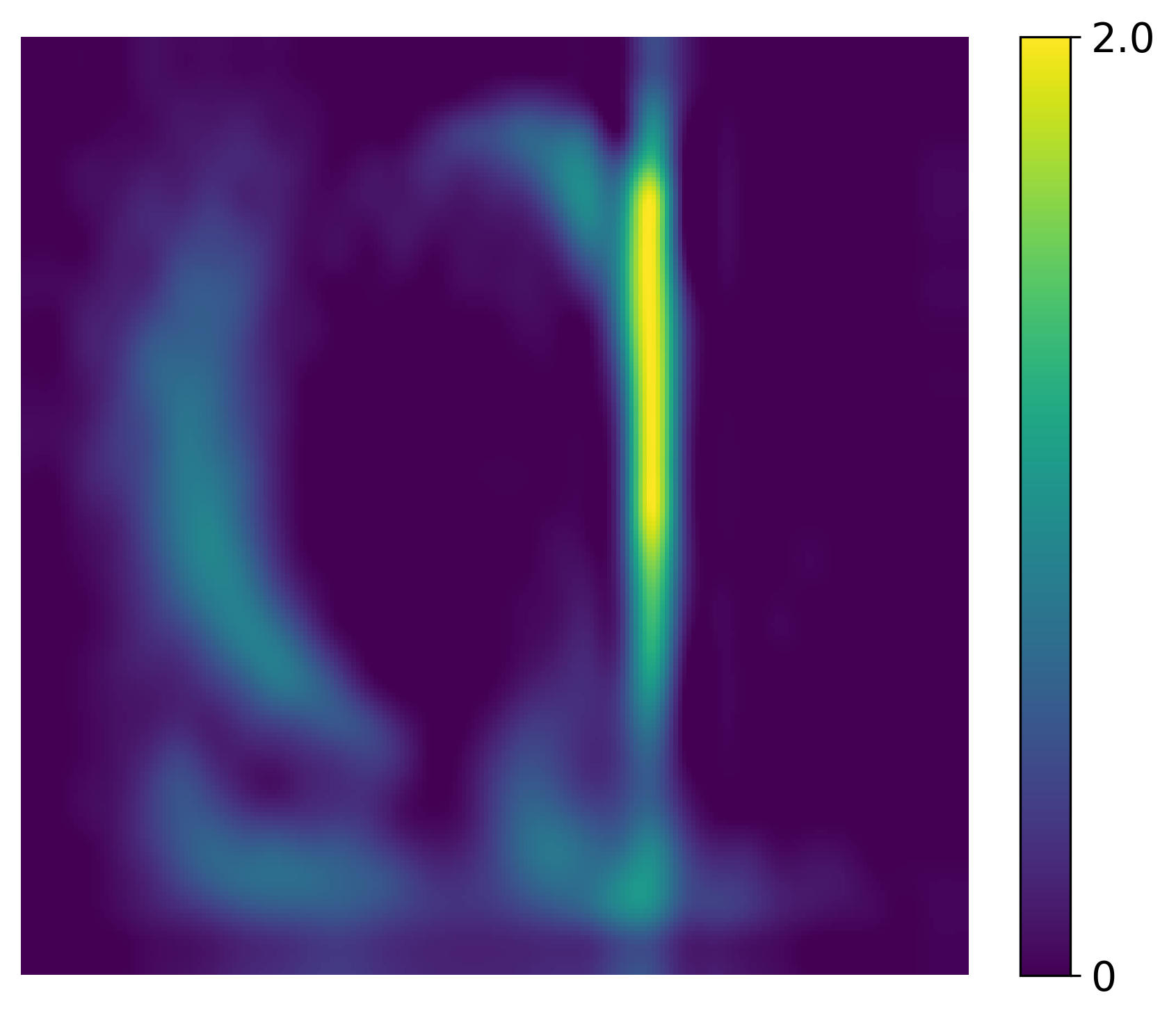}
			\caption{$\varepsilon=0.1$}
			\label{fig:eps_0.3}
		\end{subfigure}
\caption{Reconstruction results illustrating the choice of an optimal value
of the parameter $\protect\varepsilon $ in (\protect\ref{3.2}). Obviously, $%
\protect\varepsilon =0.01$ is the best one out of five values $%
\{0.001,0.01,0.03,0.05,0.1\}$. Hence, we assign the optimal value of this
parameter $\protect\varepsilon =0.01$ in all numerical tests below.}
\label{fig:epsilon_sensitivity_grid}
\end{figure}
\end{test}

As observed in Figure~\ref{fig:epsilon_sensitivity_grid}, at $\varepsilon
=0.001$, the reconstruction suffers from the numerical instability,
resulting in significant distortions. Conversely, as $\varepsilon $
increases beyond $0.03$, the images become increasingly blurred due to the
poor approximation in (\ref{3.2}). At the same time, the value $\varepsilon
=0.01$ yields both the sharpest and the most accurate reconstruction.
Therefore, we set $\varepsilon =0.01$ for the subsequent tests.

In Test 8.3 we investigate the sensitivity of reconstruction results to the
parameter $\lambda $ in the Carleman Weight Function $\varphi _{\lambda
}\left( x_{1}\right) $ in (\ref{4.4}).

\begin{test}
\label{ex:lambda}

We now investigate the choice of the optimal value of the parameter $\lambda 
$ in CWF (\ref{4.4}) by imaging an inclusion with the shape of the letter '$%
A $'. The true coefficient $a(\mathbf{x})=2$ inside of this inclusion and $a(%
\mathbf{x})=0$ outside of it, see (\ref{9.3}) and (\ref{9.4}). Other
parameters are set to the values determined from previous examples: $%
\varepsilon =0.01$, $T=4$, $N_{t}=20,$ $c=5$. To account for the data
perturbations, 1\% noise is added to the measurements, i.e. $\sigma =0.01$
in (\ref{9.10}), (\ref{9.11}). The Tikhonov regularization parameter is set
to $\alpha =3\times 10^{-5}$, and reconstructions are performed for five
values of the weight: $\lambda \in \{1,2,3,4,5\}$.

\begin{figure}[th]
\centering
\begin{subfigure}[b]{0.32\textwidth}
			\includegraphics[width=\textwidth]{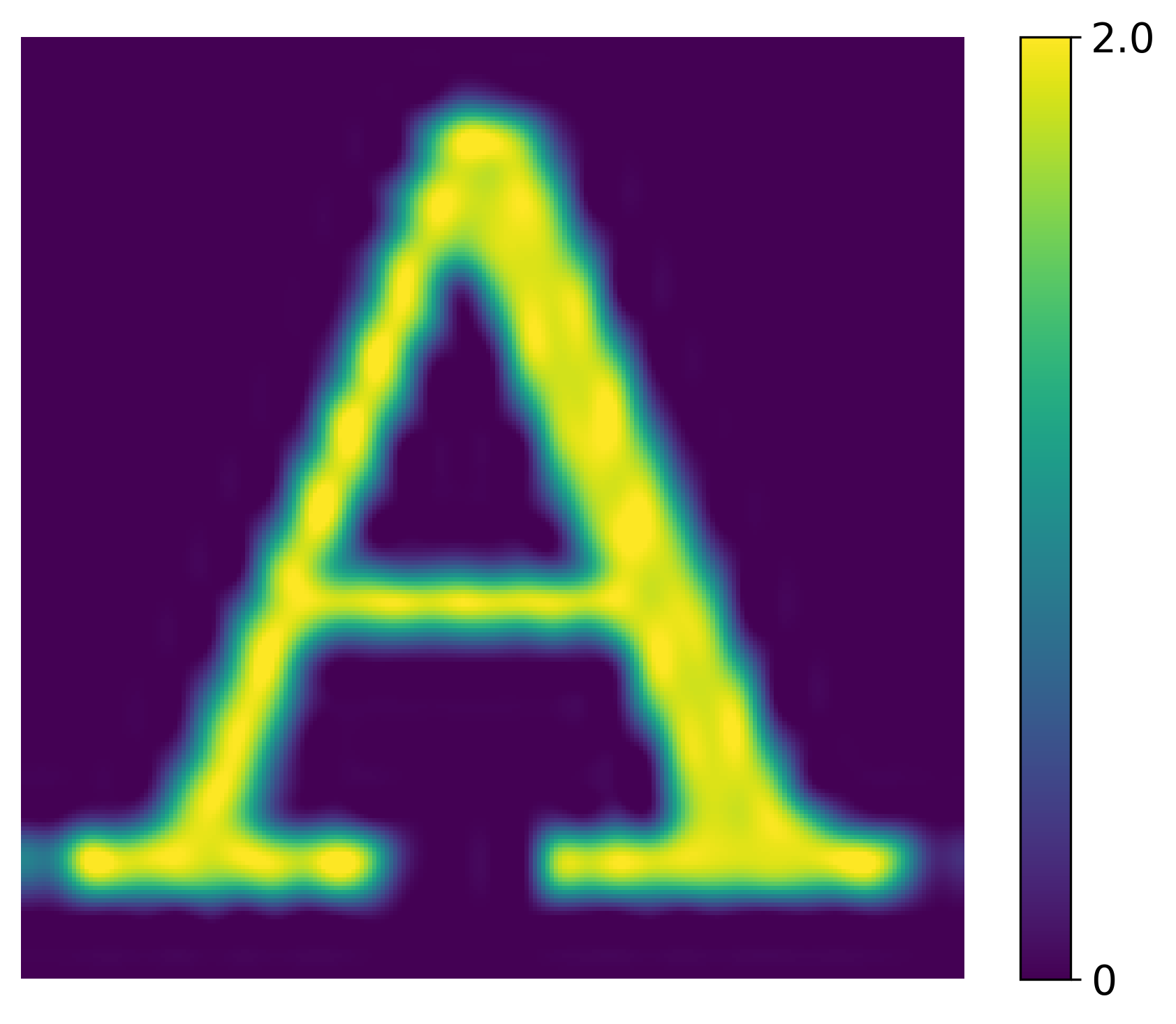}
			\caption{True image}
		\end{subfigure}
\hfill 
\begin{subfigure}[b]{0.32\textwidth}
			\includegraphics[width=\textwidth]{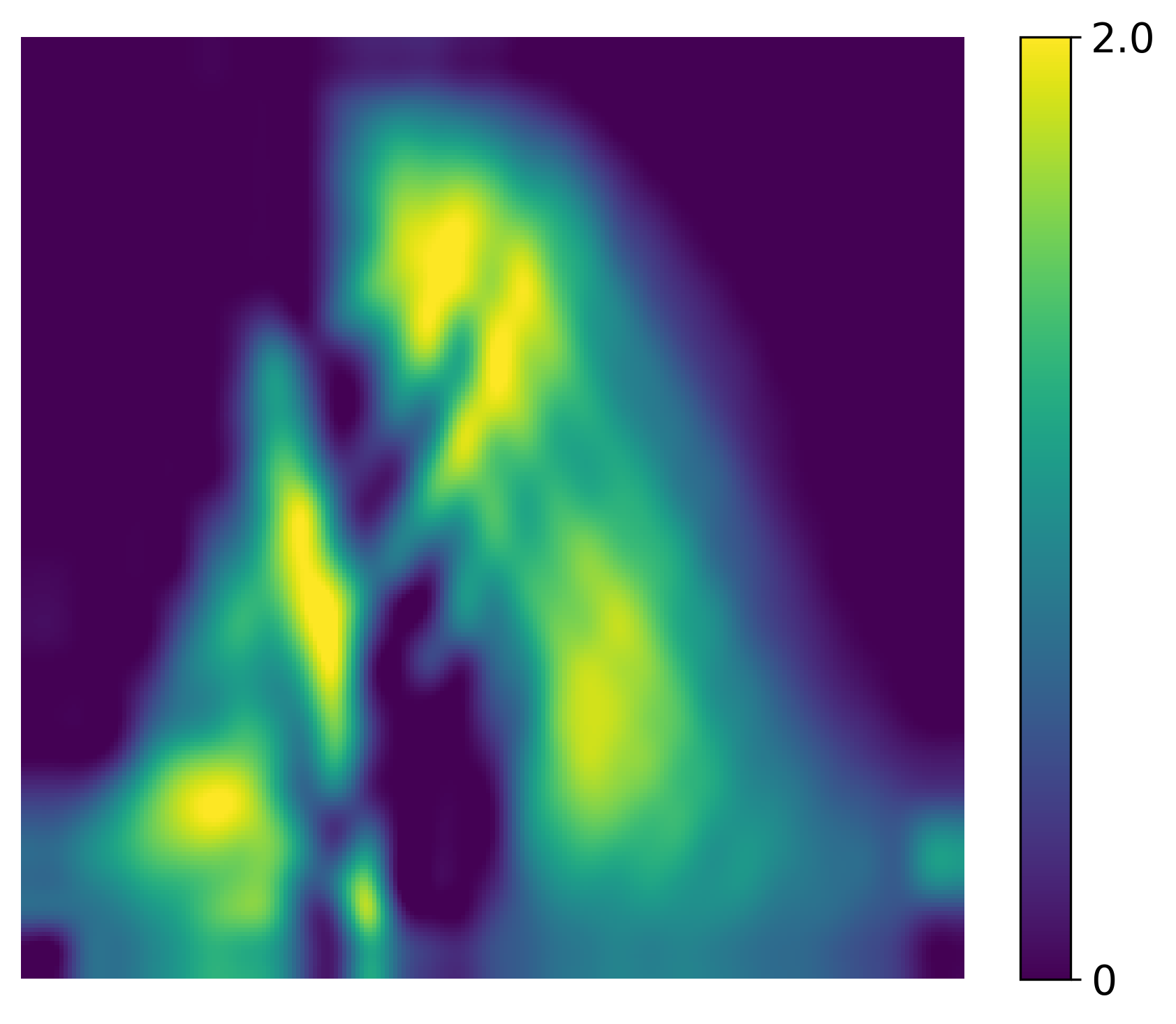}
			\caption{$\lambda=1$}
		\end{subfigure}
\hfill 
\begin{subfigure}[b]{0.32\textwidth}
			\includegraphics[width=\textwidth]{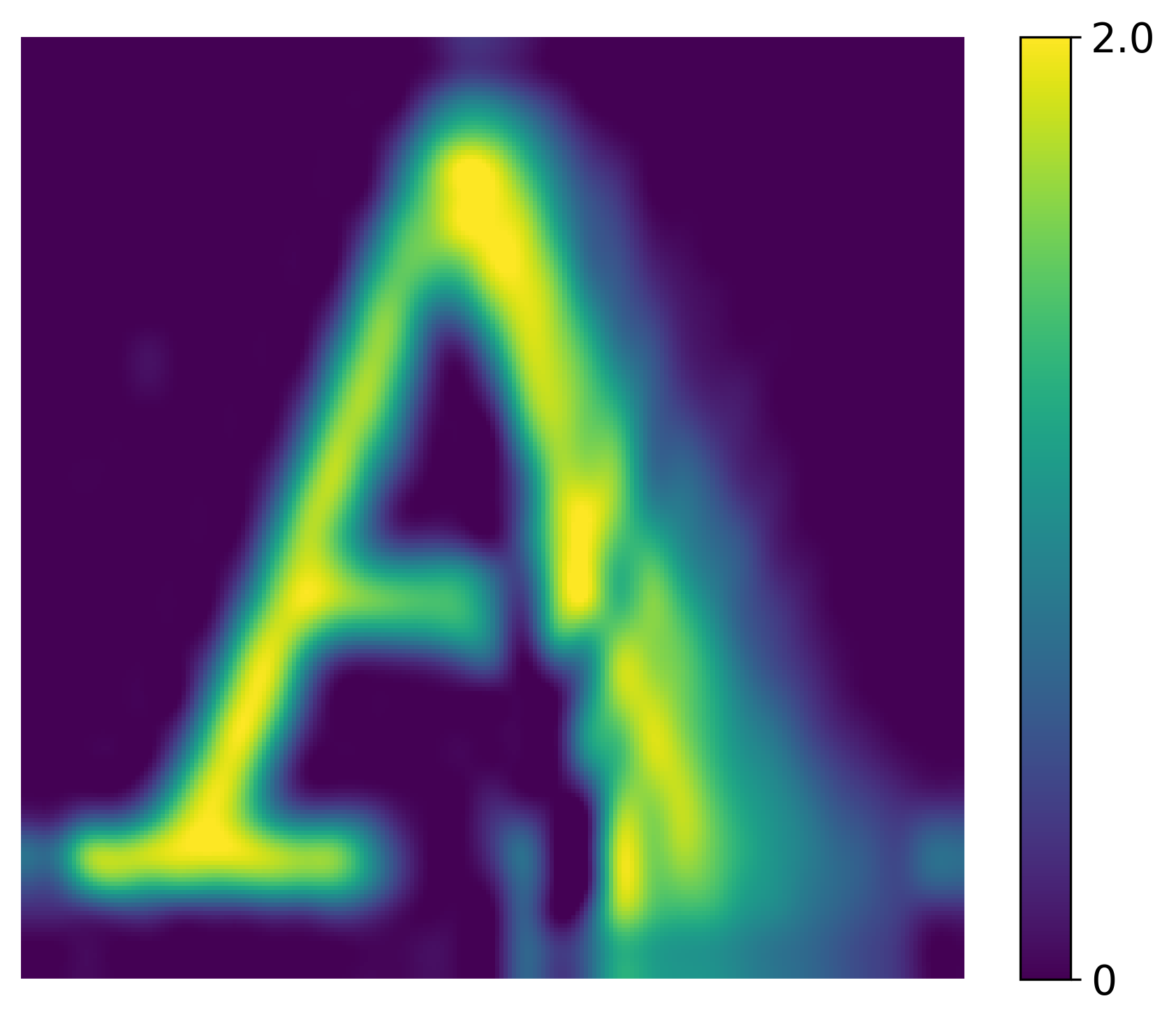}
			\caption{$\lambda=2$}
		\end{subfigure}
\par
\vspace{2mm}
\par
\phantom{a}\hfill 
\begin{subfigure}[b]{0.32\textwidth}
			\includegraphics[width=\textwidth]{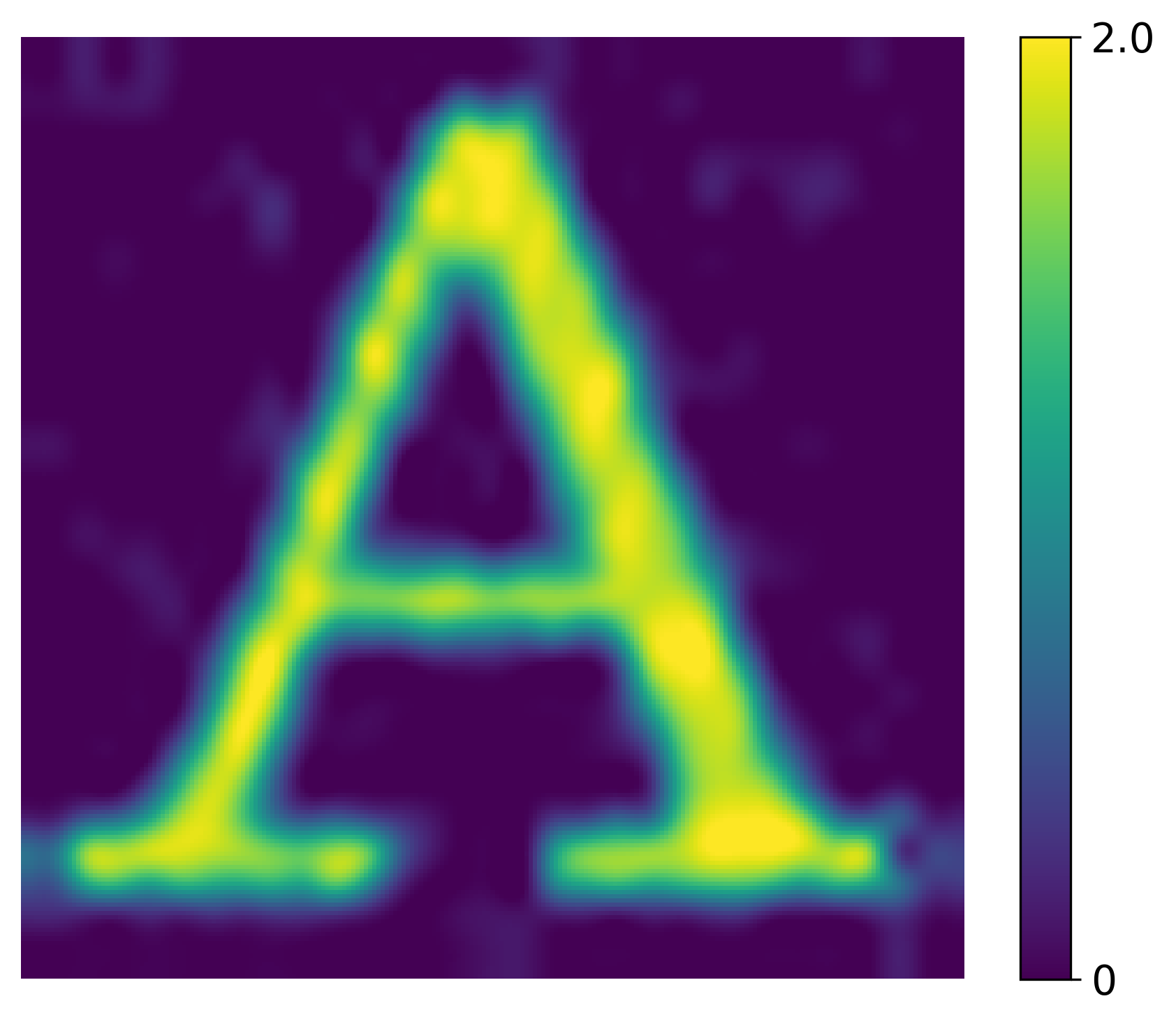}
			\caption{$\lambda=3$}
		\end{subfigure}
\hfill 
\begin{subfigure}[b]{0.32\textwidth}
			\includegraphics[width=\textwidth]{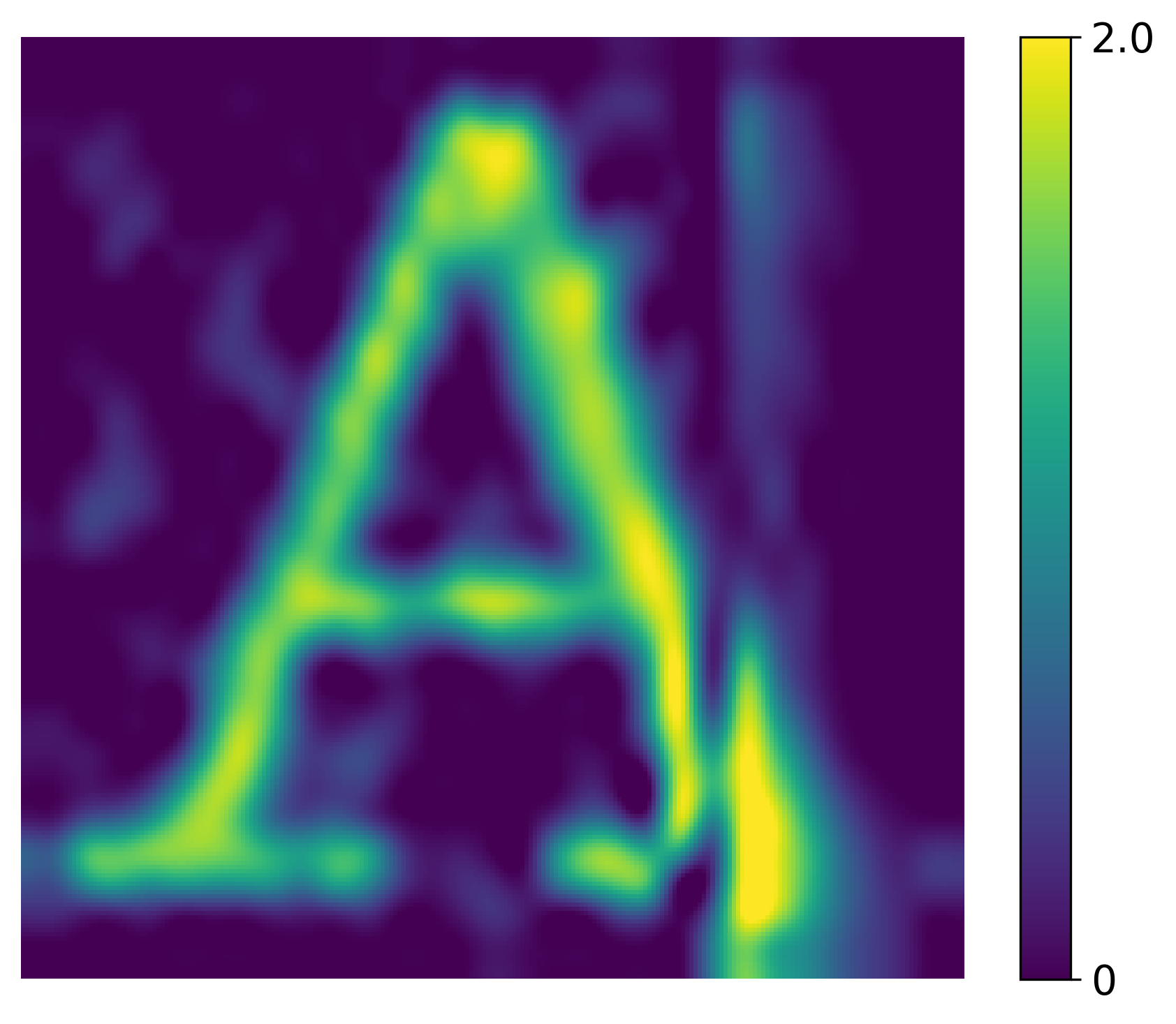}
			\caption{$\lambda=4$}
		\end{subfigure}
\hfill 
\begin{subfigure}[b]{0.32\textwidth}
			\includegraphics[width=\textwidth]{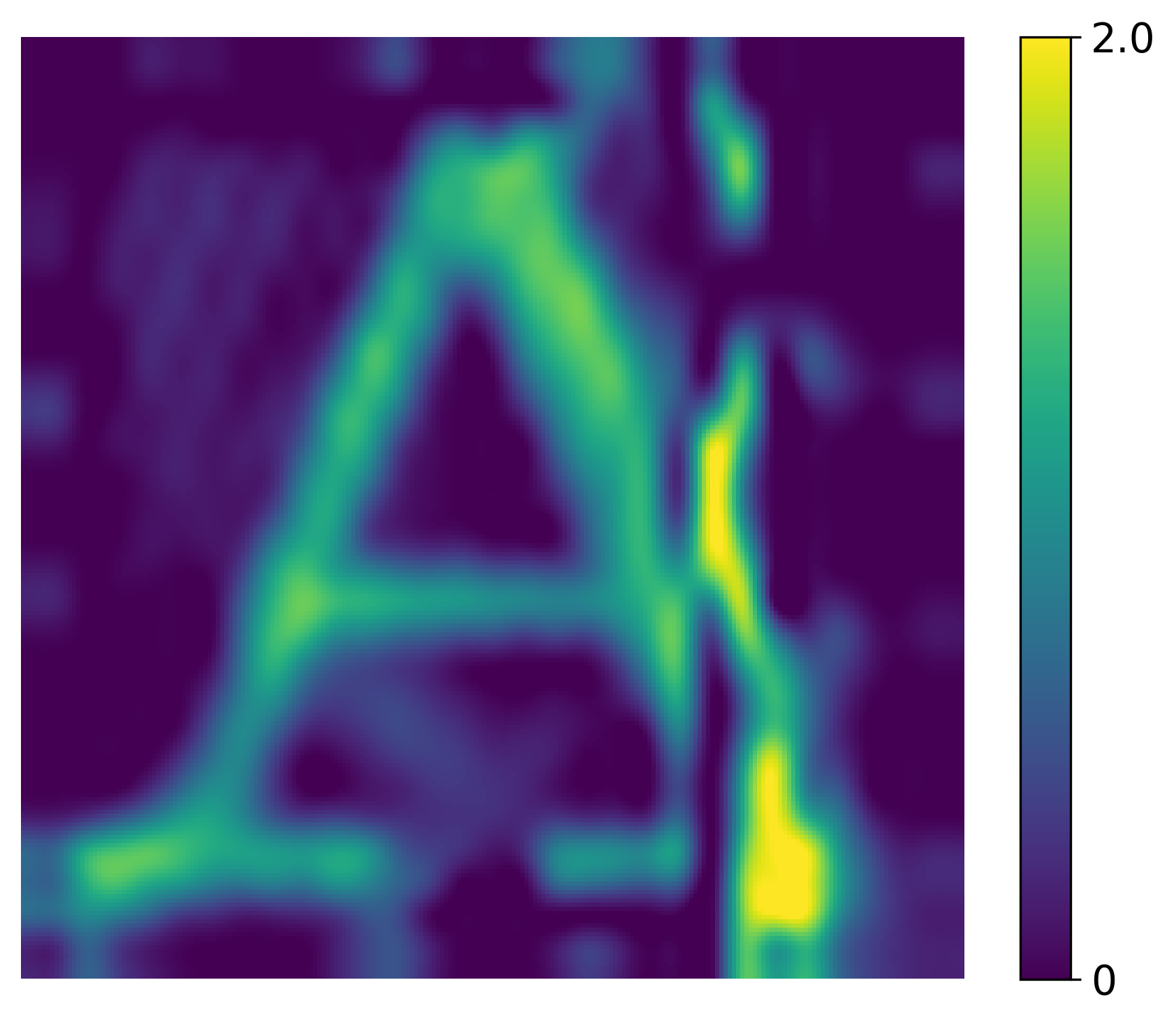}
			\caption{$\lambda=5$}
		\end{subfigure}
\par
\vspace{2mm}
\caption{Reconstruction results illustrating the choice of an optimal value
of the parameter $\protect\lambda $ in Carleman Weight Function $\protect%
\varphi _{\protect\lambda }\left( x_{1}\right) $ in (\protect\ref{4.4}).
Obviously, $\protect\lambda =3$ is the best one out of five values $%
\{1,2,3,4,5\}$. Hence, we assign the optimal value of the parameter $\protect%
\lambda =3$ in all numerical tests below.}
\label{fig:lambda_comparison}
\end{figure}
\end{test}

As shown in Figure~\ref{fig:lambda_comparison}, the reconstruction quality
varies significantly with the value of $\lambda $. For smaller values ($%
\lambda =1,2$), the recovered images appear blurred, and the boundaries of
the letter `$A$' are poorly defined. In contrast, for larger values $\lambda
=4,5$, the images exhibit noticeable artifacts and structural distortions;
specifically, the limbs of the letter `$A$' appear to be either disconnected
or broken. The value $\lambda =3$ yields the clearest result, providing a
sharp and continuous recovery of the inclusion shape with minimal background
noise. Based on this observation, we fix $\lambda =3$ for the remainder of
this study.

Based on Tests 8.1-8.3, we use the following values of the parameters in all
tests below:%
\begin{equation}  \label{9.12}
T=4,\text{ }N_{t}=20,\text{ }\varepsilon =0.01,\text{ }\alpha =3\times
10^{-5},\text{ }\lambda =3,\text{ }c=5.
\end{equation}

In the next two tests, we investigate the robustness of our method with
respect to the different noise levels and coefficient values.

\begin{test}
\label{ex:noise}

We investigate the impact of the noise level $\sigma $ in the boundary data (%
\ref{9.10}), (\ref{9.11}) via imaging an inclusion of the shape of the
letters '$SZ$'. The true coefficient is $a(\mathbf{x})=2$ inside of this
inclusion and $a(\mathbf{x})=0$ outside of it, see (\ref{9.3}) and (\ref{9.4}%
). The choice of parameters is as in (\ref{9.12}). Noise levels of $1\%$, $%
3\%$, $5\%$ are added to the measurements, i.e. $\sigma =0.01,$ $0.03,$ $%
0.05 $.
\end{test}

\begin{figure}[tbph]
\centering
\subcaptionbox{True
image}{\includegraphics[width=0.23\textwidth]{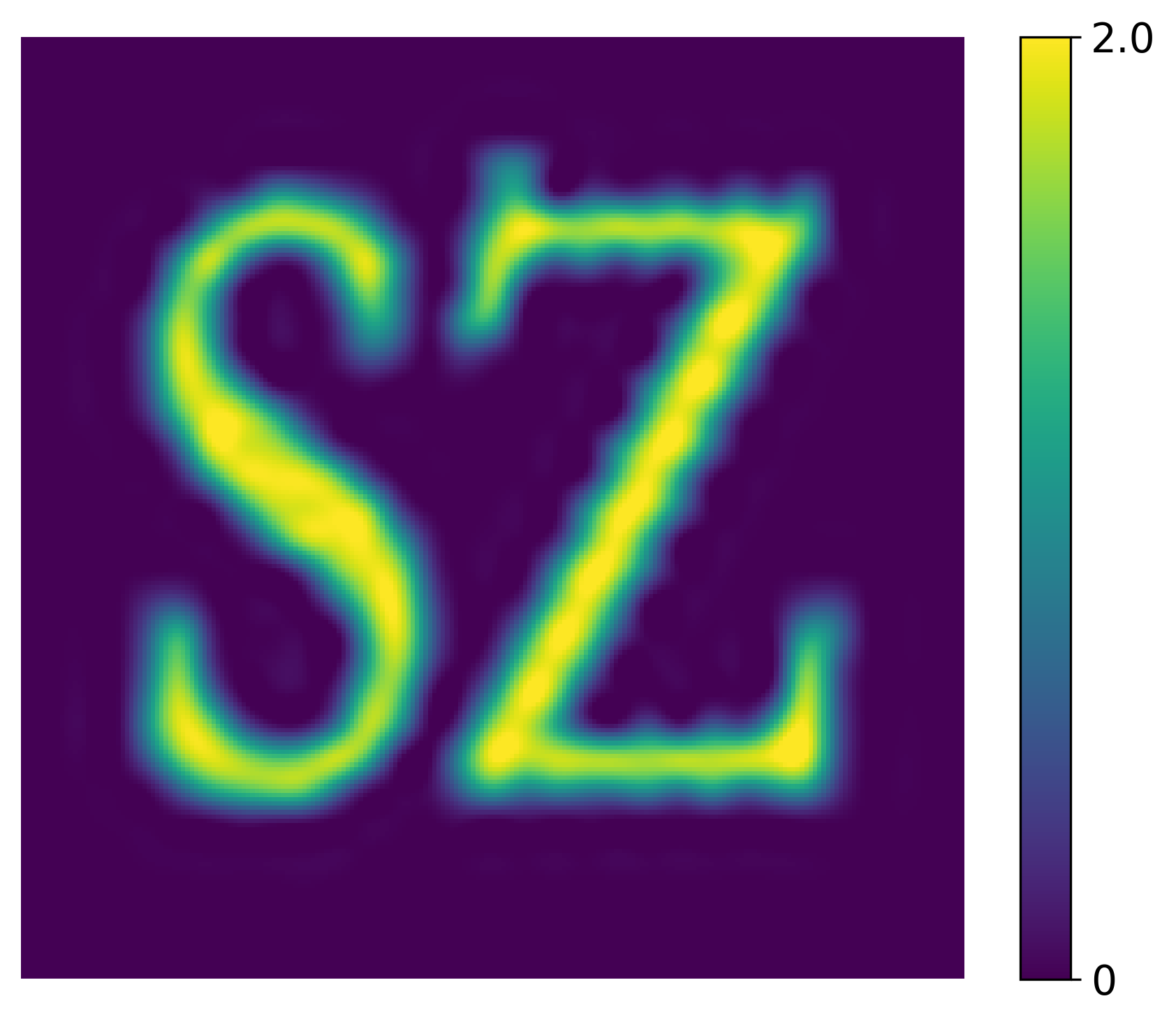}}%
\hfill 
\subcaptionbox{1\%
noise}{\includegraphics[width=0.23\textwidth]{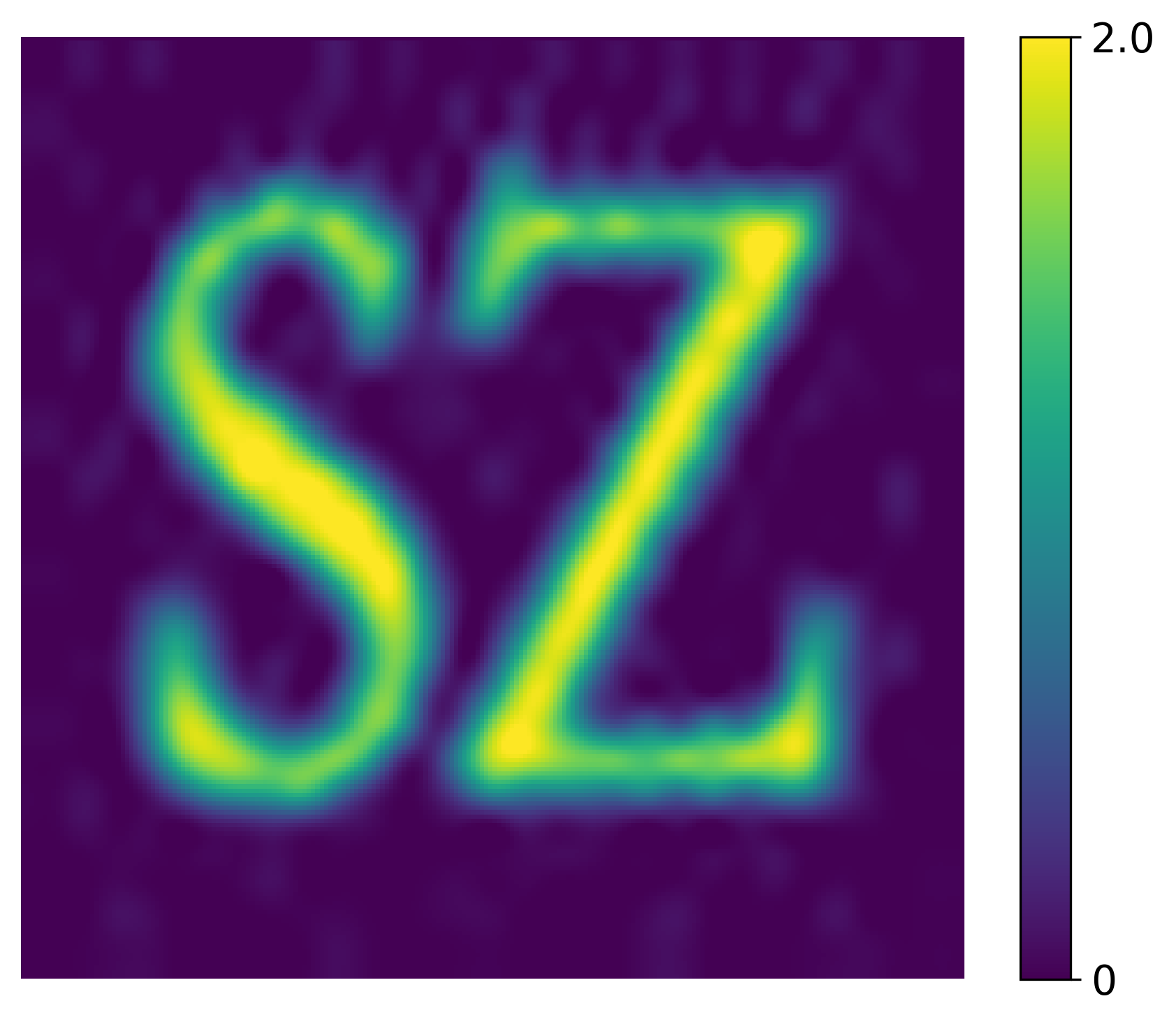}}%
\hfill 
\subcaptionbox{3\%
noise}{\includegraphics[width=0.23\textwidth]{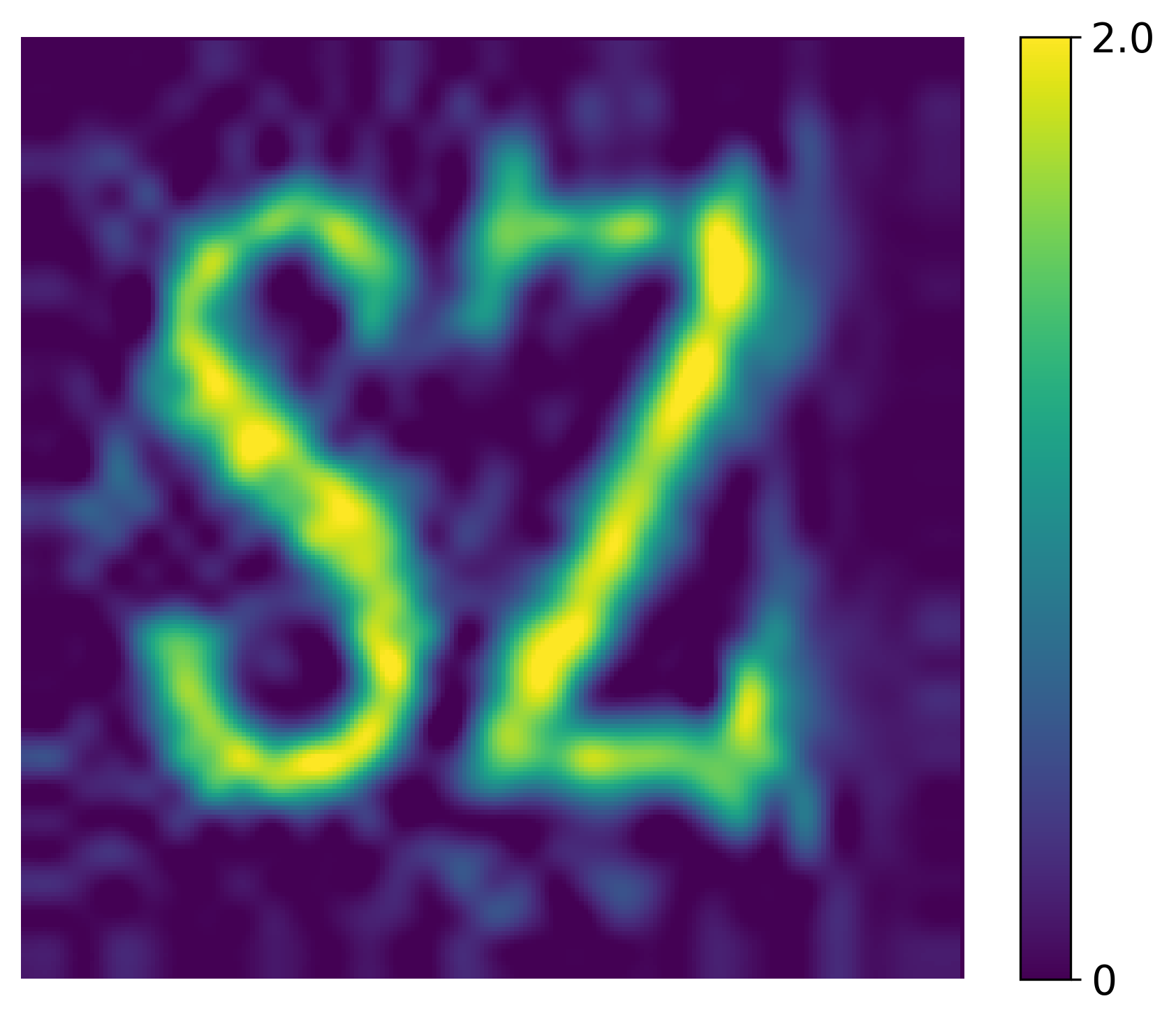}}%
\hfill 
\subcaptionbox{5\%
noise}{\includegraphics[width=0.23\textwidth]{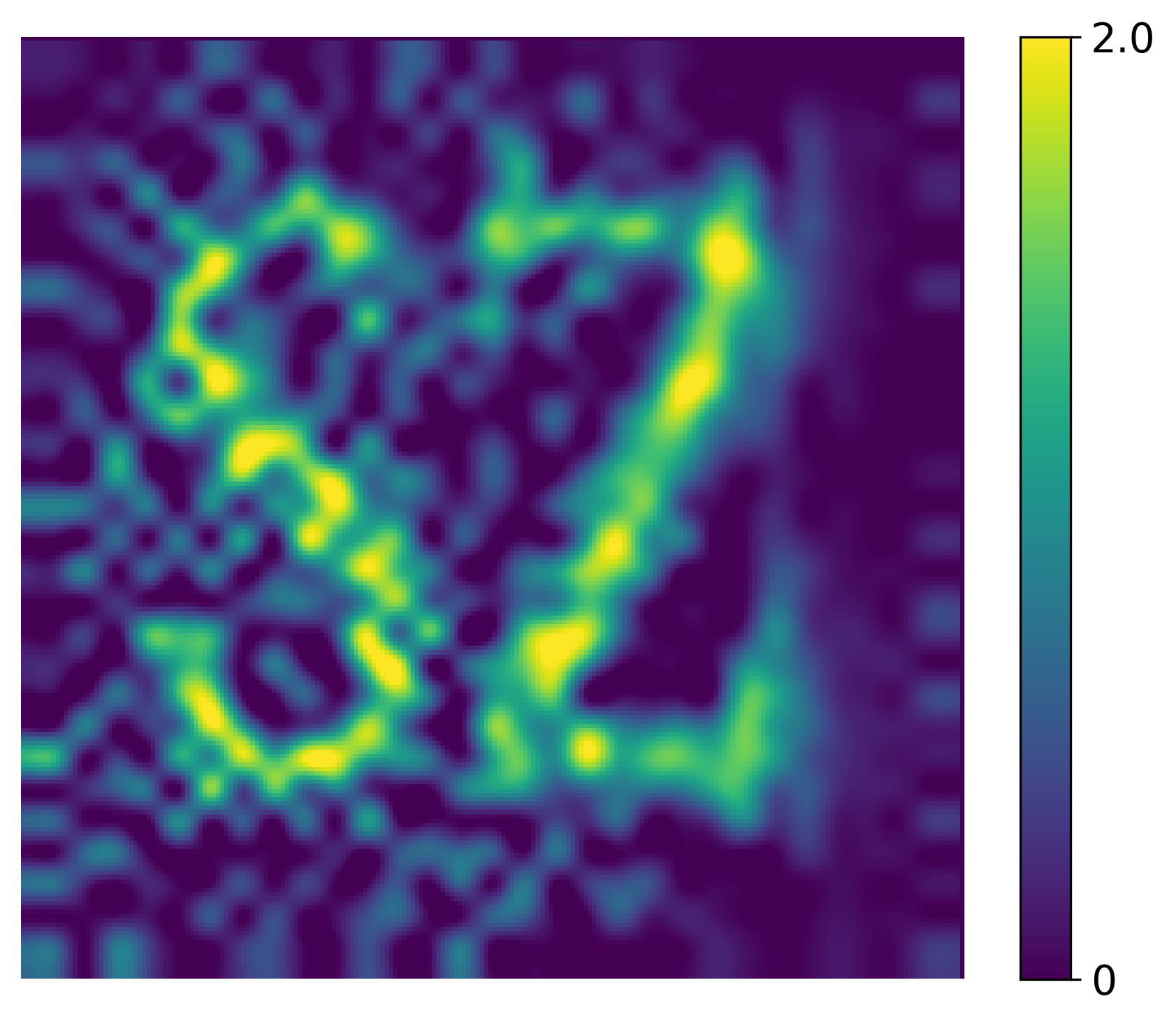}}
\caption{Reconstruction results illustrating the sensitivity to the noise
level. Even for the highest noise levels of 5\% shapes of both inclusions
are recognizable and the maximal value $\max a\left( \mathbf{x}\right) =2$
inside these letters is reconstructed accurately.}
\label{fig:noise_levels}
\end{figure}

The reconstruction results are presented in Figure~\ref{fig:noise_levels}.
As the noise level increases, we observe a corresponding increase in
background artifacts and blurring. Nevertheless, even at the highest noise
level of $5\%$, the geometry of the letters `$SZ$' remains recognizable and
the maximal value $\max a\left( \mathbf{x}\right) =2$ inside these letters
is reconstructed accurately. This demonstrates a good degree of robustness
of the method with respect to the random noise in the input data (\ref{9.10}%
), (\ref{9.11}).

Using (\ref{9.3}) and (\ref{9.4}), we now examine the influence of the
magnitude $a$ of the value of the coefficient $a\left( \mathbf{x}\right) $
inside the inclusion on the quality of the reconstruction.

\begin{test}
\label{ex:coeff_value}

This test investigates the performance of our method for different values of
the true coefficient $a(\mathbf{x})$. As in (\ref{9.3}), we set the true
coefficient to be constant within the inclusion and zero outside. The values
of the constant $a$ within inclusions are taken as in (\ref{9.4}), i.e. $a=2$%
, $3$, $5$,$10$. The shapes of inclusions are the letter '$A$' and the
letters '$SZ$'. All other parameters remain the same as the ones in (\ref%
{9.12}). In all cases, $1\%$ noise is added to the measurements. The
reconstructions are compared to assess the sensitivity of the method to the
magnitude of the unknown coefficient.
\end{test}

The results of the reconstructions are presented in Figure~\ref%
{fig:coeff_values}. The method successfully reconstructs the shapes of both
inclusions for all tested coefficient values. Furthermore, the maximal
values of the coefficient $a(\mathbf{x})$ are accurately reconstructed. It
is noteworthy that the reconstruction quality does not degrade for larger
values of $a(\mathbf{x})$. Even for the high value of $a=10$, the images
remain sharp and accurate, confirming the capability of the convexification
method to handle strongly nonlinear problems with large values of
coefficients.

Finally, we present the results for the three-dimensional case.

\begin{test}
We consider the reconstruction of $3-$Dimensional inclusions of the shapes
of the letters `$L$' and `$K$' located inside the unit cube $\Omega
=[1,2]^{3}$. The true coefficient is $a(\mathbf{x})=2$ inside the inclusions
and $a(\mathbf{x})=0$ elsewhere, see (\ref{9.3}) and (\ref{9.4}). The
spatial domain is discretized using a $20\times 20\times 20$ grid. We set
parameters as in (\ref{9.12}). We add $1\%$ random noise to the boundary
measurement data.

The reconstruction results are displayed in two separate figures for
clarity. Figure~\ref{fig:3d_results} presents the 3-Dimensional isosurface
plots of the reconstructed coefficient $a(\mathbf{x})$, illustrating the
spatial recovery of the inclusions. Figure~\ref{fig:3d_slice} provides the
corresponding X-Z cross-sectional views (front view) to demonstrate the
accuracy in the vertical plane, maintaining consistent color mapping between
the $2-$d cross-sections and $3-$d isosurfaces.

\begin{figure}[htbp]
\centering
\begin{subfigure}[b]{0.45\textwidth}
			\centering
			\includegraphics[width=\textwidth]{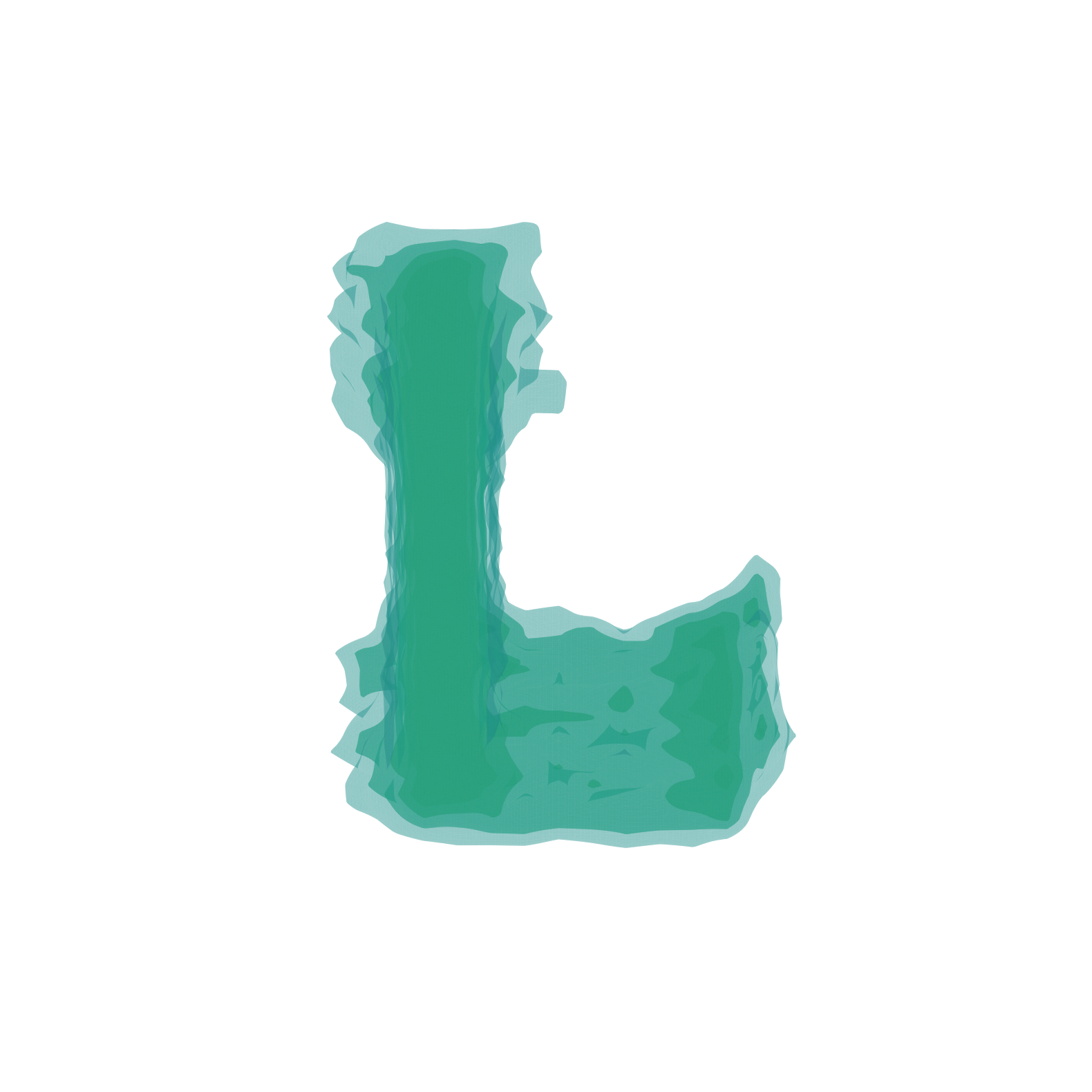} 
			\caption{True 3-Dimensional `L'}
		\end{subfigure}
\hfill 
\begin{subfigure}[b]{0.45\textwidth}
			\centering
			\includegraphics[width=\textwidth]{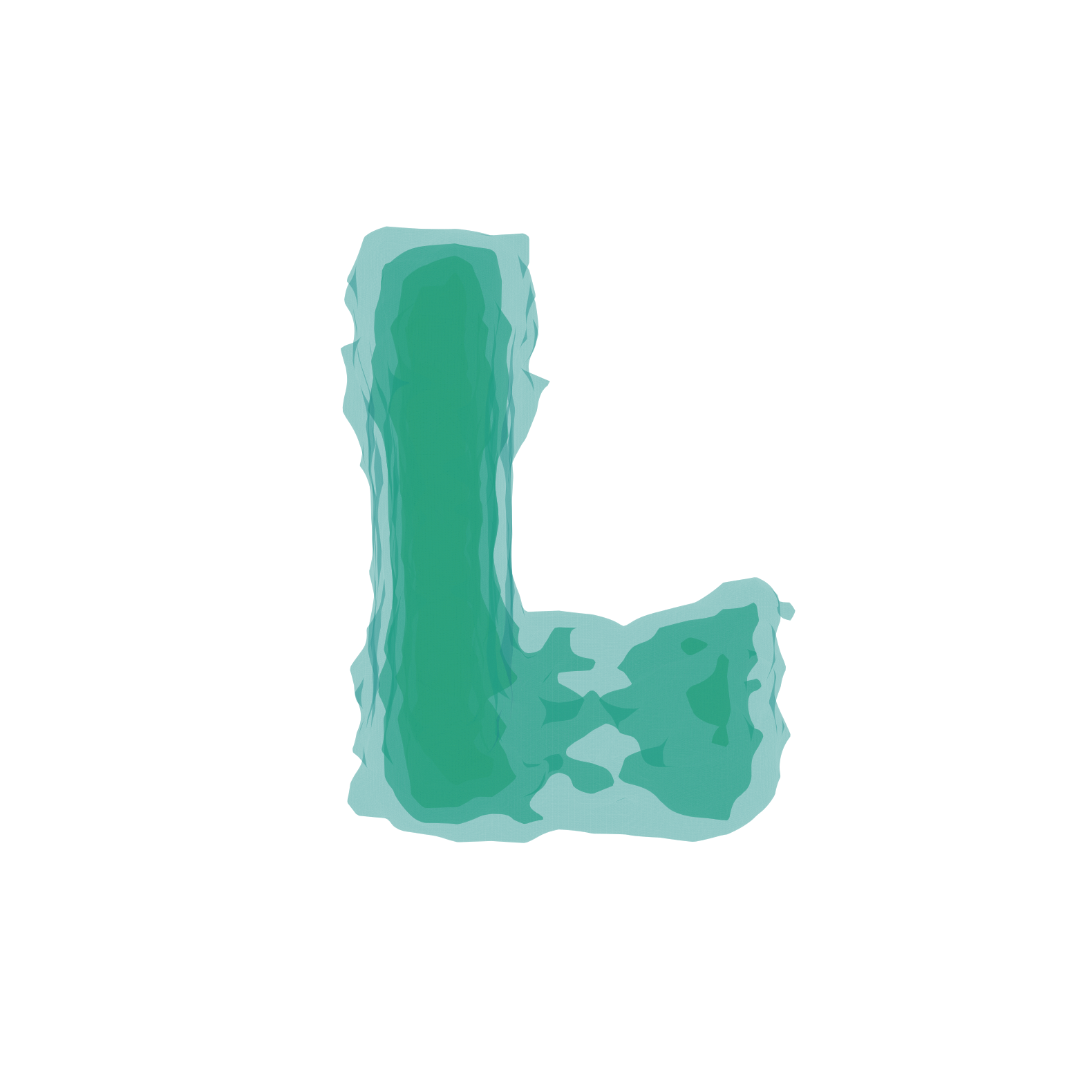}
			\caption{Reconstructed 3-Dimensional `L'}
		\end{subfigure}
\par
\vspace{1em} 
\begin{subfigure}[b]{0.45\textwidth}
			\centering
			\includegraphics[width=\textwidth]{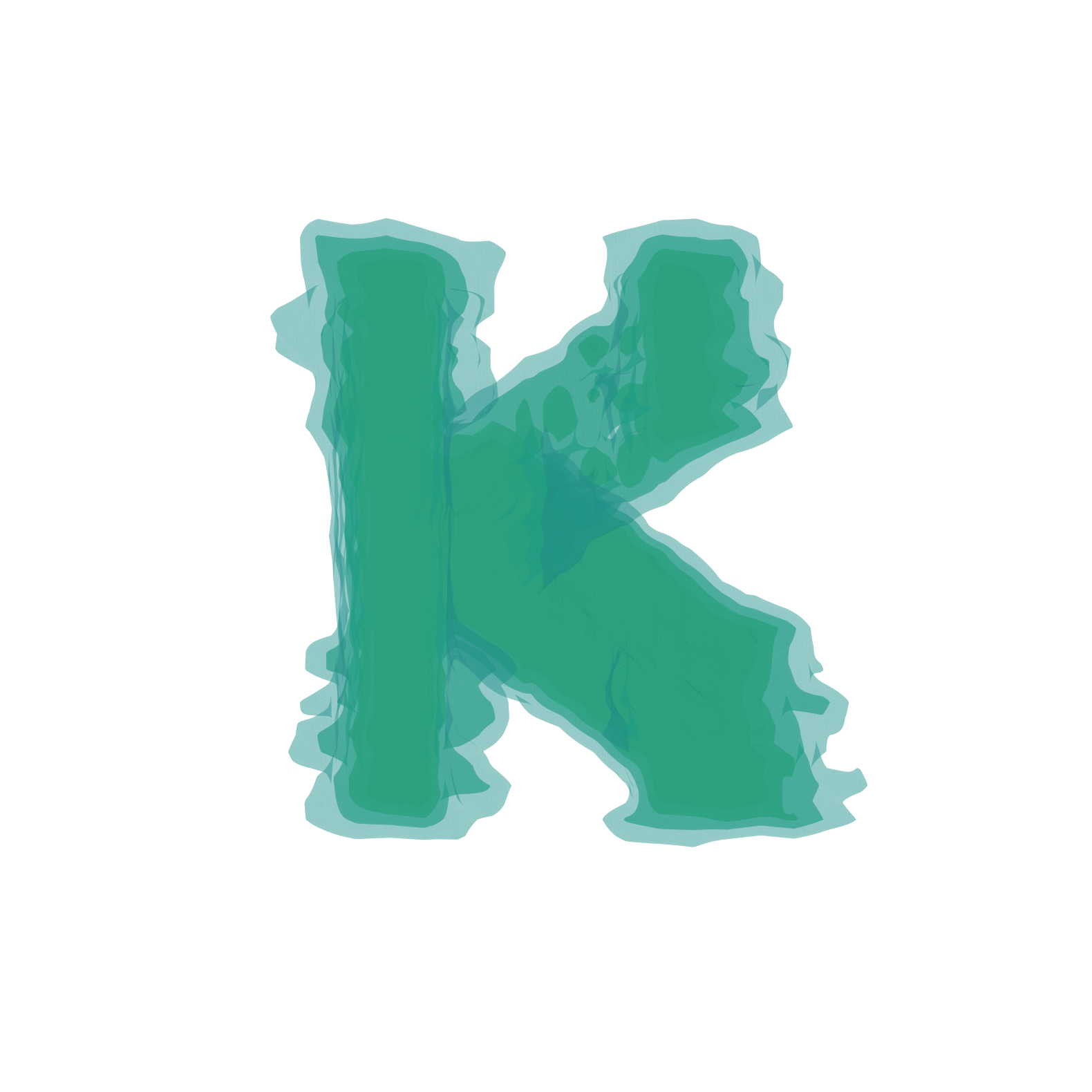} 
			\caption{True 3-Dimensional `K'}
		\end{subfigure}
\hfill 
\begin{subfigure}[b]{0.45\textwidth}
			\centering
			\includegraphics[width=\textwidth]{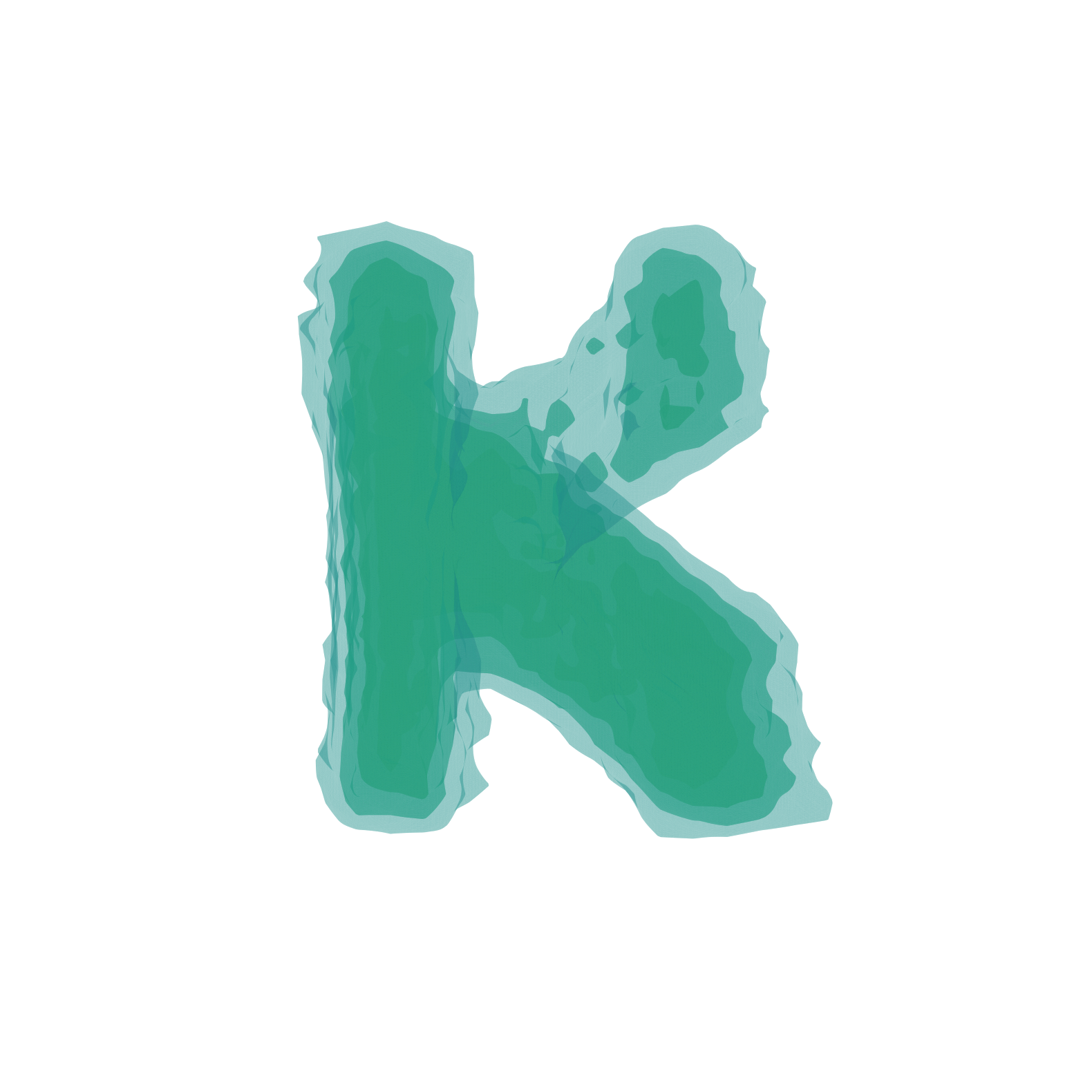} 
			\caption{Reconstructed 3-Dimensional `K'}
		\end{subfigure}
\caption{Three-dimensional reconstruction results for different inclusion
shapes (`L' and `K'). Here $a=2$, see (\protect\ref{7.3}).}
\label{fig:3d_results}
\end{figure}

\begin{figure}[tbph]
\centering
\begin{subfigure}[b]{0.45\textwidth}
			\centering
			\includegraphics[width=\textwidth]{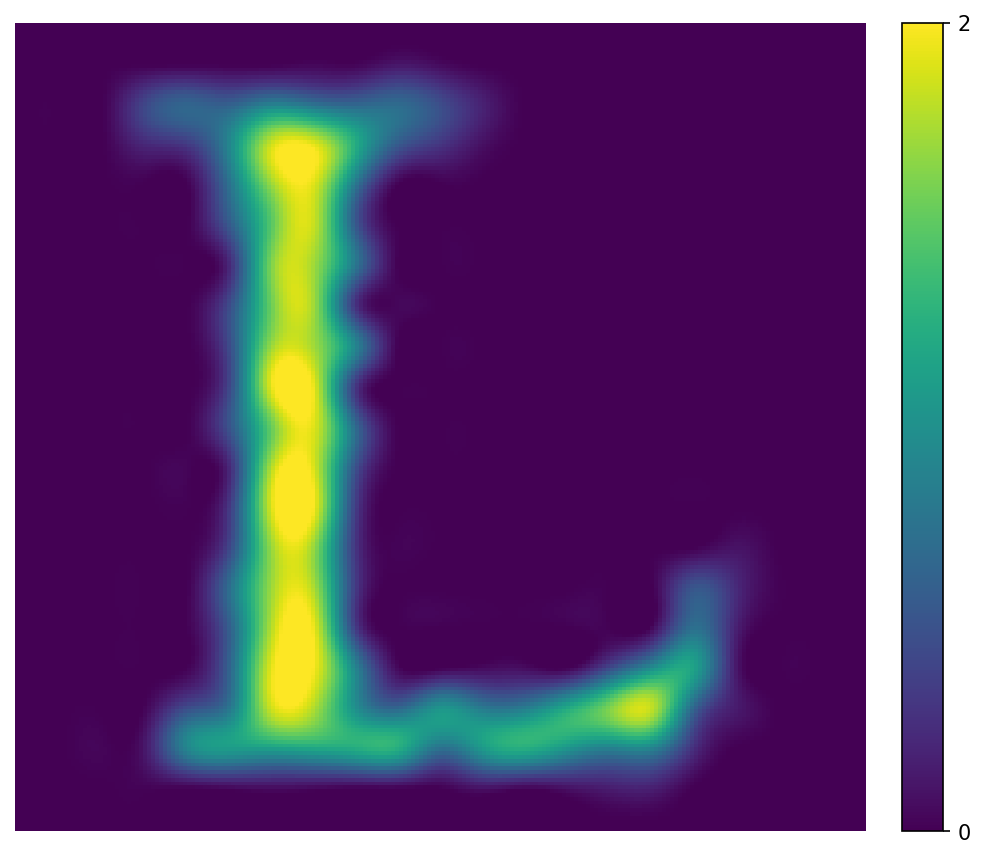} 
			\caption{True `L' in X-Z section}
		\end{subfigure}
\hfill 
\begin{subfigure}[b]{0.45\textwidth}
			\centering
			\includegraphics[width=\textwidth]{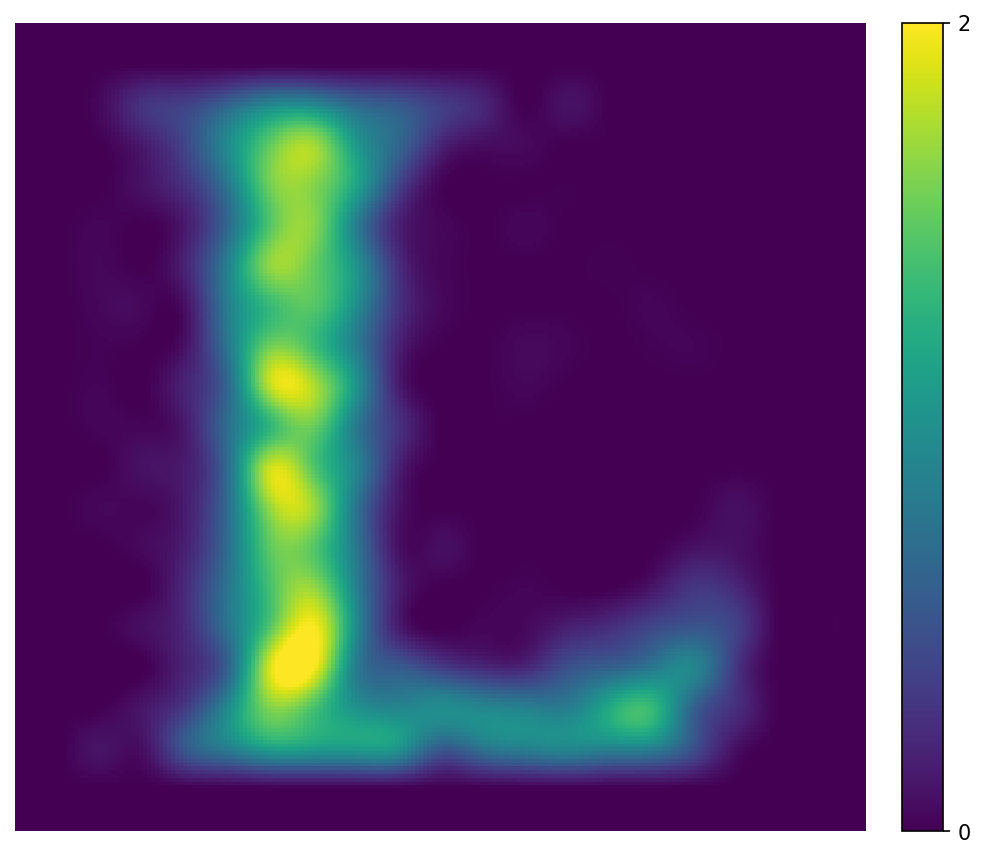} 
			\caption{Reconstructed `L' in X-Z section}
		\end{subfigure}
\par
\vspace{1em} 
\begin{subfigure}[b]{0.45\textwidth}
			\centering
			\includegraphics[width=\textwidth]{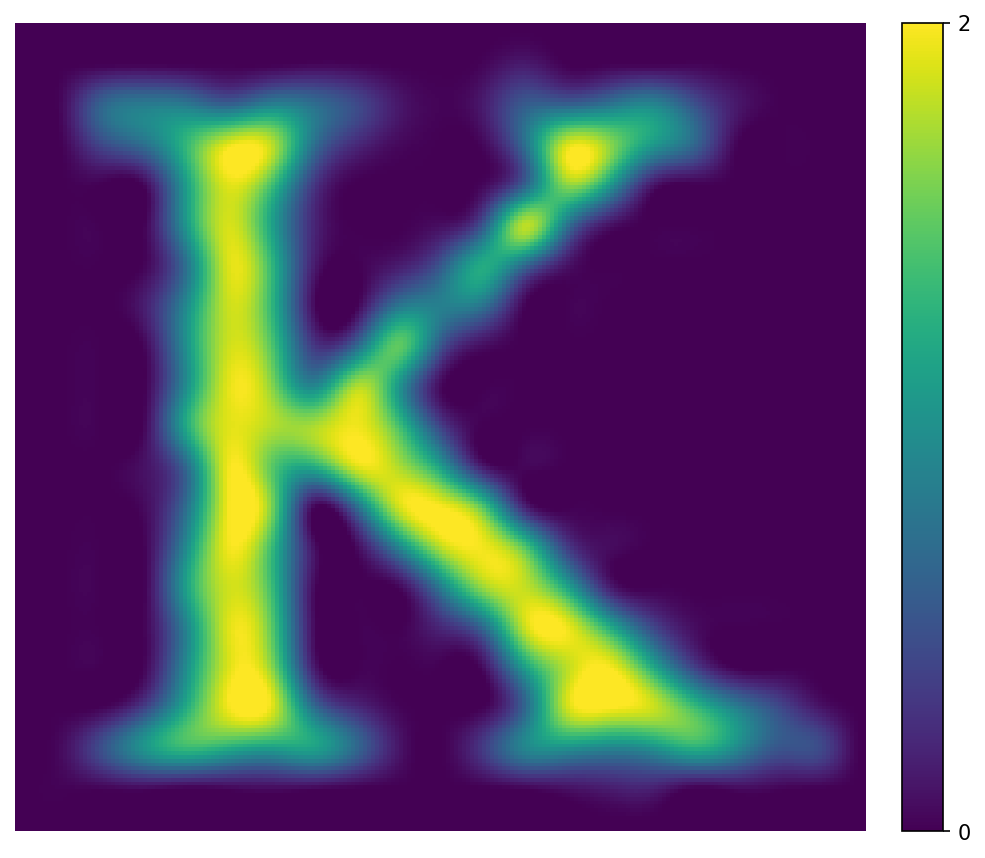}
			\caption{True `K' in X-Z section}
		\end{subfigure}
\hfill 
\begin{subfigure}[b]{0.45\textwidth}
			\centering
			\includegraphics[width=\textwidth]{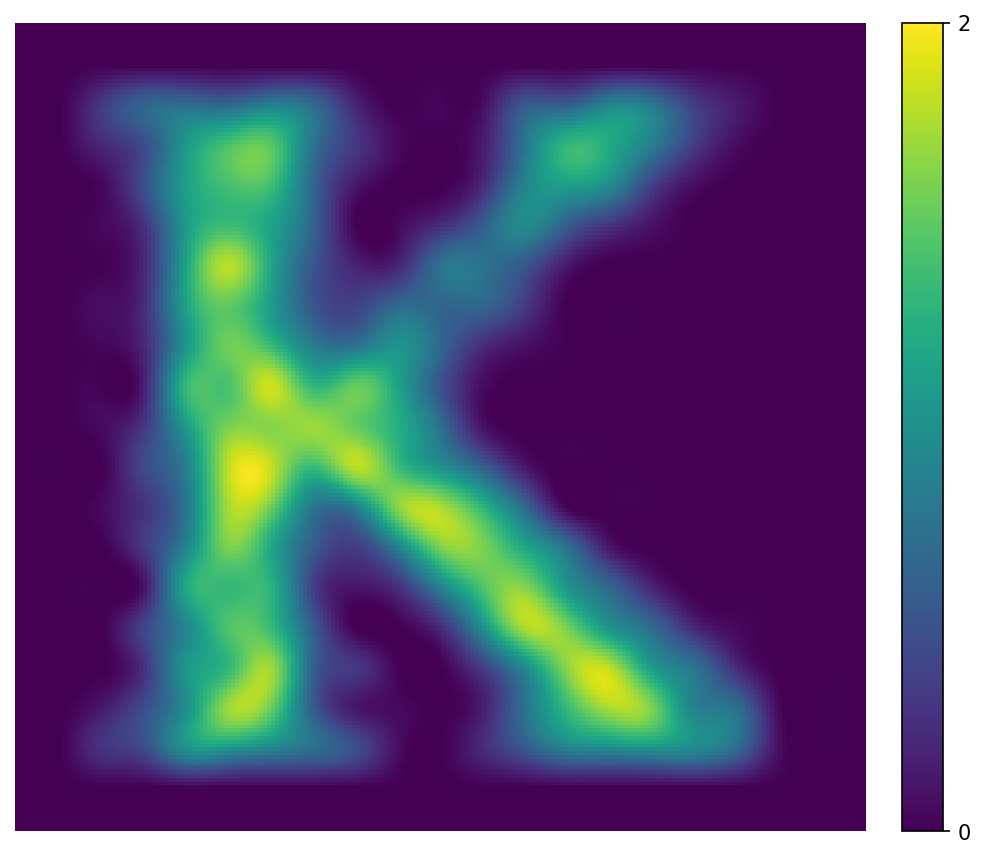} 
			\caption{Reconstructed `K' in X-Z section}
		\end{subfigure}
\caption{Three-dimensional reconstruction results: X-Z cross-sectional views
corresponding to Figure~\protect\ref{fig:3d_results}. Here $a=2$, see (%
\protect\ref{7.3}).}
\label{fig:3d_slice}
\end{figure}
\end{test}

The reconstruction results are presented in Figures~\ref{fig:3d_results} and~%
\ref{fig:3d_slice}. As observed, the algorithm successfully reconstructs the
spatial structures of both the `$L$' and `$K$' shapes with high fidelity.
The value of the number $a$ inside the inclusions is also accurately
reconstructed. The consistent accurate recoveries of these distinct 3-d
geometries strongly validates the robustness of the proposed method.

\section{Summary}

\label{sec:9}

For the first time, we have developed a globally convergent numerical method
for the $n$-Dimensional CIP for the wave-like PDE (\ref{2.9}) with the
unknown coefficient $a\left( \mathbf{x}\right) $ for the most challenging
case when the $\delta \left( \mathbf{x}\right) -$function is present in one
of the initial conditions (\ref{2.10}) with a single location of the point
source at $\left\{ 0\right\} $. First, we have applied an analog (\ref{2.12}%
) of the Laplace transform to transform the original CIP1 to the CIP2 for a
parabolic equation. Next, we have proposed a new approximate mathematical
model for CIP2.

We have developed a new version of the globally convergent convexification
numerical method for our approximate mathematical model. Global convergence
here is understood in terms of Definition 1.1. Furthermore, uniqueness
theorem is proven for that model.

We have carried out extensive numerical studies of our method both in
2-Dimensional and 3-Dimensional cases. Our computational results revealed a
very good accuracy of our reconstructions of complicated structures. We
conclude, therefore, that this reconstruction accuracy confirms a high
degree of the adequacy of our approximate mathematical model.

\section{Appendix A: Proof of Theorem 2.1}

\label{sec:10}

We use here the decomposition of the fundamental solution of problem (\ref%
{2.9}), (\ref{2.10}), which is given in \cite[Lemmata 2.2.1 and 2.2.3]{Rom1}%
. In this book the representation of the solution is in the form of infinite
series, assuming that coefficients of the corresponding hyperbolic equation
belong $C^{\infty }(\mathbb{R}^{n})$. Unlike \cite{Rom1}, we assume here
that the coefficient $a(\mathbf{x})\in C^{l}\left( \mathbb{R}^{n}\right) ,$ $%
\ell =5[n/2]+3,$ see (\ref{2.7}). Hence, we use only a finite segment of
this series with a remainder term. This representation of the solution is
given below for odd and even $n$ separately.

We start with the representation (\ref{11}) for $n=2m+1$, $m\geq 1$. In this
case formulae (\ref{12}) for functions $\alpha _{k}(\mathbf{x})$ are given
in \cite[page 32]{Rom1}. It follows from formulae (\ref{12}) that $\alpha
_{s}\in C^{\ell -2(s+m)}(\mathbb{R}^{n})$, $s=\left[ -m,S\right] ,$ and all $%
\alpha _{k}(\mathbf{x})$ are bounded in $\mathbb{R}^{n}$, since 
\begin{equation*}
\ell =5\left[ \frac{n}{2}\right] +3>2(S+m)=2\left( 2m+1\right) =2n+2.
\end{equation*}

Substituting representation (\ref{11}) in equation (\ref{2.9}) and takin
into account that $\alpha _{k}(\mathbf{x})$ are solutions of equations%
\begin{equation*}
\left. 
\begin{array}{c}
2\nabla \alpha _{-m}(\mathbf{x})\cdot \nabla \left( |\mathbf{x}|^{2}\right)
=0, \\ 
2\nabla \alpha _{s}(\mathbf{x})\cdot \nabla \left( |\mathbf{x}|^{2}\right)
+4(s-m)\alpha _{s}(\mathbf{x})=\Delta \alpha _{s-1}(\mathbf{x})+a(\mathbf{x}
)\alpha _{s-1}(\mathbf{x}),~s\in \left[ -m+1,S\right] ,%
\end{array}
\right.
\end{equation*}
we obtain that the function ${v}_{S}(\mathbf{x},t)$ is the solution of the
following problem%
\begin{equation}
\left. 
\begin{array}{c}
2\left( \partial _{t}^{2}-\Delta -a(\mathbf{x})\right) {v}_{S}=h_{S}(\mathbf{%
\ x},t),~\mathbf{x}\in \mathbb{R}^{n},~t>|\mathbf{x}|, \\ 
{v}_{S}(\mathbf{x},t)|_{t<|\mathbf{x}|}=0,%
\end{array}
\right.  \label{16}
\end{equation}
where 
\begin{equation*}
h_{S}(\mathbf{x},t)=\left( \Delta \alpha _{S}(\mathbf{x})+a(\mathbf{x}
)\alpha _{S}(\mathbf{x})\right) \theta _{S}(t^{2}-|\mathbf{x}|^{2}).
\end{equation*}%
The equality ${v}_{S}(\mathbf{x},t)=0$ for $t<|\mathbf{x}|$ follows from the
fact that the speed of sound is $1$.

Below $B_{1},B_{2}>0$ denote different numbers depending only on the number $%
a_{0}$ in (\ref{2.7}). Let $D_{t}=\{(\mathbf{x},\tau ):\mathbf{x}\in \mathbb{%
\ R}^{n},0<\tau \leq t\}$. The function $h_{S}\in H^{S}(D_{t})$ since $\ell
=2(S+m+1)+S$. Using the method of energy estimates method and condition (\ref%
{2.7}), one can prove the following estimate 
\begin{equation}
\Vert {v}_{S}\Vert _{H^{S+1}(\Omega _{t})}\leq B_{1}e^{B_{2}t}\Vert
h_{S}\Vert _{H^{S}(D_{t})},~\forall t>0,  \label{170}
\end{equation}%
where the domain $\Omega _{t}$ is defined in (\ref{10}), and numbers $C_{1}$
and $C_{2}$ depend only on $a_{0}$. Since 
\begin{equation}
\Vert h_{S}\Vert _{H^{S}(D_{t})}\leq B_{1}t^{\eta _{S}},  \label{171}
\end{equation}%
where $\eta _{S}=S+n+1$, then the following estimate holds 
\begin{equation}
\Vert {v}_{S}\Vert _{H^{S+1}(\Omega _{t})}\leq B_{1}t^{\eta
_{S}}e^{B_{2}t}\Vert h_{S}\Vert _{H^{S}(D_{t})},~\forall t>0.  \label{17}
\end{equation}

Since $2(S+1)>n+2$ then by the Sobolev embedding theorems, ${v}_{S}\in
C^{2}(\Omega _{t})$. Hence, using (\ref{171}) and (\ref{17}), we obtain 
\begin{equation}
\Vert {v}_{S}\Vert _{C^{2}(\Omega _{t})}\leq B_{1}e^{B_{2}t},~\forall t>0.
\label{18}
\end{equation}

Equation (\ref{16}) and inequality (\ref{18}) imply 
\begin{equation}
\Vert \partial _{t}^{2}{v}_{S}\Vert _{C(\Omega _{t})}\leq
B_{1}e^{B_{2}t},~\forall t>0,  \label{172}
\end{equation}%
which implies 
\begin{equation}
\Vert \partial _{t}{v}_{S}\Vert _{C(\Omega _{t})}\leq
B_{1}e^{B_{2}t},~\forall t>0.  \label{173}
\end{equation}%
Hence, we have proved that the function ${v}_{S}(\mathbf{x},t)$ is bounded
in the domain $\Omega _{t}$ together with its derivatives up to the second
order for any fixed $t>0$. In addition, we have proven that this function
grows not faster than $e^{B_{2}t}$ as $t\rightarrow \infty $, together with
its derivatives up to the second order.

Consider now the case $n=2m,m\geq 1$. In this case the representation of the
solution of problem (\ref{2.9}), (\ref{2.10}) has the form (\ref{13}), where
the function ${v}_{S}(\mathbf{x},t)$ is the solution of the Cauchy problem (%
\ref{16}) with $h_{S}(\mathbf{x},t)$ given by 
\begin{equation*}
h_{S}(\mathbf{x},t)=\left( \Delta \alpha _{S}(\mathbf{x})+a(\mathbf{x}
)\alpha _{S}(\mathbf{x})\right) \theta _{S+1/2}(t^{2}-|\mathbf{x}|^{2}).
\end{equation*}%
Since $S=m+1$, then estimates (\ref{170})-(\ref{173}) are also valid for the
function ${v}_{S}(\mathbf{x},t)$ with $\eta _{S}=S+n+3/2$. In addition, the
following estimates are valid:%
\begin{equation}
\left. 
\begin{array}{c}
\Vert v_{S}(\mathbf{x},t)\Vert _{C^{2}(B_{R})}\leq B_{1}e^{B_{2}t},~t>0, \\ 
\Vert \partial _{t}v_{S}(\mathbf{x},t)\Vert _{C^{2}(B_{R})}\leq
B_{1}e^{B_{2}t},~t>0, \\ 
\Vert \partial _{t}^{2}v_{S}(\mathbf{x},t)\Vert _{C^{2}(B_{R})}\leq
B_{1}e^{B_{2}t},~t>0.%
\end{array}
\right.  \label{21}
\end{equation}
$\square $

\section{Appendix B: Proof of Theorem 2.2}

\label{sec:11}\hfill

By (\ref{2.1}) $a\left( \mathbf{x}\right) =0$ in a small neighborhood of the
point $\left\{ \mathbf{x}=0\right\} .$ Hence, uniqueness and existence of
the solution $u\left( \mathbf{x},t\right) $ of problem (\ref{2.15}), (\ref%
{2.16}) satisfying (\ref{2.19}) easily follow from results of Chapter 4 of 
\cite{Lad}. Theorem 2.1 guarantees the existence of transformation (\ref%
{2.13}) for the function $U$ and its derivatives up to the second order.
More precisely, 
\begin{equation*}
u\left( \mathbf{x},t\right) =\mathcal{L}\left( U\right) ,u_{x_{i}}\left( 
\mathbf{x},t\right) =\mathcal{L}\left( U_{x_{i}}\right) ,\text{ }
u_{x_{i}x_{j}}\left( \mathbf{x},t\right) =\mathcal{L}\left(
U_{x_{i}x_{j}}\right) ,\text{ }u_{t}\left( \mathbf{x},t\right) =\mathcal{L}
\left( U_{tt}\right) .
\end{equation*}%
We need to prove asymptotic formulae (\ref{2.20})-(\ref{2.23}).

Let $n=2m+1$. Substituting representation (\ref{11}) in (\ref{2.13}) and
using (\ref{2.12}), we obtain%
\begin{equation*}
u(\mathbf{x},t)=\frac{e^{-\frac{|\mathbf{x}|^{2}}{4t}}}{4\sqrt{\pi t^{3}}}
\int\limits_{0}^{\infty }\left[ \sum\limits_{s=-m}^{S}\alpha _{s}(\mathbf{x}
)\theta _{s}(z)+{v}_{S}(\mathbf{x},\sqrt{z+|\mathbf{x}|^{2}})\right] e^{- 
\frac{z}{4t}}dz.
\end{equation*}%
where $S=m+1,$ and estimates (\ref{18})-(\ref{173}) hold for the function ${v%
}_{S}(\mathbf{x},t)$. Simple calculations lead to the formula 
\begin{equation}
u(\mathbf{x},t)=\frac{e^{-\frac{|\mathbf{x}|^{2}}{4t}}}{4\sqrt{\pi t^{3}}} %
\left[ w_{S}(\mathbf{x},t)+\overline{w}_{S}(\mathbf{x},t)\right] ,\quad t>0,
\label{200}
\end{equation}%
where 
\begin{equation}
w_{S}(\mathbf{x},t)=\sum\limits_{s=-m}^{S}\alpha _{s}(\mathbf{x})(4t)^{s+1},
\label{23}
\end{equation}%
\begin{equation}
\overline{w}_{S}(\mathbf{x},t)=\int\limits_{0}^{\infty }e^{-\frac{z}{4t}}{v}
_{S}\left( \mathbf{x},\sqrt{z+|\mathbf{x}|^{2}}\right) \,dz.  \label{24}
\end{equation}%
Using estimates (\ref{18})-(\ref{173}), we obtain%
\begin{equation}
\left. 
\begin{array}{c}
\overline{w}_{S}(\mathbf{x},t)=O(t),~ \\ 
\partial _{x_{i}}\overline{w}_{S}(\mathbf{x},t)=O(t),~i\in \left[ 1,n\right]
, \\ 
\partial _{x_{i}x_{j}}\overline{w}_{S}(\mathbf{x},t)=O(t),~i,j\in \left[ 1,n %
\right] ,~ \\ 
\partial _{t}\overline{w}_{S}(\mathbf{x},t)=O(t), \\ 
\mathbf{x}\in B_{M},\text{ }t\rightarrow 0^{+},%
\end{array}
\right.  \label{7.1}
\end{equation}%
for any fixed $M>0$. Indeed, to prove the formula in the first line of (\ref%
{7.1}), we use (\ref{18}). Hence,%
\begin{equation*}
\left\vert \overline{w}_{S}(\mathbf{x},t)\right\vert \leq
\int\limits_{0}^{\infty }e^{-\frac{z}{4t}}\left\vert {v}_{S}\left( \mathbf{x}
,\sqrt{z+|\mathbf{x}|^{2}}\right) \right\vert \,dz\leq
B_{1}\int\limits_{0}^{\infty }e^{-\frac{z}{4t}}e^{B_{2}\sqrt{z+|\mathbf{x}
|^{2}}}\,dz\leq
\end{equation*}%
\begin{equation}
\leq 4tB_{1}e^{B_{2}M}\int\limits_{0}^{\infty }e^{-y+2B_{2}\sqrt{ty}
}\,dy=O(t),~\mathbf{x}\in B_{M},~t\rightarrow 0^{+}.  \label{7.2}
\end{equation}%
Estimates in the second and third lines of (\ref{7.1}) of the first and
second derivatives of the function $\overline{w}_{S}(\mathbf{x},t)$ with
respect to $x_{i}$ can be done similarly. However the estimate of the
derivative of the function $\overline{w}_{S}(\mathbf{x},t)$ with respect to $%
t$ in the fourth line of (\ref{7.1}) requires more explanations.

First, we note that 
\begin{equation}
\int\limits_{0}^{\infty }e^{-\frac{s}{4t}}{v}_{S}\left( \mathbf{x},\sqrt{z+| 
\mathbf{x}|^{2}}\right) \,dz=(4t)^{2}\int\limits_{0}^{\infty }e^{-\frac{z}{%
4t }}\partial _{z}^{2}\left[ v_{S}\left( \mathbf{x},\sqrt{z+|\mathbf{x}|^{2}}
\right) \right] \,dz.  \label{7.3}
\end{equation}%
Hence, 
\begin{equation*}
\partial _{t}\int\limits_{0}^{\infty }e^{-\frac{z}{4t}}{v}_{S}\left( \mathbf{%
\ x},\sqrt{z+|\mathbf{x}|^{2}}\right) \,dz=4\int\limits_{0}^{\infty }z\,e^{- 
\frac{s}{4t}}\partial _{z}^{2}\left[ v_{S}\left( \mathbf{x},\sqrt{z+|\mathbf{%
\ x}|^{2}}\right) \right] \,dz.
\end{equation*}%
Hence, using (\ref{172}) and (\ref{7.2}), we obtain%
\begin{equation*}
\left\vert \partial _{t}\int\limits_{0}^{\infty }e^{-\frac{z}{4t}}{v}
_{S}\left( \mathbf{x},\sqrt{z+|\mathbf{x}|^{2}}\right) \,dz\right\vert \leq
4\int\limits_{0}^{\infty }z\,e^{-\frac{s}{4t}}\left\vert \partial _{z}^{2} %
\left[ v_{S}\left( \mathbf{x},\sqrt{z+|\mathbf{x}|^{2}}\right) \right]
\right\vert \,dz\leq
\end{equation*}%
\begin{equation*}
\leq B_{1}\int\limits_{0}^{\infty }\,e^{-\frac{s}{4t}}ze^{B_{2}\sqrt{z+| 
\mathbf{x}|^{2}}}dz=O\left( t\right) ,~\mathbf{x}\in B_{M},~t\rightarrow
0^{+},
\end{equation*}%
which proves the estimate in the fourth line of (\ref{7.1}).

We now work with the function $w_{S}(\mathbf{x},t)$ in (\ref{200}), (\ref{23}%
). First, we note that by (\ref{12}) 
\begin{equation*}
\frac{e^{-\frac{|\mathbf{x}|^{2}}{4t}}}{4\sqrt{\pi t^{3}}}\alpha _{-m}( 
\mathbf{x})(4t)^{1-m}=\frac{1}{(2\sqrt{\pi t})^{n}}\,e^{-\frac{|\mathbf{x}
|^{2}}{4t}}=u_{0}(\mathbf{x},t).
\end{equation*}%
Hence, we can rewrite representation (\ref{200}) in the form 
\begin{equation}
u(\mathbf{x},t)=u_{0}(\mathbf{x},t)\left[ \widehat{w}_{S}(\mathbf{x},t)+2\pi
^{m}(4t)^{m-1}\overline{w}_{S}(\mathbf{x},t)\right] ,\quad t>0,  \label{26}
\end{equation}%
where 
\begin{equation}
\begin{array}{l}
\widehat{w}_{S}(\mathbf{x},t)=1+\sum\limits_{s=-m+1}^{S}\widehat{\alpha }
_{s}(\mathbf{x})(4t)^{s+m}, \\ 
\widehat{\alpha }_{s}(\mathbf{x})=\frac{1}{4}\int\limits_{0}^{1}z^{s+m-1}(
\Delta \alpha _{s-1}(\xi )+a(\xi )\alpha _{s-1}(\xi ))|_{\xi =z\mathbf{x}
}\,dz,~\>s\in \left[ -m+1,S\right] .%
\end{array}
\label{27}
\end{equation}%
Representation (\ref{2.20}) with $\kappa =1$ follows from (\ref{26}) and (%
\ref{27}). Differentiating (\ref{26}) with respect to $t$ and using (\ref%
{7.1}) , we obtain%
\begin{equation*}
u_{t}(\mathbf{x},t)=\partial _{t}u_{0}(\mathbf{x},t)[1+O(t)]+
\end{equation*}

\begin{equation*}
+u_{0}(\mathbf{x},t)\left[ O(t)+2\pi ^{m}\partial _{t}\left(
(4t)^{m-1}\int\limits_{0}^{\infty }e^{-\frac{z}{4t}}{v}_{S}\left( \mathbf{x}%
, \sqrt{z+|\mathbf{x}|^{2}}\right) \,dz\right) \right] ,~\mathbf{x}\in
B_{M}, \text{ }t\rightarrow 0.
\end{equation*}
Using (\ref{7.3}), we obtain 
\begin{eqnarray*}
&&\partial _{t}\left( (4t)^{m-1}\int\limits_{0}^{\infty }e^{-\frac{z}{4t}}{v}
_{S}\left( \mathbf{x},\sqrt{z+|\mathbf{x}|^{2}}\right) \,dz\right) = \\
&&\qquad =4(m+1)(4t)^{m}\int\limits_{0}^{\infty }e^{-\frac{z}{4t}}\partial
_{z}^{2}\left[ v_{S}\left( \mathbf{x},\sqrt{z+|\mathbf{x}|^{2}}\right) %
\right] \,dz= \\
&&\qquad +4(4t)^{m-1}\int\limits_{0}^{\infty }z\,e^{-\frac{z}{4t}}\partial
_{z}^{2}\left[ v_{S}\left( \mathbf{x},\sqrt{z+|\mathbf{x}|^{2}}\right) %
\right] \,dz=O(t^{m}),~\mathbf{x}\in B_{M},~t\rightarrow 0^{+}.
\end{eqnarray*}%
Hence, 
\begin{equation*}
u_{t}(\mathbf{x},t)=\partial _{t}u_{0}(\mathbf{x},t)[1+O(t)],~\mathbf{x}\in
B_{M},~t\rightarrow 0^{+}.
\end{equation*}%
It proves relation (\ref{2.21}) with $\kappa =1$.

Differentiating (\ref{26}) with respect to $x_{i}$ and using again (\ref{7.1}%
), we obtain 
\begin{equation*}
u_{x_{i}}(\mathbf{x},t)=\partial _{x_{i}}u_{0}(\mathbf{x},t)[1+O(t)]+u_{0}(%
\mathbf{x},t)\,O(t),~t\rightarrow 0^{+}.
\end{equation*}%
Also, 
\begin{eqnarray*}
u_{0}(\mathbf{x},t) &=&-\partial _{x_{i}}u_{0}(\mathbf{x},t)\frac{2t}{x_{i}}=
\\
&=&\partial _{x_{i}}u_{0}(\mathbf{x},t)\,O(t),~\mathbf{x}\in G_{\rho
}(M),~t\rightarrow 0^{+}.
\end{eqnarray*}%
Hence, 
\begin{equation*}
u_{x_{i}}(\mathbf{x},t)=\partial _{x_{i}}u_{0}(\mathbf{x},t)[1+O(t)],~%
\mathbf{x}\in G_{\rho }(M),~t\rightarrow 0^{+}.
\end{equation*}%
The latter equality is the same as the one in (\ref{2.22}) with $\kappa =1$.
The proof of equality (\ref{2.23}) with $\kappa =1$ is similar.

Thus, Theorem 2.2 is proven for the case when $n=2m+1,m\geq 1$, $\ell
=5m+3=5[n/2]+3$.

Consider now case $n=2m,$ $m\geq 1$. Then (\ref{2.12}) and (\ref{2.13}) lead
to%
\begin{equation*}
u(\mathbf{x},t)=\frac{e^{-\frac{|\mathbf{x}|^{2}}{4t}}}{4\sqrt{\pi t^{3}}}
\int\limits_{0}^{\infty }e^{-\frac{s}{4t}}\left[ \sum\limits_{s=-m}^{S}
\alpha _{s}(\mathbf{x})\theta _{s+1/2}(z)+{v}_{S}\left( \mathbf{x},\sqrt{z+| 
\mathbf{x}|^{2}}\right) \right] dz.
\end{equation*}
This equality can be represented in the form similar to (\ref{200}): 
\begin{equation}
u(\mathbf{x},t)=\frac{e^{-\frac{|\mathbf{x}|^{2}}{4t}}}{4\sqrt{\pi t^{3}}} %
\left[ 2\sqrt{\pi t}\,w_{S}(\mathbf{x},t)+\overline{w}_{S}(\mathbf{x},t) %
\right] ,\quad t>0,  \label{28}
\end{equation}%
where functions $w_{S}(\mathbf{x},t)$ and $\overline{w}_{S}(\mathbf{x},t)$
are determined by formulae (\ref{23}) and (\ref{24}), respectively.

If $S=m+1$ then estimates (\ref{170})-(\ref{173}) are also valid for ${v}%
_{S}(\mathbf{x},t)$ with $\eta _{S}=S+n+3/2,$ and estimates (\ref{21}) hold
as well.

Note that in this case 
\begin{equation*}
\frac{e^{-\frac{|\mathbf{x}|^{2}}{4t}}}{4\sqrt{\pi t^{3}}}\,\left( 2\sqrt{
\pi t}\right) \alpha _{-m}(4t)^{1-m}=\frac{e^{-\frac{|\mathbf{x}|^{2}}{4t}}}{
(4{\pi t})^{n/2}}=u_{0}(\mathbf{x},t).
\end{equation*}%
Hence, representation (\ref{28}) we can be written as: 
\begin{equation}
u(\mathbf{x},t)=u_{0}(\mathbf{x},t)\left[ \widehat{w}_{K}(\mathbf{x},t)+\pi
^{m-1/2}(4t)^{m-3/2}\bar{w}_{K}(\mathbf{x},t)\right] ,\quad t>0,  \label{29}
\end{equation}%
where $\widehat{w}_{K}(\mathbf{x},t)$ is given by the first line of (\ref{27}%
). Note that representation (\ref{29}) differs from representation (\ref{26}%
) by the term $(4t)^{m-3/2}$, where the power $m-3/2$ is less by $1/2,$ as
compares with ${m-1}$ in (\ref{26}). It is important if $m=1,$ then $n=2.$
only. If $m>1$, then the validity of Theorem 2.2 with $\kappa =1$ follows
from representation (\ref{29}) and inequalities (\ref{170})-(\ref{173}), (%
\ref{21}) just as it was in the above case when $n$ is odd. However, if $%
n=2, $ then $m=1$, and we need to put $\kappa =1/2$ in estimates (\ref{2.20}%
)-(\ref{2.23}) .

Thus, Theorem 2.2 is proven in both cases: odd and even $n$. $\square $

\end{document}